\documentclass[10pt,a4paper]{article}

\usepackage{amsmath,amssymb,amsfonts,amsthm}
\usepackage{stmaryrd}
\usepackage{xfrac}
\usepackage{multirow}
\usepackage{colortbl}
\usepackage[table]{xcolor}
\usepackage{enumitem}
\usepackage{subcaption}
\usepackage{lscape}
\usepackage{fullpage}
\usepackage{color}

\usepackage{authblk}
\usepackage[numbers,sort&compress]{natbib}
\renewcommand{\cite}[1]{\citep{#1}}

\usepackage{ifthen}
\usepackage{graphicx, caption}
\usepackage{cancel}
\usepackage{bm}
\graphicspath{ 
  {images/}{../images/} 
}
\usepackage{lipsum}
\usepackage{varwidth}

\usepackage[ruled,boxed]{algorithm2e}
\SetKwFor{ParFor}{parfor}{do}{end}
\usepackage{multicol}
\usepackage{algpseudocode}
\usepackage{epsfig}

\usepackage{float}

\usepackage{tikz}
\usetikzlibrary{arrows}
\tikzstyle{vertex} = [circle, draw, scale=1.00]
\tikzstyle{pertex} = [circle, draw, scale=0.50]
\tikzstyle{square} = [        draw, scale=1.20]
\tikzstyle{pquare} = [        draw, scale=0.60]

\definecolor{gray          }{gray}{0.50}
\definecolor{lightgray     }{gray}{0.75}
\definecolor{lightlightgray}{gray}{0.90}

\usepackage[pdftex=true,colorlinks=true,citecolor=red,linkcolor=blue]{hyperref}

  \newcommand{\im  }   {\imath               }

\newcommand{\bmat}[1]{ \boldsymbol{\mathcal{ #1}} }
\newcommand{\bvec}[1]{ \mathbf{#1} }

\theoremstyle{definition}

\title{Multidirectionnal sweeping preconditioners with non-overlapping checkerboard domain decomposition for Helmholtz problems}
\author[1,3]{R. Dai}
\author[2]{A. Modave}
\author[1]{J.-F. Remacle}
\author[3]{C. Geuzaine}

\affil[1]{\footnotesize IMMC, Universit\'e catholique de Louvain, 1348 Louvain-la-Neuve, Belgium,
\authorcr \href{ruiyang.dai@uclouvain.be}{ruiyang.dai@uclouvain.be},
	  \href{jean-francois.remacle@uclouvain.be}{jean-francois.remacle@uclouvain.be} }
\affil[2]{\footnotesize POEMS, CNRS, Inria, ENSTA Paris, Institut Polytechnique de Paris, 91120 Palaiseau, France,
\authorcr \href{axel.modave@ensta-paris.fr}{axel.modave@ensta-paris.fr} }
\affil[3]{\footnotesize Universit\'e de Li\`ege, Institut Montefiore B28, 4000 Li\`ege, Belgium,
\authorcr \href{ruiyang.dai@uclouvain.be}{ruiyang.dai@uclouvain.be},
          \href{cgeuzaine@uliege.be}{cgeuzaine@uliege.be} }

\date{}

\begin{document}

\maketitle

\begin{abstract}
This paper explores a family of generalized sweeping preconditionners for Helmholtz problems with non-overlapping checkerboard partition of the computational domain.
The domain decomposition procedure relies on high-order transmission conditions and cross-point treatments, which cannot scale without an efficient preconditioning technique when the number of subdomains increases.
With the proposed approach, existing sweeping preconditioners, such as the symmetric Gauss-Seidel and parallel double sweep preconditioners, can be applied to checkerboard partitions with different sweeping directions (\textit{e.g.}~horizontal and diagonal).
Several directions can be combined thanks to the flexible version of GMRES, allowing for the rapid transfer of information in the different zones of the computational domain, then accelerating the convergence of the final iterative solution procedure.
Several two-dimensional finite element results are proposed to study and to compare the sweeping preconditioners, and to illustrate the performance on cases of increasing complexity.
\end{abstract}

\setlength{\parskip}{4pt}

\section{Introduction}
Time-harmonic wave simulations are of interest in many scientific
and engineering disciplines. For example, in radar or sonar imaging and wireless
communications, the wavelength of the signal is usually several orders of
magnitude smaller than the size of the domains of interest.  Similarly, in
seismic imaging, wave fields in complex geological media show a wide range of
space-varying wavenumbers, caused by large variations in the velocity profile.
Solving such time-harmonic problems numerically using finite
element-type methods is notoriously difficult because it leads to (extremely)
large indefinite linear systems~\cite{ernst2012difficult}, especially in the high-frequency regime.  One the one hand,
sparse direct solvers exhibit poor scalability w.r.t. memory and computational
time for such linear systems, in particular for three-dimensional problems. On
the other hand, most iterative methods that have proved successful for elliptic
problems become inefficient when applied to problems with highly oscillatory
solutions, and no robust and scalable preconditioner currently
exists~\cite{Lahaye2017}.

Parallel iterative solvers and parallel preconditionners, called generically \textit{``domain decomposition methods''} (DDMs), are currently intensively studied for time-harmonic problems.
These methods rely on the parallel solution of subproblems of smaller sizes, amenable to sparse direct solvers.
In a finite element context, there are two types of DDMs: overlapping DDMs in which the
meshes of two adjacent subdomains overlap by at least one element; and
non-overlapping DDMs where adjacent subdomains only communicate through an
artificial (lower dimensional) interface.
The latter include \textit{e.g.}~non-overlapping
Schwarz methods \cite{Despres1991}, FETI algorithms \cite{de1998non,
 farhat2005feti} and the method of polarized traces \cite{zepeda2016method}.

In this work, we focus on non-overlapping domain decomposition solvers with optimized
transmission conditions, which are well suited for time-harmonic wave problems
\cite{Despres1991,gander2002optimized,boubendir2012quasi}. After
non-overlapping Schwarz methods were introduced by Lions \cite{Lions1990} for the
Laplace equation and proven to converge for the Helmholtz equation by
Despr\'{e}s \cite{Despres1991}, considerable efforts have been made to develop
efficient transmission conditions to improve the rate of convergence for
DDMs. The optimal convergence is obtained by using as transmission condition on
each interface the Dirichlet-to-Neumann (DtN) map related to the complementary
of the subdomain of interest. This DtN map is a nonlocal operator and is thus in
practice very expensive to compute. Optimized Schwarz methods were introduced in
\cite{gander2002optimized,gander2006optimized}, where the nonlocal DtN is approximated by first or
second order polynomial approximations, with coefficients obtained by
optimization on simple geometries. Later, quasi-optimal optimized Schwarz
methods were proposed in \cite{boubendir2012quasi} based on rational
approximations related to those used in high-order absorbing boundary conditions
(HABC). Conditions based on second-order operators \cite{nicolopoulos2019formulations}, perfectly matched layers (PML) \cite{stolk2013rapidly, vion2014double, astaneh2016two} and non-local approaches \cite{lecouvez2014quasi,collino2020exponentially} have also been investigated.

Even with optimal transmission conditions however, as is expected for a
one-level method, the number of iterations of the DDM will grow as the number of
subdomains increases. A solution for certain classes of problems is to add a
component to the algorithm that is known in the DDM community as a \textit{``coarse
grid''} \cite{farhat1998unified}, in effect a second-level to enable longer-range
information exchange than the local sharing (from one subdomain to its
neighbors) of the one-level DDM.
Nevertheless, the design of robust coarse grids is very challenging for wave-type problems because of the highly oscillatory behavior of the solution, and several approaches are currently investigated in the community (see \textit{e.g.}~\cite{conen2014coarse,bootland2020comparison} and references herein).
As an alternative approach, sweeping preconditioners have been proposed and studied for convection-diffusion problems in the 90's
\cite{nataf1993use,nataf1997convergence}.
They have recently garnered a lot of interest for
high-frequency Helmholtz problems~\cite{engquist2011sweeping, stolk2013rapidly, vion2014double,
zepeda2016method, gander2019class, bouziani2021overlapping, taus2020sweeps}, promising a number
of DDM iterations that is quasi independent of the number of subdomains.
In particular, they are very effective for waveguide and open cavity configurations, which exhibit a unique natural direction in which information is transferred.
However, they
have two major drawbacks: they rely on intrinsically sequential operations (they
are related to a LU-type factorization of the underlying iteration operator)
and they are naturally only suited for layered-type domain decompositions (where
the layered structure allows to explicit the LU factorization as a double sweep
across the subdomains).

In this work, we explore a family of generalized sweeping preconditionners where
sweeps can be done in several directions for non-overlapping domain decomposition solvers with
\textit{``checkerboard''} domain partition. This contribution relies on the
availability of transmission conditions able to deal with the cross-points
(\textit{i.e.}~points where more than two subdomains meet) arising in such configurations.
We consider a domain decomposition solver with high-order Pad\'e-type transmission
conditions \cite{boubendir2012quasi} and a cross-point treatment proposed in \cite{modave2020non}.
Sweeping preconditioners are derived in a systematic manner, based on the explicit representation of the iteration matrix in the case of checkerboard decompositions.
The sweeps can be performed in Cartesian and diagonal directions, and several sweeping directions can be combined by using the flexible version of GMRES \cite{Saad1993,saad2003iterative}, which allows to change the preconditioner at each iteration.
For applicative cases, the resulting preconditioners provide an effective way to rapidly transfer information in the different zones of the computational domain, then accelerating the convergence of iterative solution procedure with GMRES.
Our approach is related to the recent work on L-sweeps preconditioners \cite{taus2020sweeps} and diagonal sweeping technique \cite{leng2020diagonal}, where sweeping strategies are proposed in the context of the method of polarized traces and the source transfer method, respectively.
Here, the preconditioners are proposed for non-overlapping domain decomposition solvers with high-order transmission conditions, and several directions can be combined thanks to the use of flexible GMRES.



The paper is organized as follows. In Section \ref{sec:ddmHelmholtz}, we present
the non-overlapping domain decomposition algorithm with high-order transmission
condition and cross-points treatment for checkerboard partitions. The
matrix representation of the iteration operator is derived and studied in
Section \ref{Matrix Form of the Iteration Operator}. It is used in
Section \ref{Parallel Sweeping Preconditioners} to derive sweeping
preconditioners.  In Section \ref{Benchmarks}, we study and compare the 
preconditioned domain decomposition algorithms with two-dimensional finite element benchmarks.
The efficiency of these methods is demonstrated on numerical models and configurations of increasing complexity.

\section{Domain decomposition algorithm for the Helmholtz equation}
\label{sec:ddmHelmholtz}

To describe our approach, we consider the two-dimensional scattering problem of
an incident acoustic plane wave by a sound-soft obstacle of boundary
$\Gamma^{\text{sca}}$.  The numerical simulations are performed in a rectangular
computational domain $\Omega$ of boundary
$\partial\Omega=\Gamma^{\text{sca}}\bigcup\Gamma^{\infty}$, with $\Gamma^{\infty}$
the external (artificial) boundary (see
Figure~\ref{fig:illustration_partition_and_numbering}, left).  We seek the
scattered field $u(\mathbf{x})$ that verifies
\begin{equation}
  \left\{
  \begin{aligned}
    - \Delta u - \kappa^2u &= 0, && \text{in } \Omega, \\
    \partial_{\bm{n}} u + \mathcal{B} u &=  0, && \text{on } \Gamma^{\infty}, \\
    u &= -u^{\text{inc}}, && \text{on } \Gamma^{\text{sca}},
  \end{aligned}
  \right.
  \label{eqn:origPbm}
\end{equation}
where $\kappa$ is the wavenumber, $u^{\text{inc}}$ is the incident
wave, $\partial_{\bm{n}}$ is the exterior normal derivative and $\mathcal{B}$ is
an impedance operator to be defined.  We take the convention that the
time-dependence of the fields is $e^{-\im\omega t}$, where $\omega$ is the
angular frequency and $t$ is the time.  The impedance operator
corresponding to a Sommerfeld absorbing boundary condition (ABC) on
$\Gamma^{\infty}$ is $\mathcal{B} = -\im\kappa$.


\subsection{Domain decomposition algorithm on a checkerboard partition}
\label{sec:ddStandard}

\begin{figure}[!tb]
  \centering \small
  \captionsetup{font=small, labelfont=bf}
  \begin{tikzpicture}[scale=1.00]

    \def\length{3}
    \def\xI{0}
    \def\yI{0}
    \def\xJ{5.0}
    \def\yJ{0}
    \def\xK{10.0}
    \def\yK{0}

    \coordinate (n1) at (\xI        , \yI        );
    \coordinate (n2) at (\xI+\length, \yI        );
    \coordinate (n3) at (\xI+\length, \yI+\length);
    \coordinate (n4) at (\xI        , \yI+\length);

    \coordinate (n5) at (\xJ        , \yJ        );
    \coordinate (n6) at (\xJ+\length, \yJ        );
    \coordinate (n7) at (\xJ+\length, \yJ+\length);
    \coordinate (n8) at (\xJ        , \yJ+\length);

    \coordinate (n9)  at (\xK        , \yK        );
    \coordinate (n10) at (\xK+\length, \yK        );
    \coordinate (n11) at (\xK+\length, \yK+\length);
    \coordinate (n12) at (\xK        , \yK+\length);

    \draw[very thick, black] (n1)--(n2);
    \draw[very thick, black] (n2)--(n3);
    \draw[very thick, black] (n3)--(n4);
    \draw[very thick, black] (n4)--(n1);
    \draw[very thick, black] (\xI+0.5, \yI+0.5) circle (0.3);
    \node[] at (\xI+\length/2, \yI+\length/2) {\LARGE$\Omega$};
    \node[] at (\xI+1.25, \yI+0.45) {\Large$\Gamma^{\mathrm{sca}}$};
    \node[] at (\xI        -0.50, \yI+\length/2) {\Large$\Gamma^{\infty}$};

    \draw[very thick, black] (n5)--(n6);
    \draw[very thick, black] (n6)--(n7);
    \draw[very thick, black] (n7)--(n8);
    \draw[very thick, black] (n8)--(n5);
    \draw[very thick, black] (\xJ+0.5, \yJ+0.5) circle (0.3);
    \draw[very thick, black] (\xJ, \yJ+\length/3*1)--(\xJ+\length, \yJ+\length/3*1);
    \draw[very thick, black] (\xJ, \yJ+\length/3*2)--(\xJ+\length, \yJ+\length/3*2);
    \draw[very thick, black] (\xJ+\length/3*1, \yJ)--(\xJ+\length/3*1, \yJ+\length);
    \draw[very thick, black] (\xJ+\length/3*2, \yJ)--(\xJ+\length/3*2, \yJ+\length);
    \node[] at (\xJ+\length/6*1-1, \yJ+\length/6*1) {\large$\Omega_1$};
    \draw[->, thick, black] (\xJ+\length/22, \yJ+\length/12*3) to
	 [out=165,in=15, looseness=1](\xJ-\length/12, \yJ+\length/12*3);
    \node[] at (\xJ+\length/6*1, \yJ+\length/6*3) {\large$\Omega_2$};
    \node[] at (\xJ+\length/6*1, \yJ+\length/6*5) {\large$\Omega_3$};
    \node[] at (\xJ+\length/6*3, \yJ+\length/6*1) {\large$\Omega_4$};
    \node[] at (\xJ+\length/6*3, \yJ+\length/6*3) {\large$\Omega_5$};
    \node[] at (\xJ+\length/6*3, \yJ+\length/6*5) {\large$\Omega_6$};
    \node[] at (\xJ+\length/6*5, \yJ+\length/6*1) {\large$\Omega_7$};
    \node[] at (\xJ+\length/6*5, \yJ+\length/6*3) {\large$\Omega_8$};
    \node[] at (\xJ+\length/6*5, \yJ+\length/6*5) {\large$\Omega_9$};

    \draw[very thick, black] (n9)--(n10);
    \draw[very thick, black] (n10)--(n11);
    \draw[very thick, black] (n11)--(n12);
    \draw[very thick, black] (n12)--(n9);

    \node[] at (\xK+\length/2   , \yK+\length/2) {\LARGE$\Omega_9$};
    \node[] at (\xK        -0.60, \yK+\length/2) {\Large$\Gamma_{9, 6}$};
    \node[] at (\xK+\length+0.75, \yK+\length/2) {\Large$\Gamma_{9,-3}$};
    \node[] at (\xK+\length/2, \yK        -0.70) {\Large$\Gamma_{9, 8}$};
    \node[] at (\xK+\length/2, \yK+\length+0.70) {\Large$\Gamma_{9,-4}$};

  \end{tikzpicture}
  \caption{\small Configuration of the problem \textit{(left)}, illustration of the checkerboard partition \textit{(middle)} and notation for the edges of the subdomain $\Omega_I$ \textit{(right)}.}
  \label{fig:illustration_partition_and_numbering}
\end{figure}
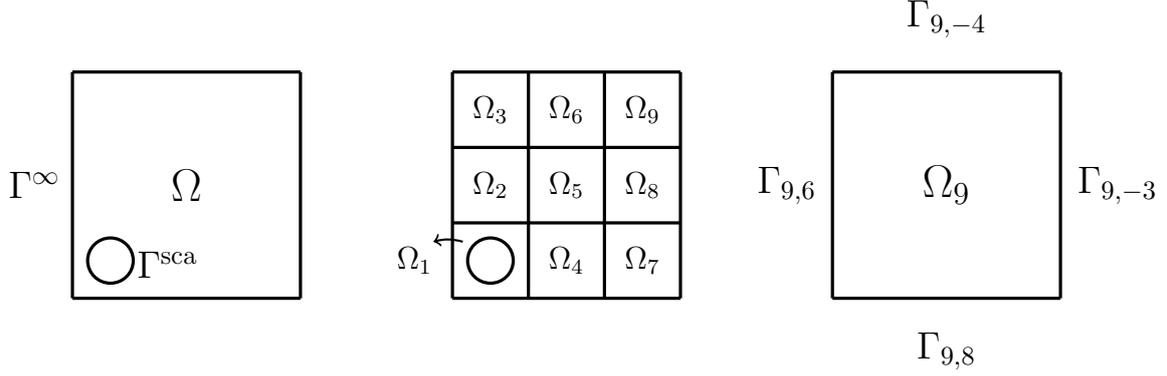

We consider a checkerboard partition of the domain $\Omega$, that consists in a
lattice of rectangular non-overlapping subdomains $\Omega_I$ $(I= 1 \dots N_{\text{dom}})$ with $N_r$ rows and $N_c$
columns (then, $N_{\text{dom}} = N_r \times N_c$). 
For each rectangular subdomain $\Omega_I$, there are four edges, which are on the left, on the bottom, on the right, and on the top of the subdomain, respectively (see Figure~\ref{fig:illustration_partition_and_numbering}, right).
and we define the set
\[
  D_I^{\infty} := \left\{ J\in\{ -1, -2, -3, -4\}
  \text{ such that $\Gamma_{IJ}$} \neq \varnothing \right\},
\]
where the superscripts $-1$, $-2$, $-3$, and $-4$ correspond to the non-empty edges belonging to $\partial \Omega$ on the left, on the bottom, on the right, and on the top of the subdomain, respectively.
The union of the edges of $\Omega_I$ then reads
$$
  \textstyle \left( \bigcup_{J \in D_I} \Gamma_{IJ} \right) \bigcup \left( \bigcup_{J \in D_I^{\infty}} \Gamma_{IJ} \right),
$$
where the set $D_I^{\infty} := \left\{ J\in\{1,\dots,N_{\text{dom}}\} \text{ such that $J \leq I$ and $\Gamma_{IJ}$} \neq \varnothing \right\}.$
To simplify the presentation, we assume that the obstacle is included in only
one subdomain, the boundary of which is the union of the four edges and the boundary of the obstacle (see Figure~\ref{fig:illustration_partition_and_numbering}, middle).

Each edge of one subdomain $\Omega_I$ is either a \textit{boundary edge} if it
belongs to the boundary of the global domain
$(J \in D_I^{\infty})$ or an \textit{interface edge} if there is
a neighboring subdomain on the other side of the edge
$(J \in D_I)$.  In this checkerboard partition,
there are corners where at least two edges meet.  Each corner of a subdomain is
an \textit{interior cross-point} (point that belongs to four subdomains), a
\textit{boundary cross-point} (point that belongs to two subdomains and to the
exterior border $\partial\Omega$) or a corner of the domain $\Omega$.

With these definitions, the non-overlapping domain decomposition algorithm can
be set up as follows.  For each subdomain $\Omega_I$, we seek the solution
$u_I(\mathbf{x})$ of the subproblem
\begin{equation}
  \left\{
  \begin{aligned}
    - \Delta u_I - \kappa^2 \, u_I &= 0,
      && \text{in } \Omega_I, \\
	  \partial_{\bm{n}_{IJ}} u_I + \mathcal{B}_{IJ} \, u_I &= 0,
      && \text{on each } \Gamma_{IJ}, \forall J \in D_I^{\infty},\\
	  \partial_{\bm{n}_{IJ}} u_I + \mathcal{B}_{IJ} \, u_I &= g_{IJ},
      && \text{on each } \Gamma_{IJ}, \forall J \in D_I,\\
    u_I &= -u_{\text{inc}}, && \text{on } \partial\Omega_I\cap\Gamma_{\text{sca}},
  \end{aligned}
  \right.
  \label{eqn:subPbm}
\end{equation}
where $\bm{n}_{IJ}$ is the outward unit normal to the edge $\Gamma_{IJ}$,
$\mathcal{B}_{IJ}$ is an impedance operator and $g_{IJ}$ is a transmission
variable defined on $\Gamma_{IJ}$.  The second and third equations in
\eqref{eqn:subPbm} are boundary and transmission conditions, respectively.  

For
a given boundary edge $\Gamma_{IJ}\subset\partial\Omega$, we must have
$\mathcal{B}_{IJ}=\mathcal{B}$ to ensure the equivalence between all the
subproblems and the original problem.  If
$\Gamma_{IJ} \not\subset\partial\Omega$, there is some flexibility in the
choice of $\mathcal{B}_{IJ}$.  The transmission variable is defined as
\begin{equation}
  g_{IJ} := \partial_{\bm{n}_{IJ}} \, u_{J} + \mathcal{B}_{IJ} \, u_{J},
\end{equation}
where $u_J$ is the solution of the neighboring subdomain $\Omega_J$.  The
transmission conditions defined on both sides of the interface enforce the continuity of
the solution across the interface.  Assuming that the impedance operators used
on both sides of the shared interface edge $\Gamma_{IJ}=\Gamma_{Ji}=\partial\Omega_I\cap\partial\Omega_J$ are the
same (\textit{i.e.}~$\mathcal{B}_{IJ} = \mathcal{B}_{JI}$), the transmission
variables defined on this edge verify
\begin{equation}
  g_{IJ} = - g_{JI} + 2 \mathcal{B}_{JI} \, u_J,
  \label{eqn:interface}
\end{equation}
where $g_{JI}$ is the transmission variable defined on the edge $\Gamma_{JI}$
of $\Omega_J$.

The non-overlapping optimized Schwarz domain decomposition algorithm consists in
solving subproblems associated to all the subdomains (equation
\eqref{eqn:subPbm}) concurrently and updating the transmissions variables using
\eqref{eqn:interface} in an interative process.  At each iteration $n+1$, the
update formula of a transmission variable living on an interface edge
$\Gamma_{IJ}$ of a subdomain $\Omega_I$ reads
\begin{equation}
  g_{IJ}^{(n+1)} = - g_{JI}^{(n)} + 2 \mathcal{B}_{JI} \, u_J^{(n)},
  \label{eq:updateTransVar}
\end{equation}
where $u_J^{(n)}$ is the solution of the neighboring subdomain $\Omega_J$ at the
iteration $n$.
The update of all the transmission variables can be recast as
one application of the iteration operator $\boldsymbol{\mathcal{A}}$ defined by
\begin{equation}
  \mathbf{g}^{(n+1)} = \boldsymbol{\mathcal{A}} \mathbf{g}^{(n)} + \mathbf{b},
\end{equation}
where $\mathbf{g}^{(n)}$ is the set of all transmission variables defined on the
interface edges and $\mathbf{b}$ is given by the source term.  It is well known
that this algorithm can be seen as a Jacobi scheme applied to the linear system
\begin{equation}
  (\boldsymbol{\mathcal{I}} - \boldsymbol{\mathcal{A}}) \mathbf{g} = \mathbf{b},
\end{equation}
where $\boldsymbol{\mathcal{I}}$ is the identity operator.  In order to
accelerate the convergence of the procedure, this system can be solved with
Krylov subspace iterative methods, such as GMRES.


\subsection{Transmission operators}
\label{subsec:transmissionOperators}

The convergence rate of the non-overlapping DDMs strongly depends on the
impedance operator used in the transmission conditions.  The optimal
transmission operator corresponds to the non-local Dirichlet-to-Neumann (DtN)
map related to the complementary of each subdomain.  Since the cost of computing
the exact DtN is prohibitive, strategies based on approximate DtN operators
started to be investigated in the late 80's and early 90's (see
\textit{e.g.}~\cite{hagstrom1988numerical, nataf1994optimal}).  For Helmholtz
problems, Despr\'es \cite{despres1991domain, benamou1997domain} used a
Robin-type operator, which is a coarse approximation of the exact DtN operator.
Improved methods with optimized second-order transmission operators have next
been introduced in \cite{piacentini1998improved, gander2002optimized}.  More
recently, domain decomposition approaches with improved convergence rates have
been proposed by using transmission conditions based on high-order absorbing
boundary conditions (HABCs) \cite{boubendir2012quasi, boubendir2018non,
  kim2015optimized, marsic2020convergence}, perfectly matched layers (PMLs)
\cite{schadle2007additive, stolk2013rapidly, vion2014double, astaneh2016two} and
nonlocal operators \cite{lecouvez2014quasi, stupfel2010improved,
  collino2020exponentially}.  As for ABCs, transmission boundary conditions
related to HABCs and PMLs represent a good compromise between the basic
impedance conditions (which lead to suboptimal convergence) and nonlocal
approaches (which are expensive to compute).

In this work, we use the DtN operator associated to the Pad\'e-type HABC as the
impedance operator in the transmission conditions, following
\cite{boubendir2012quasi}.  For an edge $\Gamma_{IJ}$, the operator can be written
as 
\begin{equation}
  \mathcal{B}_{IJ} = -\im\kappa\alpha \left[1 + \frac{2}{M} \sum_{i=1}^{N} c_i \left(1 - \alpha^2 (c_i+1)\left[(\alpha^2 c_i + 1) + \partial_{\bm{\tau}\bm{\tau}}/\kappa^2\right]^{-1}\right)\right], \quad\text{on } \Gamma_{IJ},
  \label{eqn:padeOp}
\end{equation}
where $\partial_{\bm{\tau}}$ is the tangential derivative and we have $\alpha=e^{\im\phi/2}$,
$c_i=\tan^2(i\pi/M)$ and $M=2N+1$.  This Pad\'e-type impedance operator is
obtained by approximating the exact non-local DtN map associated to the exterior
half-plane problem (see
\textit{e.g.}~\cite{engquist1977absorbing,milinazzo1997rational}).  The symbol
of the non-local operator exhibits a square-root which is replaced with the
$(2N+1)^\text{th}$-order Pad\'e approximation after a $\phi$-rotation of the
branch-cut.  The performance of the obtained operator depends on the number of
terms $N$ and the angle of rotation $\phi$.  The particular parameters $N=0$ and
$\phi=0$ leads to the basic ABC operator $\mathcal{B}_{IJ}=-\im\kappa$.  See
\textit{e.g.} \cite{kechroud2005numerical,modave2020corner} for further details.

For the effective implementation of the transmission condition, the application
of the Pad\'e-type impedance operator on a field is written in such a way that
it involves only differential operators.  Following an approach first used by
Lindman \cite{lindman1975free} for ABCs, we introduce $N$ auxiliary fields
governed by auxiliary equations on
interface edge $\Gamma_{IJ}$.  The application of the Pad\'e-type impedance
operator is then written as
\begin{equation}
  \mathcal{B}_{IJ} u_I
  = B\Big(u_I,\{\varphi_{IJ,i}\}_{i=1\dots N}\Big)
  := -\im\kappa\alpha\left[u_I + \frac{2}{M} \sum_{i=1}^{N} c_i \left(u_I + \varphi_{IJ,i}\right)\right],
  \quad\text{on } \Gamma_{IJ},
  \label{eqn:padeImpOp}
\end{equation}
with the auxiliary fields $\{\varphi_{IJ,i}\}_{i=1\dots N}$ defined only on the edge and governed by the auxiliary equations
\begin{align}
  - \partial_{\bm{\tau}\bm{\tau}} \varphi_{IJ,i}
        - \kappa^2 \big((\alpha^2 c_{i}+1)\varphi_{IJ,i} + \alpha^2 (c_{i}+1)u_I\Big) &= 0,
  \quad \text{on } \Gamma_{IJ},
  \label{eqn:padeAuxEqn}
\end{align}
for $i=1 \dots N$.  The linear multivariate function $B$ is introduced to
simplify the expressions in the remainder of the paper.  When this operator is
used in a boundary condition for polygonal domains, a special treatment is
required to preserve the accuracy of the solution at the corners.  In the case
of right-angle corners, an approach based on compatibility relations reveals to
be very efficient \cite{modave2020corner}.


\subsection{Cross-point treatment}
\label{subsec:crossPointTreatmentHABC}

When the Pad\'e-type impedance operator \eqref{eqn:padeImpOp} is used in the boundary conditions and the interface conditions of the subproblem \eqref{eqn:subPbm}, a special treatment is required at the corners of the subdomain $\Omega_I$.
Indeed, the auxiliary fields governed by equation \eqref{eqn:padeAuxEqn} on the edges require boundary conditions at the extremities of the edges, which corners of the subdomain.
In this work, we use a cross-point treatment based on compatibility relations developped for the right-angle case, first proposed in \cite{modave2020non}.

\begin{figure}[!tb]
  \centering\small
  \begin{tikzpicture}[scale=1.00]

    \def\length{3}
    \def\xI{0}
    \def\yI{0}
    \def\xJ{6}
    \def\yJ{0}
    \def\xK{6+\length}
    \def\yK{0}

    \coordinate (n5) at (\xJ        , \yJ        );
    \coordinate (n6) at (\xJ+\length, \yJ        );
    \coordinate (n7) at (\xJ+\length, \yJ+\length);
    \coordinate (n8) at (\xJ        , \yJ+\length);
    \filldraw[black] (n5) circle (2pt);
    \filldraw[black] (n6) circle (2pt);
    \filldraw[black] (n7) circle (2pt);
    \filldraw[black] (n8) circle (2pt);

    \coordinate (n9)  at (\xK        , \yK        );
    \coordinate (n10) at (\xK+\length, \yK        );
    \coordinate (n11) at (\xK+\length, \yK+\length);
    \coordinate (n12) at (\xK        , \yK+\length);
    \filldraw[black] (n9)  circle (2pt);
    \filldraw[black] (n10) circle (2pt);
    \filldraw[black] (n11) circle (2pt);
    \filldraw[black] (n12) circle (2pt);
    \node[] at (\xJ+\length/2, \yJ+\length/2) {\large\color{black}$\Omega_I$};
    \node[] at (\xK+\length/2, \yK+\length/2) {\large\color{black}$\Omega_J$};

    \draw[dashed, black] (-0.75+\xJ, \yJ+\length/3*0)--(+0.75+\xJ+\length, \yJ+\length/3*0);
    \draw[dashed, black] (-0.75+\xJ, \yJ+\length/3*3)--(+0.75+\xJ+\length, \yJ+\length/3*3);
    \draw[dashed, black] (\xJ+\length/3*0, -0.75+\yJ)--(\xJ+\length/3*0, +0.75+\yJ+\length);
    \draw[->, dashed, black] (\xJ+\length/3*3, -0.75+\yJ)--(\xJ+\length/3*3, +0.75+\yJ+\length);
    \draw[very thick, black] (n5)--(n6);
    \draw[very thick, black] (n6)--(n7);
    \draw[very thick, black] (n7)--(n8);
    \draw[very thick, black] (n8)--(n5);

    \draw[->, dashed, black] (-0.75+\xK, \yK+\length/3*0)--(+0.75+\xK+\length, \yK+\length/3*0);
    \draw[    dashed, black] (-0.75+\xK, \yK+\length/3*3)--(+0.75+\xK+\length, \yK+\length/3*3);
    \draw[dashed, black]     (\xK+\length/3*3, -0.75+\yK)--(\xK+\length/3*3, +0.75+\yK+\length);
    \draw[very thick, black] (n9)--(n10);
    \draw[very thick, black] (n10)--(n11);
    \draw[very thick, black] (n11)--(n12);

    \node[above] at (+0.75+\xK+\length, \yK+\length/3*0) {$x$};
    \node[above] at (\xJ+\length/3*3, +0.75+\yJ+\length) {$y$};
    \node[below, xshift=0.3cm] at (n6) {$O$};

    \filldraw[black]        (n12) circle (2pt);
    \filldraw[red!80!black] (n6)  circle (3pt);

    \draw[->, thick, blue!80!black] (\xJ+\length, \yJ+\length/2) to (\xJ+\length+0.5, \yJ+\length/2);
    \node[above] at (\xJ+\length+0.5, \yJ+\length/2) {\color{blue!80!black}$x$};

    \node[] at(\xJ+\length, \yJ+\length/6)   {\large\color{red!80!black}$P_I^{x y}=P_J^{x y}$ };
    \node[] at(\xJ+\length, \yJ+\length/6*5) {\large\color{black}$\Gamma_{I}^{x}=\Gamma_{J}^{x}$ };

    \node[left] at (\xJ+\length/2+0.5, -0.4) {\large\color{black}$\Gamma_{I}^{y}$ };
    \node[left] at (\xK+\length/2+0.5, -0.4) {\large\color{black}$\Gamma_{J}^{y}$ };
    \node[left] at (\xJ+\length/2+0.5, \length+0.4) {\large\color{black}$\Gamma_{I}^{y'}$ };
    \node[left] at (\xK+\length/2+0.5, \length+0.4) {\large\color{black}$\Gamma_{J}^{y'}$ };

  \end{tikzpicture}
  \caption{Configuration with two subdomains.}
  \label{fig:shared_egdes_points}
\end{figure}
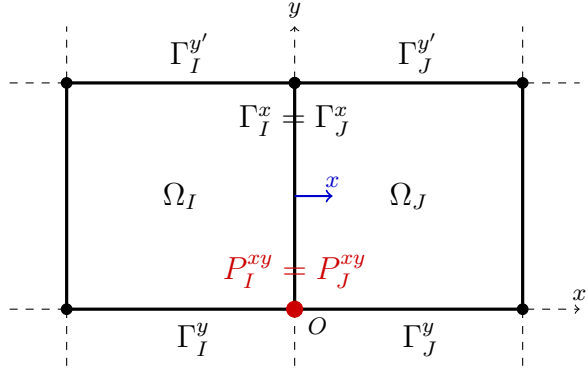

To present the cross-point treatment, we consider the neighboring subdomains $\Omega_I$ and $\Omega_J$ represented on Figure \ref{fig:shared_egdes_points}.
For the sake of shortness, we describe the methodology in the case where transmission conditions are prescribed on all the edges of both subdomains (\textit{i.e.}~they all are interface edges) and the Pad\'e-type impedance operator is used with the same parameters on all the edges.

The subdomains share the interface edge $\Gamma_I^{x} = \Gamma_J^{x}$ ($\Gamma_I^{x} := \Gamma_{IJ}$  and $\Gamma_J^{x} := \Gamma_{JI}$) and the interior cross-points $P_I^{x y} = P_J^{x y}$ and $P_I^{x y'} = P_J^{x y'}$.
On the shared interface edge $\Gamma_I^{x} = \Gamma_J^{x}$, we have the transmission conditions
\begin{align}
  \partial_{x} \, u_I + B\Big(u_I,\{\varphi_{I,i}^x\}_{i=1\dots N}\Big) &= g_{I}^x,
  \quad \text{on } \Gamma_I^{x}, \label{eqn:transCondLeft} \\
 -\partial_{x} \, u_J + B\Big(u_J,\{\varphi_{J,i}^x\}_{i=1\dots N}\Big) &= g_{J}^x,
  \quad \text{on } \Gamma_J^{x}, \label{eqn:transCondRight}
\end{align}
where $u_I$ and $u_J$ are the main fields defined on $\Omega_I$ and $\Omega_J$, respectively.
The auxiliary fields $\{\varphi_{I,i}^x\}_i := \{\varphi_{IJ,i}\}_i$ and $\{\varphi_{J,i}^x\}_i := \{\varphi_{JI,i}\}_i$ defined on the shared interface are governed by equation \eqref{eqn:padeAuxEqn}.
The first set of auxiliary fields is associated to the subproblem defined on $\Omega_I$, and the second set is associated to the one defined on $\Omega_J$.
By using the impedance operator \eqref{eqn:padeImpOp} in equation \eqref{eqn:interface} on both sides of the interface, we have that the transmission variables $g_{I}^x$ and $g_{J}^x$ verify
\begin{align}
  g_{I}^x &= -g_{J}^x + 2B\Big(u_J,\{\varphi_{J,i}^x\}_{i=1\dots N}\Big),
  \quad \text{on } \Gamma_I^{x}, \label{eqn:transVarInt} \\
  g_{J}^x &= -g_{I}^x + 2B\Big(u_I,\{\varphi_{I,i}^x\}_{i=1\dots N}\Big),
  \quad \text{on } \Gamma_J^{x}.
\end{align}
With these transmission variables, the transmission conditions \eqref{eqn:transCondLeft} and \eqref{eqn:transCondRight} enforce weakly the continuity of the main field across the shared interface.

The cross-point treatment consists in enforcing weakly the continuity of auxiliary fields at cross-points.
More precisely, only auxiliary fields defined on edges that are aligned are continuous.
For instance, the auxiliary fields $\{\varphi_{I,j}^{y}\}_j$ defined on $\Gamma_I^{y}$ and the auxiliary fields $\{\varphi_{J,j}^{y}\}_j$ defined on $\Gamma_J^{y}$ (\textit{i.e.}~defined on the upper edges of $\Omega_I$ and $\Omega_J$, respectively, see Figure \ref{fig:shared_egdes_points}) must be equal at $P_I^{x y} = P_J^{x y}$.
Following the approach detailed in \cite{modave2020non}, transmission conditions with specific impedance operators are used to enforce weakly the continuity.
At the cross-point, we use the transmission conditions
\begin{align}
  \partial_{x}\varphi_{I,j}^{y}+ B\Big(\varphi_{I,j}^{y},\{\psi_{I,ij}^{x y}\}_{i=1\dots N}\Big) &= g_{I,j}^{x y},
  \quad \text{at } P_I^{x y}, \\
 -\partial_{x}\varphi_{J,j}^{y}+ B\Big(\varphi_{J,j}^{y},\{\psi_{J,ij}^{x y}\}_{i=1\dots N}\Big) &= g_{J,j}^{x y},
  \quad \text{at } P_J^{x y},
  \label{eqn:applyingOperator_N_onCrossPoint}
\end{align}
for $j=1 \dots N$, with the scalar variables $\{\psi_{I,ij}^{x y}\}_{ij}$ and $\{\psi_{J,ij}^{x y}\}_{ij}$ defined as
\begin{align}
  \psi_{I,ij}^{x y}
  &= -\big[\alpha^2 (c_{j}+1) \varphi_{I,i}^{x} + \alpha^2 (c_i+1) \varphi_{I,j}^{y}\big]\Big/\big[\alpha^2 c_i + \alpha^2 c_{j} + 1\big],
  \quad \text{at } P_I^{x y}, \label{eqn:crossLeft} \\
  \psi_{J,ij}^{x y}
  &= -\big[\alpha^2 (c_{j}+1) \varphi_{J,i}^{x} + \alpha^2 (c_i+1) \varphi_{J,j}^{y}\big]\Big/\big[\alpha^2 c_i + \alpha^2 c_{j} + 1\big],
  \quad \text{at } P_J^{x y}, \label{eqn:crossRight}
\end{align}
for $i,j=1\dots N$.
At the cross-point, the new transmission variables $\{g_{I,j}^{x}\}_j$ and $\{g_{J,j}^{x}\}_j$ verify
\begin{align}
  g_{I,j}^{x y} &= - g_{J,j}^{x y} + 2 B\Big(\varphi_{J,j}^{y}, \{\psi_{J,ij}^{x y}\}_{i=1\dots N}\Big), \quad \text{at } P_I^{x y}, \label{eqn:transVarPnt} \\
  g_{J,j}^{x y} &= - g_{I,j}^{x y} + 2 B\Big(\varphi_{I,j}^{y}, \{\psi_{I,ij}^{x y}\}_{i=1\dots N}\Big), \quad \text{at } P_J^{x y},
\end{align}
for $j=1 \dots N$.
In a nutshell, the same transmission condition with the multivariate function $B$ is used on the interface to couple the main fields (equations \eqref{eqn:transCondLeft}-\eqref{eqn:transCondRight}) and at the cross-points to couple the auxiliary fields (equations \eqref{eqn:crossLeft}-\eqref{eqn:crossRight}).
The scalar variables defined at the corners of the subdomains introduce a coupling of auxiliary fields living on adjacent edges of each subdomain.
This strategy can be adapted rather straightforwardly to deal with boundary cross-points, where interface edges and boundary edges with boundary condition meet.
For further details, we refer to \cite{modave2020non}.

The iterative domain decomposition algorithm is very similar to the algorithm described at the end of section \ref{sec:ddStandard}.
At each iteration, subproblems associated to the subdomains are solved concurrently, and transmission variables are updated.
Here, the subproblem associated to $\Omega_I$ consists in finding the main field verifying system \eqref{eqn:subPbm} and auxiliary fields verifying equations similar to \eqref{eqn:padeAuxEqn} on each interface edge.
The transmission variables are associated to interfaces edges and cross-points.
They are updated with formulas similar to
\begin{equation}
  g_{I}^{x\,(n+1)} = -g_{J}^{x\,(n)} + 2 B\Big(u_J^{(n)}, \{\varphi_{J,i}^{x\,(n)}\}_{i=1\dots N}\Big),
  \quad \text{on } \Gamma_I^{x}
  \label{eqn:gDDMinterf}
\end{equation}
and
\begin{equation}
  g_{I,j}^{x y\,(n+1)} = -g_{J,j}^{x y\,(n)} + 2 B\Big(\varphi_{J,i}^{y\,(n)}, \{\psi_{J,ij}^{x y\,(n)}\}_{i=1\dots N}\Big),
  \quad \text{on } P_I^{x y},
  \label{eqn:gDDMcrossp}
\end{equation}
which are obtained by rewritting equations \eqref{eqn:transVarInt} and \eqref{eqn:transVarPnt} similarly to the general update formula \eqref{eq:updateTransVar}.
Here, $g_{IJ}$ can be defined as
\begin{equation}
\arraycolsep=3.6pt\def\arraystretch{1.5}
g_{IJ} :=
\left[
\begin{array}{c}
  g_{I  }^{x   } \\
  g_{I,j}^{x y } \\
  g_{I,j}^{x y'} \\
\end{array}
\right],
\end{equation}
where $g_{I,j}^{x y'}$ is the transmission variable at $P_I^{xy'}$.
One can consider equations \eqref{eqn:transVarInt} and \eqref{eqn:transVarPnt}, and the formula at $P_I^{xy'}$ that is the same to \eqref{eqn:transVarPnt}, which leads to the following formula:
\begin{equation}
\arraycolsep=3.6pt\def\arraystretch{1.5}
g_{IJ}^{(n+1)} = -g_{JI}^{(n)} +
\left[
\begin{array}{c}
 2 B\Big(u_J^{(n)},               \{\varphi_{J,i}^{x \,(n)}\}_{i=1\dots N}\Big) \\
 2 B\Big(\varphi_{J,i}^{y \,(n)}, \{\psi_{J,ij}^{x y'\,(n)}\}_{i=1\dots N}\Big) \\
 2 B\Big(\varphi_{J,i}^{y'\,(n)}, \{\psi_{J,ij}^{x y'\,(n)}\}_{i=1\dots N}\Big) \\
\end{array}
\right].
\end{equation}
This formula is similar to the general update formula \eqref{eq:updateTransVar}.
Following the approach explained
in section \ref{sec:ddStandard}, all the transmission variables can be merged
into a global vector $\mathbf{g}^{(n+1)}$, and the global process can be recast
as one application of an iterative operator $\boldsymbol{\mathcal{A}}$ on the
vector.  At each iteration $n+1$, the whole process can be seen as one step of
the Jacobi method to solve the linear system
$(\boldsymbol{\mathcal{I}}-\boldsymbol{\mathcal{A}})\mathbf{g}=\mathbf{b}$,
which could be solved with a Krylov subspace iterative method.  Here, the main
difference with most of the works is that the global vector includes
transmission variables associated to both interfaces and cross-points.

\section{Algebraic structure of the interface problem}
\label{Matrix Form of the Iteration Operator}

In this section, we analyze the algebraic structure of the global interface problem, which can be written in an abstract form as
\begin{align}
  \boldsymbol{\mathcal{F}}\mathbf{g} := (\boldsymbol{\mathcal{I}} - \boldsymbol{\mathcal{A}}) \mathbf{g} = \mathbf{b},
  \label{eqn:matrix:abstractSys}
\end{align}
where $\boldsymbol{\mathcal{A}}$ is the iteration matrix, $\mathbf{g}$ is the set of all transmission variables and $\mathbf{b}$ is given by the source term.
The global matrix $\boldsymbol{\mathcal{F}}:=\boldsymbol{\mathcal{I}}-\boldsymbol{\mathcal{A}}$ can be represented as a $N_{\text{dom}} \times N_{\text{dom}}$ sparse block matrix, which each block corresponds to the coupling between the unknowns of two subdomains.
The nature of the blocks is discussed \ref{sec:matrix:localStruct} and the sparse structure of this global block matrix is analyzed in subsection \ref{sec:matrix:globalStruct}.


\subsection{Identification of the blocks}
\label{sec:matrix:localStruct}

Using the block representation, the abstract system \eqref{eqn:matrix:abstractSys} can be rewritten as
\begin{equation}
  \sum_{J=1}^{N_{\text{dom}}} \boldsymbol{\mathcal{F}}_I^J \mathbf{g}_J = \mathbf{b}_I, \quad I=1\dots N_{\text{dom}},
  \label{eqn:matrix:matrixSys}
\end{equation}
where the vectors $\mathbf{g}_I$ and $\mathbf{b}_I$ contain all the transmission variables and the source terms, respectively, for the subdomain $\Omega_I$.
The block $\boldsymbol{\mathcal{F}}_I^J$ corresponds to a coupling between the transmission variables of the subdomains $\Omega_I$ and $\Omega_J$.
The blocks corresponding to subdomains that are not neighbours (\textit{i.e} which do not share any interface edge) are cancelled because there is no direct coupling between the corresponding variables.
Since there are at most four neighbouring subdomains for each subdomains, there are at most four off-diagonal blocks in each line and each column of the global block matrix.

For studying the blocks, we consider a setting with transmission conditions based on the basic impedance operator for the sake of simplicity.
In that case, all the transmission variables are associated to the interface edges.
Since there are two transmission variables per interface edge (one for each neighboring subdomain), the total size of the vectors $\mathbf{g}$ and $\mathbf{b}$ is twice the number of interface edges.
The size of the blocks $\mathbf{g}_I$ and $\mathbf{b}_I$ in these vectors corresponds to the number of interface edges for the subdomain $\Omega_I$.
The block can then be rectangular if the neighbouring subdomains $\Omega_I$ and $\Omega_J$ have different numbers of interface edges.
To simplify the presentation, we assume hereafter that $\Omega_I$ and $\Omega_J$ do not touch the exterior border of the main domain (\textit{i.e.}~each of them has four interface edges and $\boldsymbol{\mathcal{F}}_I^J$ is a $4\times4$ matrix).


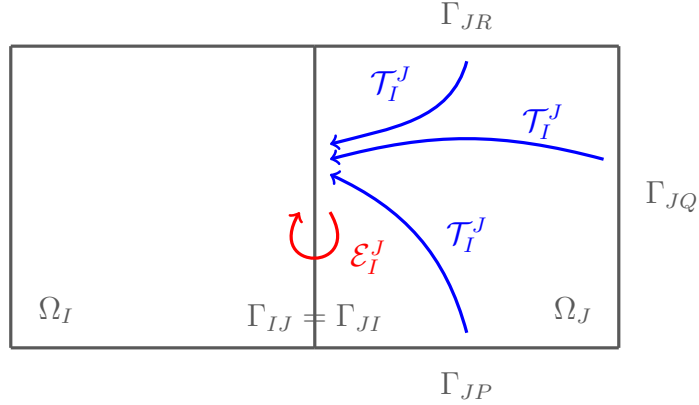
\begin{figure}[!tb]
  \centering
  \captionsetup{font=small, labelfont=bf}
  \small
  \begin{tikzpicture}[scale=1.00]
    
    \coordinate (n0) at (2 ,0) ;
    \coordinate (n1) at (6 ,0) ;
    \coordinate (n2) at (6 ,4) ;
    \coordinate (n3) at (2 ,4) ;
    \coordinate (n4) at (6 ,0) ;
    \coordinate (n5) at (10,0) ;
    \coordinate (n6) at (10,4) ;
    \coordinate (n7) at (6 ,4) ;
    
    \draw[very thick, black!65] (n0)--(n1);
    \draw[very thick, black!65] (n2)--(n3);
    \draw[very thick, black!65] (n3)--(n0);
    \draw[very thick, black!65] (n4)--(n5);
    \draw[very thick, black!65] (n5)--(n6);
    \draw[very thick, black!65] (n6)--(n7);
    \draw[very thick, black!65] (n7)--(n4);
    
    \node[] at( 6,  0.4) {\large\color{black!65} $\Gamma_{IJ}=\Gamma_{JI}$ };
    \node[] at( 8, -0.5) {\large\color{black!65} $\Gamma_{JP}$ };
    \node[] at( 10.7, 2) {\large\color{black!65} $\Gamma_{JQ}$ };
    \node[] at( 8,  4.4) {\large\color{black!65} $\Gamma_{JR}$ };
    
    \node[] at (2.6, 0.5) {\large\color{black!65} $\Omega_{I}$};
    \node[] at (9.4, 0.5) {\large\color{black!65} $\Omega_{J}$};
    
    \draw[->, very thick, blue] ( 8.0, 3.8)  to [out=-105,in=  15, looseness=1] ( 6.2, 2.7);
    \draw[->, very thick, blue] ( 9.8, 2.5)  to [out= 165,in=  15, looseness=1] ( 6.2, 2.5);
    \draw[->, very thick, blue] ( 8.0, 0.2)  to [out= 105,in=- 25, looseness=1] ( 6.2, 2.3);
    \draw[->, very thick, red ] ( 6.2, 1.8)  to [out= -60,in=-120, looseness=6] ( 5.8, 1.8);
    
    \node[blue] at( 7.0, 3.5) {\large $\mathcal{T}_{I}^{J} $};
    \node[blue] at( 9.0, 3.0) {\large $\mathcal{T}_{I}^{J} $};
    \node[blue] at( 8.0, 1.5) {\large $\mathcal{T}_{I}^{J} $};
    \node[red ] at( 6.7, 1.2) {\large $\mathcal{E}_{I}^{J} $};
    
  \end{tikzpicture}
  \caption{Illustration of the coupling introduced by the self-coupling operator (in red) and the transfer operator (in blue) for a configuration with neighboring subdomains $\Omega_I$ and $\Omega_J$ with the shared interface edge $\Gamma_{IJ}=\Gamma_{JI}$.}
  \label{fig:illustration_OT&OE}
\end{figure}

Every line of the system corresponds to a relation similar to equation \eqref{eqn:interface}.
Let us consider the line for the transmission variable $g_{IJ}$, which corresponds to the relation
\begin{equation}
g_{IJ} = - g_{JI} + 2 \, \mathcal{B}_{JI} u_J, \quad \text{on $\Gamma_{IJ}$},
  \label{eqn:line:GammaIF}
\end{equation}
where $J$ is such that $\Omega_I$ and $\Omega_J$ are neighbouring subdomains with the shared interface edge $\Gamma_{IJ}=\Gamma_{JI}$.
By the linearity of the problem, the solution $u_J$ can be slip into two contributions, $u_J=v_J+w_J$.
The field $v_J$ is the solution of subproblem \eqref{eqn:subPbm} for $\Omega_J$ where the right-hand-side term of the Dirichlet boundary condition on $\partial\Omega_J\cap\Gamma_{\text{sca}}$ is cancelled.
The field $w_J$ is the solution of subproblem \eqref{eqn:subPbm} for $\Omega_J$ where the right-hand-side terms of the transmission conditions are cancelled.
Equation \eqref{eqn:line:GammaIF} can then be rewritten as
\begin{equation}
  g_{IJ} = - g_{JI} + 2 \, \mathcal{B}_{JI}  v_J + b_{JI}, \quad \text{on $\Gamma_{IJ}$},
  \label{eqn:line:GammaIF_vw}
\end{equation}
where $b_{JI} := 2 \mathcal{B}_{JI} w_J$ depends only on the data of the problem, with is the incident plane wave is the present case.
In order to exhibit dependences between transmission variables, we decompose the field $v_J$ into several contributions.
For every interface edge $\Gamma_{JK}$ (with $K \in D_J$), we introduce the field $v_{JK}(g_{JK})$ as the solution of subproblem \eqref{eqn:subPbm} for $\Omega_J$ with the transmission variable $g_{JK}$ prescribed on $\Gamma_{JK}$ and where all the other transmission variables and the right-hand-side term of the Dirichlet condition are cancelled.
By linearity, the solution of $v_J$ can then be written as
\begin{equation}
v_J = \sum_{K \in D_J} v_{JK}(g_{JK}),
\end{equation}
Using this decomposition into equation \eqref{eqn:line:GammaIF_vw} gives
\begin{equation}
g_{IJ} + g_{JI} - 2 \sum_{K \in D_J} \mathcal{B}_{JI} v_{JK}(g_{JK}) = b_{JI}, \quad \text{on $\Gamma_{IJ}$}.
\end{equation}
Defining the self-coupling and transfer operators as
\begin{equation}
\begin{aligned}
  \mathcal{E}_{I}^{J} :
    &&
    g_{JI} & \ \longmapsto
    & g_{JI} - 2 \mathcal{B}_{JI} v_{JI}(g_{JI})
    & \ := \ \mathcal{E}_{I}^{J} \, g_{JI}, \\
  \mathcal{T}_{I}^{J} :
    &&
    g_{JK} & \ \longmapsto
    & - 2 \mathcal{B}_{JI} v_{JK}(g_{JK})
    & \ := \ \mathcal{T}_{I}^{J}\, g_{JK},
    && \text{for $K \in D_J, K \neq I$,}
\end{aligned}
\label{eqn:OT&OE}
\end{equation}
we finally have the representation
\begin{equation}
  g_{IJ}  + \mathcal{E}_{I}^{J} g_{JI} + \sum_{K \neq I} \mathcal{T}_{I}^{J} g_{JK} =
  b_{JI}, \quad \text{on $\Gamma_{IJ}$}.
  \label{eqn:interf:representation}
\end{equation}
The self-coupling operator $\mathcal{E}_{I}^{J}$ introduces a coupling between the transmission variables living on the same interface edge $\Gamma_{IJ}=\Gamma_{JI}$, while the transfer operators introduces a coupling between $g_{IJ}$ and the transmission variables living on the other interface edges of $\Omega_J$ (i.e.\ any $\Gamma_{JK}\neq\Gamma_{IJ}$).
These couplings are illustrated in Figure \ref{fig:illustration_OT&OE}.

Thanks to the representation in equation \eqref{eqn:interf:representation}, the elements of the matrix and the right-hand side of the global system can be identified.
The right-hand side of \eqref{eqn:interf:representation} is an element of $\mathbf{b}_I$.
Looking at the first term in the left-hand side, we straightforwardly have that the blocks on the diagonal of the global matrix are identity matrices.
Finally, the second and third terms correspond to elements in the off-diagonal block $\boldsymbol{\mathcal{F}}_I^J$.
The other elements of this block are equal to zero because the relation \eqref{eqn:interf:representation} corresponding to the other edges of $\Omega_I$ do not involve transmission variables of $\Omega_J$.
For instance, there are the neighbouring subdomains $\Omega_I$ and $\Omega_J$, and four neighbouring subdomains of $\Omega_J$ are $\Omega_I$, $\Omega_P$, $\Omega_Q$, $\Omega_R$ (see Figure \ref{fig:illustration_OT&OE}).
Assuming that the shared interface edge is $\Gamma_{IJ}=\Gamma_{JI}$, the block $\boldsymbol{\mathcal{F}}_I^J$ and the vector $\mathbf{g}_J$ read
\begin{equation} 
\arraycolsep=3.6pt\def\arraystretch{1.5}
\boldsymbol{\mathcal{F}}_{I}^{J} =
\left[
\begin{array}{cccc}
  0 & 0 & 0 & 0 \\
  0 & 0 & 0 & 0 \\
  \mathcal{E}_{I}^{J} &
  \mathcal{T}_{I}^{J} &
  \mathcal{T}_{I}^{J} &
  \mathcal{T}_{I}^{J} \\
  0 & 0 & 0 & 0 \\
\end{array}
\right]
\quad\quad\text{and}\quad\quad
\mathbf{g}_J =
\left[
\begin{array}{c}
  g_{JI} \\
  g_{JP}   \\
  g_{JQ} \\
  g_{JR} \\
\end{array}
\right].
\label{eqn:TransmissionMatrx}
\end{equation}
Using equation \eqref{eqn:interf:representation}, one has the relation
\begin{equation} 
\arraycolsep=3.6pt\def\arraystretch{1.5}
\left[
\begin{array}{cccc}
  0, & 0, & g_{IJ}, & 0
\end{array}
\right]^\top
+ \ 
\boldsymbol{\mathcal{F}}_{I}^{J} \ \mathbf{g}_J \ = \ 
\left[
\begin{array}{cccc}
  0, & 0, & b_{JI}, & 0
\end{array}
\right]^\top.
\end{equation}
This example corresponds to the illustration in Figure \ref{fig:illustration_OT&OE}.
In the general case with subdomains touching the exterior border of the main domain, this matrix can be rectangular with numbers of lines and columns between two and four.
Each block always exhibits only one non-zero line, with one self-coupling operator and between one and three transfer operators.


\subsection{Block matrix forms for the global system}
\label{sec:matrix:globalStruct}

\begin{figure}[!tb]
\centering
\begin{tabular}{cc}
\captionsetup{font=small, labelfont=bf}
  \begin{subfigure}[!b]{0.45\textwidth}
  \centering
  \includegraphics[scale = 0.20]{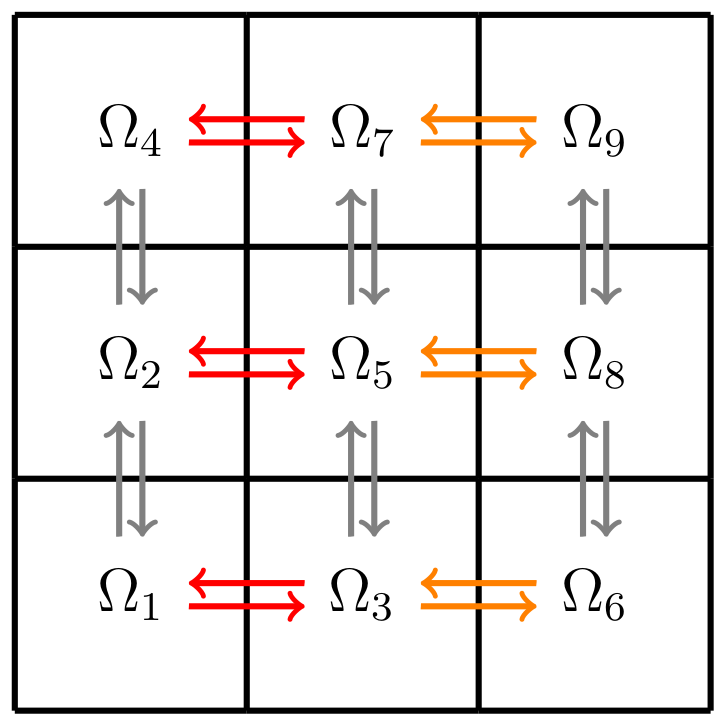} \\
  \caption{Column-type arrangement}
  \label{fig:illustration_CTA}
  \end{subfigure}
&
  \begin{subfigure}[!b]{0.45\textwidth}
  \centering
  \includegraphics[scale = 0.20]{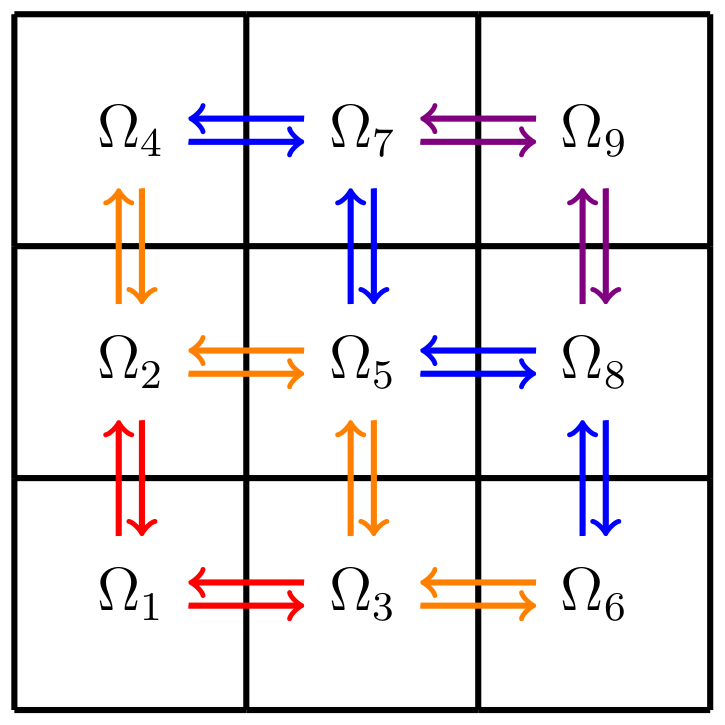} \\
  \caption{Diagonal-type arrangement}
  \label{fig:illustration_DTA}
  \end{subfigure}
\end{tabular}
\caption{Illustration of the subdomains arranged in two different manners with a $3\times3$ checkerboard partition. The colored arrows indicate interactions between groups of subdomains.}
\label{fig:illustration_numbering}
\end{figure}

The sparse structure of the global block matrix consists of all blocks $\boldsymbol{\mathcal{F}}_{I}^{J}$ and identity blocks $\boldsymbol{\mathbf{I}}_{I}$.
With one-dimensional domain partitions, the matrix is block tridiagonal, 
which was leveraged to device efficient sweeping preconditioners (see \textit{e.g.}~\cite{stolk2013rapidly, vion2014double, astaneh2016two}).
With checkerboard partitions, the matrix can also be block tridiagonal if the blocks are arranged correctly.
In order to clearly present this tridiagonal structure which will illustrate horizontal sweeps and diagonal sweeps in the next section, 
the subdomains are arranged in columns and diagonals, which are illustrated in Figure \ref{fig:illustration_numbering} for a $3\times 3$ checkerboard partition.

We first analyze the structure of the global system obtained with the column-type arrangement.
For the $3\times3$ checkerboard partition, the system can be written as
\begin{equation} 
\arraycolsep=3.6pt\def\arraystretch{1.8}
\left[
\begin{array}{ccccccccc}
\cline{1-6}
  \multicolumn{1}{|c}{ \mathbf{I}_1 } &
  { \color{gray} \boldsymbol{\mathcal{F}}^{2}_{1} } &
  &
  \multicolumn{1}{|c}{ \color{red} \boldsymbol{\mathcal{F}}^{3}_{1} } &
  &
  &
  \multicolumn{1}{|c}{ }
\\
  \multicolumn{1}{|c}{ \color{gray} \boldsymbol{\mathcal{F}}_{2}^{1} } &
  \mathbf{I}_2 &
  { \color{gray} \boldsymbol{\mathcal{F}}^{4}_{2} } &
  \multicolumn{1}{|c}{ } &
  { \color{red} \boldsymbol{\mathcal{F}}^{5}_{2} } &
  &
  \multicolumn{1}{|c}{ }
\\ 
  \multicolumn{1}{|c}{ } &
  { \color{gray} \boldsymbol{\mathcal{F}}^{2}_{4} } &
  { \mathbf{I}_4 } &
  \multicolumn{1}{|c}{ } &
  &
  { \color{red} \boldsymbol{\mathcal{F}}^{7}_{4} } &
  \multicolumn{1}{|c}{ }
\\ \cline{1-9}
  \multicolumn{1}{|c}{ \color{red} \boldsymbol{\mathcal{F}}^{1}_{3} } &
  &
  &
  \multicolumn{1}{|c}{ \mathbf{I}_3 } &
  { \color{gray} \boldsymbol{\mathcal{F}}^{5}_{3} } &
  &
  \multicolumn{1}{|c}{ \color{orange} \boldsymbol{\mathcal{F}}^{6}_{3} } &
  &
  \multicolumn{1}{c|}{ }
\\
  \multicolumn{1}{|c}{ } &
  { \color{red} \boldsymbol{\mathcal{F}}^{2}_{5} } &
  &
  \multicolumn{1}{|c}{ \color{gray} \boldsymbol{\mathcal{F}}^{3}_{5} } &
  \mathbf{I}_5 &
  { \color{gray} \boldsymbol{\mathcal{F}}^{7}_{5} } &
  \multicolumn{1}{|c}{ } &
  { \color{orange} \boldsymbol{\mathcal{F}}^{8}_{5} } &
  \multicolumn{1}{c|}{ }
\\
  \multicolumn{1}{|c}{ } &
  &
  { \color{red} \boldsymbol{\mathcal{F}}^{4}_{7} } &
  \multicolumn{1}{|c}{ } &
  { \color{gray} \boldsymbol{\mathcal{F}}^{5}_{7} } &
  { \mathbf{I}_7 } &
  \multicolumn{1}{|c}{ } &
  &
  \multicolumn{1}{c|}{ \color{orange} \boldsymbol{\mathcal{F}}^{9}_{7} }
\\ \cline{1-9}
  &
  &
  &
  \multicolumn{1}{|c}{ \color{orange} \boldsymbol{\mathcal{F}}^{3}_{6} } &
  &
  &
  \multicolumn{1}{|c}{ \mathbf{I}_6 } &
  { \color{gray} \boldsymbol{\mathcal{F}}^{8}_{6} } &
  \multicolumn{1}{c|}{ }
\\
  &
  &
  &
  \multicolumn{1}{|c}{ } &
  { \color{orange} \boldsymbol{\mathcal{F}}^{5}_{8} } &
  &
  \multicolumn{1}{|c}{ \color{gray} \boldsymbol{\mathcal{F}}^{6}_{8} } &
  \mathbf{I}_8 &
  \multicolumn{1}{c|}{ \color{gray} \boldsymbol{\mathcal{F}}^{9}_{8} }
\\ 
  &
  &
  &
  \multicolumn{1}{|c}{ } &
  &
  { \color{orange} \boldsymbol{\mathcal{F}}^{7}_{9} }  &
  \multicolumn{1}{|c}{ } &
  { \color{gray} \boldsymbol{\mathcal{F}}^{8}_{9} }  &
  \multicolumn{1}{c|}{ \mathbf{I}_9 }
\\ \cline{4-9}
\end{array}
\right]
\left[
\begin{array}{c}
\mathbf{g}_{1}     \\
\mathbf{g}_{2}     \\
\mathbf{g}_{4}     \\ \cline{1-1}
\mathbf{g}_{3}     \\
\mathbf{g}_{5}     \\
\mathbf{g}_{7}     \\ \cline{1-1}
\mathbf{g}_{6}     \\
\mathbf{g}_{8}     \\
\mathbf{g}_{9}     \\
\end{array}
\right]
=
\left[
\begin{array}{c}
\mathbf{b}_{1}     \\
\mathbf{b}_{2}     \\
\mathbf{b}_{4}     \\ \cline{1-1}
\mathbf{b}_{3}     \\
\mathbf{b}_{5}     \\
\mathbf{b}_{7}     \\ \cline{1-1}
\mathbf{b}_{6}     \\
\mathbf{b}_{8}     \\
\mathbf{b}_{9}     \\
\end{array}
\right],
\label{eqn:matrix_CTA}
\end{equation}
where $\mathbf{I}_I$ is the identity matrix associated to the subdomain $\Omega_I$.
The global matrix can be rewritten as a $3\times3$ block tridiagonal matrix, where each large block corresponds to interactions between subdomains belonging to two given columns of the domain partition.
Each large block contains $3\times3$ small blocks corresponding to interactions between the subdomains of both columns.
The limits of these large blocks are drawed in equation \eqref{eqn:matrix_CTA}. 
The off-diagonal large blocks, which correspond to interactions between two different columns, are block diagonal.
For a general $N_r \times N_c$ partition of the domain, the structure remains the same.
The global matrix remains a block tridiagonal matrix with $N_c\times N_c$ large blocks, each large block contains $N_r \times N_r$ small blocks, and the off-diagonal large blocks remains block diagonal.

With the diagonal-type arrangement, the structure of the global system can be written as
\begin{equation} 
\arraycolsep=3.6pt\def\arraystretch{1.8}
\left[
\begin{array}{ccccccccc}
  \cline{1-3}
  \multicolumn{1}{|c}{ \mathbf{I}_1 } &
  \multicolumn{1}{|c}{ \color{red} \boldsymbol{\mathcal{F}}^2_1 } &
  { \color{red} \boldsymbol{\mathcal{F}}^3_1 } &
  \multicolumn{1}{|c}{ } &
  &
  &
  &
\\ \cline{1-6}
  \multicolumn{1}{|c}{ \color{red} \boldsymbol{\mathcal{F}}^1_2 } &
  \multicolumn{1}{|c}{ \mathbf{I}_2 } &
  & 
  \multicolumn{1}{|c}{ \color{orange} \boldsymbol{\mathcal{F}}^4_2 } &
  { \color{orange} \boldsymbol{\mathcal{F}}^5_2 } &
  &
  \multicolumn{1}{|c}{ }
\\ 
  \multicolumn{1}{|c}{ \color{red} \boldsymbol{\mathcal{F}}^1_3 } &
  \multicolumn{1}{|c}{ } &
  \mathbf{I}_3 &
  \multicolumn{1}{|c}{ } &
  { \color{orange} \boldsymbol{\mathcal{F}}^5_3 } &
  { \color{orange} \boldsymbol{\mathcal{F}}^6_3 } &
  \multicolumn{1}{|c}{ }
\\ \cline{1-8}
  &
  \multicolumn{1}{|c}{ \color{orange} \boldsymbol{\mathcal{F}}^{2}_{4} } &
  &
  \multicolumn{1}{|c}{ \mathbf{I}_{4} } &
  &
  &
  \multicolumn{1}{|c}{ \color{blue} \boldsymbol{\mathcal{F}}^{7}_{4} } &
  &
  \multicolumn{1}{|c}{ }
\\
  &
  \multicolumn{1}{|c}{ \color{orange} \boldsymbol{\mathcal{F}}^{2}_{5} } &
  { \color{orange} \boldsymbol{\mathcal{F}}^{3}_{5} } &
  \multicolumn{1}{|c}{ } &
  \mathbf{I}_{5} &
  &
  \multicolumn{1}{|c}{ \color{blue} \boldsymbol{\mathcal{F}}^{7}_{5} } &
  { \color{blue} \boldsymbol{\mathcal{F}}^{8}_{5} } &
  \multicolumn{1}{|c}{ }
\\
  &
  \multicolumn{1}{|c}{ } &
  { \color{orange} \boldsymbol{\mathcal{F}}^{3}_{6} }  &
  \multicolumn{1}{|c}{ } &
  &
  \mathbf{I}_{6} &
  \multicolumn{1}{|c}{ } &
  { \color{blue} \boldsymbol{\mathcal{F}}^{8}_{6} } &
  \multicolumn{1}{|c}{ }
\\ \cline{2-9}
  &
  &
  &
  \multicolumn{1}{|c}{ \color{blue} \boldsymbol{\mathcal{F}}^{4}_{7} } &
  { \color{blue} \boldsymbol{\mathcal{F}}^{5}_{7} } &
  &
  \multicolumn{1}{|c}{ \mathbf{I}_7 } &
  &
  \multicolumn{1}{|c|}{ \color{violet} \boldsymbol{\mathcal{F}}^{9}_{7} }
\\
  &
  &
  &
  \multicolumn{1}{|c}{ } &
  { \color{blue} \boldsymbol{\mathcal{F}}^{5}_{8} } &
  { \color{blue} \boldsymbol{\mathcal{F}}^{6}_{8} } &
  \multicolumn{1}{|c}{ } &
  \mathbf{I}_8 &
  \multicolumn{1}{|c|}{ \color{violet} \boldsymbol{\mathcal{F}}^9_8 }
\\ \cline{4-9}
  &
  &
  &
  &
  &
  &
  \multicolumn{1}{|c}{ \color{violet} \boldsymbol{\mathcal{F}}^7_9 } &
  { \color{violet} \boldsymbol{\mathcal{F}}^8_9 } &
  \multicolumn{1}{|c|}{ \mathbf{I}_9 }
\\ \cline{7-9}
\end{array}
\right]
\left[
\begin{array}{c}
\mathbf{g}_{1} \\ \cline{1-1}
\mathbf{g}_{2} \\
\mathbf{g}_{3} \\ \cline{1-1}
\mathbf{g}_{4} \\
\mathbf{g}_{5} \\
\mathbf{g}_{6} \\ \cline{1-1}
\mathbf{g}_{7} \\
\mathbf{g}_{8} \\ \cline{1-1}
\mathbf{g}_{9} \\
\end{array}
\right]
=
\left[
\begin{array}{c}
\mathbf{b}_{1} \\ \cline{1-1}
\mathbf{b}_{2} \\
\mathbf{b}_{3} \\ \cline{1-1}
\mathbf{b}_{4} \\
\mathbf{b}_{5} \\
\mathbf{b}_{6} \\ \cline{1-1}
\mathbf{b}_{7} \\
\mathbf{b}_{8} \\ \cline{1-1}
\mathbf{b}_{9} \\
\end{array}
\right].
\label{eqn:matrix_DTA}
\end{equation}
The matrix of the system can again be written as a block tridiagonal matrix with large blocks corresponding to interactions between two groups of subdomains.
Here, each group corresponds to the subdomains on a given diagonal of the domain partition (see Figure \ref{fig:illustration_DTA}).
In the matrix, the diagonal large blocks are identity matrices because subdomains belonging to the same group are never neighbours.
The off-diagonal large block are rectangular with different sizes because the groups contain different numbers of subdomains.
For a general $N_r\times N_c$ partition of the domain, there are $N_r+N_c-1$ groups of subdomains and the matrix of the system can be still be written as a block tridiagonal matrix. 

Whatever the subdomains are grouped by column, row or diagonal, the global system can be represented with a bloc tridiagonal matrix.
For convenience, we introduce the general representation
\begin{equation} 
\arraycolsep=3.6pt\def\arraystretch{2.2}
\left[
\begin{array}{cccccc}
  \boldsymbol{\mathcal{F}}^{[1]}_{[1]} & { \color{red} \boldsymbol{\mathcal{F}}^{[2]}_{[1]} } \\
  { \color{red} \boldsymbol{\mathcal{F}}^{[1]}_{[2]} } & \boldsymbol{\mathcal{F}}^{[2]}_{[2]} & { \color{orange} \boldsymbol{\mathcal{F}}^{[3]}_{[2]} } \\
  & { \color{orange} \boldsymbol{\mathcal{F}}^{[2]}_{[3]} } & \boldsymbol{\mathcal{F}}^{[3]}_{[3]} & { \color{blue} \boldsymbol{\mathcal{F}}^{[4]}_{[3]} } \\
  & & { \color{blue} \boldsymbol{\mathcal{F}}^{[3]}_{[4]} } & \boldsymbol{\mathcal{F}}^{[4]}_{[4]} & { \color{violet} \boldsymbol{\mathcal{F}}^{[5]}_{[4]} } \\
  & & & { \color{violet} \boldsymbol{\mathcal{F}}^{[4]}_{[5]} } & \boldsymbol{\mathcal{F}}^{[5]}_{[5]}
\end{array}
\right]
\left[
\begin{array}{c}
  \mathbf{g}_{[1]} \\
  \mathbf{g}_{[2]} \\
  \mathbf{g}_{[3]} \\
  \mathbf{g}_{[4]} \\
  \mathbf{g}_{[5]}
\end{array}
\right]
=
\left[
\begin{array}{c}
  \mathbf{b}_{[1]} \\
  \mathbf{b}_{[2]} \\
  \mathbf{b}_{[3]} \\
  \mathbf{b}_{[4]} \\
  \mathbf{b}_{[5]}
\end{array}
\right],
\label{eqn:bigTridiagSys}
\end{equation}
where the vectors $\mathbf{g}_{[I]}$ and $\mathbf{b}_{[I]}$ are associated to one group of subdomains and each block $\boldsymbol{\mathcal{F}}_{[I]}^{[J]}$ corresponds to the coupling between the transmission variables of two groups.
Each block corresponds to one box in equations \eqref{eqn:matrix_CTA} and \eqref{eqn:matrix_DTA}.

\section{Sweeping preconditioners for the interface problem}
\label{Parallel Sweeping Preconditioners}

In this section, we present sweeping preconditioners to accelerate the solution of the interface problem $\bmat{F} \bvec{g} = \bvec{b}$ with standard iterative schemes based on Krylov subspaces.
To be efficient, a preconditioner $\widetilde{\bmat{F}}$ must be designed such that solving the preconditioned problem $\widetilde{\bmat{F}}\!\:^{-1} \bmat{F} \bvec{g} = \widetilde{\bmat{F}}\!\:^{-1} \bvec{b}$ is faster than solving the unpreconditioned problem, and applying the inverse of the preconditioner on any vector is affordable.
With sweeping preconditioners, applying the inverse of $\widetilde{\bmat{F}}$ on a vector corresponds to solving subproblems in a certain order to transfer information following the natural path taken by propagative waves.

Sweeping preconditioners have been proposed for layered partitions of the domain \textit{e.g.}~in \cite{nataf1994optimal,stolk2013rapidly,vion2014double,stolk2017improved,vion2018improved}.
With this kind of partition, the global matrix can be written with a block tridiagonal representation, which each block on the diagonal is an identity matrix and each off-diagonal block is associated to the coupling between two neighboring layers.
Thanks to this structure, the lower and upper triangular parts of the global matrix, which are used in standard Gauss-Seidel and SOR preconditioners, can be explicitly inverted.
Applying the inverse of the lower and upper triangular matrices simply corresponds to solving subproblems following forward and backward sweeps over the subdomains, respectively.
This approach has been used in \cite{vion2018improved} to design various sweeping preconditioners for layered partitions.

We propose an extension of sweeping preconditioners for checkerboard partitions.
Ideas used in \cite{vion2018improved} are applied here by considering the block tridiagonal representation of the global matrix shown in equation \eqref{eqn:bigTridiagSys}.
Although the sweeps are performed over the groups of subdomains, the sweeping directions don't depend on the arrangement of the subdomains.
Column-type and diagonal-type arrangement of the subdomains are only used to illustrate horizontal and diagonal sweeps.
For the sake of clarity, we introduce the groups of subdomains $\Omega_{[S]}$, with $S=1\dots N_{\text{gr}}$, which correspond to columns or diagonals in Figure \ref{fig:illustration_numbering}.
Block symmetric Gauss-Seidel (SGS) and parallel double sweep (DS) preconditioners are described in section \ref{section:SGS} and \ref{section:DS}.
Computational aspects and extensions are discussed in section \ref{section:comput}.



\subsection{Block Symmetric Gauss-Seidel (SGS) preconditioner}
\label{section:SGS}

The general block Symmetric Gauss-Seidel (SGS) preconditioner reads
\begin{equation}
  \widetilde{\bmat{F}}_{\text{SGS}} = (\bmat{D}+\mathring{\bmat{L}}) \bmat{D}^{-1} (\bmat{D}+\mathring{\bmat{U}}),
\end{equation}
where $\bmat{D}$, $\mathring{\bmat{L}}$ and $\mathring{\bmat{U}}$ are respectively the diagonal part, the strictly lower triangular part and the strictly upper triangular part of the block matrix $\mathcal{\bmat{F}}$ represented in equation  \eqref{eqn:bigTridiagSys}.
Assuming there is no coupling between subdomains of the same group, the diagonal part is an identity matrix, $\bmat{D}=\bmat{I}$.
The preconditioner can then be rewritten as $\widetilde{\bmat{F}}_{\text{SGS}}=\bmat{L}\bmat{U}$, with the lower triangular matrix $\bmat{L}=\bmat{I}+\mathring{\bmat{L}}$ and the upper triangular matrix $\bmat{U}=\bmat{I}+\mathring{\bmat{U}}$, which reads
\begin{equation} 
\arraycolsep=3.6pt\def\arraystretch{1.8}
\bmat{L} =
\left[
\begin{array}{cccccc}
  \bmat{I} \\
  \bmat{F}^{[1]}_{[2]} & \bmat{I} \\
  & \bmat{F}^{[2]}_{[3]} & \bmat{I} \\
  & & \ddots & \ddots \\
  & & & \bmat{F}^{[N_{\text{gr}}-1]}_{[N_{\text{gr}}]} & \bmat{I}
\end{array}
\right],
\quad\quad
\bmat{U} =
\left[
\begin{array}{cccccc}
  \bmat{I} & \bmat{F}^{[2]}_{[1]} \\
  & \bmat{I} & \bmat{F}^{[3]}_{[2]} \\
  & & \ddots & \ddots \\
  & & & \bmat{I} & \bmat{F}^{[N_{\text{gr}}]}_{[N_{\text{gr}}-1]} \\
  & & & & \bmat{I}
\end{array}
\right]
.
\end{equation}
To compute explicitly the inverse of the preconditioner, we introduce the matrices $\bmat{L}_{([I],[J])}$ and $\bmat{U}_{([I],[J])}$ that contain the unit diagonal and only one off-diagonal block with the element $\bmat{F}^{[J]}_{[I]}$.
We have
\begin{equation}
  \widetilde{\bmat{F}}_{\text{SGS}}
  = \bmat{L} \:
    \bmat{U}
  = \Big(
    \bmat{L}_{([2],[1])} \:
    \bmat{L}_{([3],[2])} \:
    \cdots \:
    \bmat{L}_{([N_{\text{gr}}],[N_{\text{gr}}-1])} \:
    \Big)
    \Big(
    \bmat{U}_{([N_{\text{gr}}-1],[N_{\text{gr}}])} \:
    \cdots \:
    \bmat{U}_{([2],[3])} \:
    \bmat{U}_{([1],[2])}
    \Big)
\end{equation}
and
\begin{equation}
  \widetilde{\bmat{F}}\!\:^{-1}_{\text{SGS}}
  = \bmat{U}^{-1} \:
    \bmat{L}^{-1}
  = \Big(
    \bmat{U}_{([1],[2])}^{-1} \:
    \bmat{U}_{([2],[3])}^{-1} \:
    \cdots \:
    \bmat{U}_{([N_{\text{gr}}-1],[N_{\text{gr}}])}^{-1} \:
    \Big)
    \Big(
    \bmat{L}_{([N_{\text{gr}}],[N_{\text{gr}}-1])}^{-1} \:
    \cdots \:
    \bmat{L}_{([3],[2])}^{-1} \:
    \bmat{L}_{([2],[1])}^{-1}
    \Big).
\end{equation}
We can easily see that the inverse of the matrices $\bmat{L}_{([I],[J])}$ and $\bmat{U}_{([I],[J])}$ are the same matrices, but with the opposite sign on the off-diagonal block.
Applying the inverse of the preconditioner can be then computed using Algorithm \ref{alg:GeneralSGS}.
The procedure is written more explicitly in Algorithm \ref{alg:DetailPrecond}.

\begin{algorithm}[!tb]
\caption{Application of the SGS preconditioner: $\mathbf{r} \leftarrow \widetilde{\bmat{F}}\!\:^{-1}_{\text{SGS}} \, \mathbf{r}$.}
\label{alg:GeneralSGS}
\setlength{\parskip}{2pt}
  \smallskip
  \texttt{$//$ Forward sweep} \ (application of $\bmat{L}^{-1}$) \\
  \For{$S=1:(N_{\text{gr}}-1)$}{
    $\mathbf{r}_{S+1} \leftarrow \mathbf{r}_{S+1} - \bmat{F}_{[S+1]}^{[S]} \mathbf{r}_{S}$
  }
  \smallskip
  \texttt{$//$ Backward sweep} \ (application of $\bmat{U}^{-1}$) \\
  \For{$S=N_{\text{gr}}:2$}{
    $\mathbf{r}_{S-1} \leftarrow \mathbf{r}_{S-1} - \bmat{F}_{[S-1]}^{[S]} \mathbf{r}_{S}$
  }
\end{algorithm}

\begin{algorithm}[tb!]
\caption{Application of the SGS or DS preconditioner: $\mathbf{r} \leftarrow \boldsymbol{\widetilde{\mathcal{F}}}\!\:^{-1} \, \mathbf{r}$.}
\label{alg:DetailPrecond}
\setlength{\parskip}{2pt}
  \smallskip
  \texttt{$//$ Forward sweep} \\
  \For{$S=1:(N_{\text{gr}}-1)$}{
    \ParFor{each $I$ such that $\Omega_I\subset\Omega_{[S]}$}{
      \smallskip
      Configure boundary data for subproblem on $\Omega_I$: \\
      \quad\quad $u_D \leftarrow 0$ on $\partial\Omega_I\cap\Gamma_D$ \\
      \quad\quad $g_{IK} \leftarrow r_{IK}$ on each interface edge $\Gamma_{IK}\not\subset\partial\Omega_{[S+1]}$ \\
      \quad\quad $g_{IK} \leftarrow r_{IK}$ on each interface edge $\Gamma_{IK}\subset\partial\Omega_{[S+1]}$ \hfill \textit{\color{blue}(If SGS prec.)} \\
      \quad\quad $g_{IK} \leftarrow 0$ on each interface edge $\Gamma_{IK}\subset\partial\Omega_{[S+1]}$ \hfill \textit{\color{blue}(If DS prec.)} \\
      \smallskip
      Compute $u_I$ by solving subproblem on $\Omega_I$. \\
      \smallskip
      Update transmission data for subproblems of the next group: \\
      \quad\quad $r_{JI} \leftarrow r_{JI}-r_{IJ}+2\mathcal{B}_{IJ} u_I$ on each interface edge $\Gamma_{IJ}\subset\partial\Omega_{[S+1]}$ \hfill \textit{\color{blue}(If SGS prec.)} \\
      \quad\quad $r_{JI} \leftarrow r_{JI} + 2\mathcal{B}_{IJ} u_I$ on each interface edge $\Gamma_{IJ}\subset\partial\Omega_{[S+1]}$ \hfill \textit{\color{blue}(If DS prec.)} \\
      \hfill with $J$ and $I$ such that $\Omega_J\subset\Omega_{[S+1]}$ and $\Gamma_{IJ}=\Gamma_{JI}$
    }
  }
  \smallskip
  \texttt{$//$ Backward sweep} \\
  \For{$S=N_{\text{gr}}:2$}{
    \ParFor{each $I$ such that $\Omega_I\subset\Omega_{[S]}$}{
      \smallskip
      Configure boundary data for subproblem on $\Omega_I$: \\
      \quad\quad $u_D \leftarrow 0$ on $\partial\Omega_I\cap\Gamma_D$ \\
      \quad\quad $g_{IK} \leftarrow r_{IK}$ on each interface edge $\Gamma_{IK}\not\subset\partial\Omega_{[S-1]}$ \\
      \quad\quad $g_{IK} \leftarrow r_{IK}$ on each interface edge $\Gamma_{IK}\subset\partial\Omega_{[S-1]}$ \hfill \textit{\color{blue}(If SGS prec.)} \\
      \quad\quad $g_{IK} \leftarrow 0$ on each interface edge $\Gamma_{IK}\subset\partial\Omega_{[S-1]}$ \hfill \textit{\color{blue}(If DS prec.)} \\
      \smallskip
      Compute $u_I$ by solving subproblem on $\Omega_I$. \\
      \smallskip
      Update transmission data for subproblems of the previous group: \\
      \quad\quad $r_{JI} \leftarrow r_{JI}-r_{IJ}+2\mathcal{B}_{IJ} u_I$ on each interface edge $\Gamma_{IJ}\subset\partial\Omega_{[S-1]}$ \hfill \textit{\color{blue}(If SGS prec.)} \\
      \quad\quad $r_{JI} \leftarrow r_{JI} + 2\mathcal{B}_{IJ} u_I$ on each interface edge $\Gamma_{IJ}\subset\partial\Omega_{[S-1]}$ \hfill \textit{\color{blue}(If DS prec.)} \\
      \hfill with $J$ and $I$ such that $\Omega_J\subset\Omega_{[S-1]}$ and $\Gamma_{IJ}=\Gamma_{JI}$
    }
  }
\end{algorithm}

\begin{figure}[!tb]
\centering
\captionsetup{font=small, labelfont=bf}
\small
\begin{subfigure}{\textwidth}
  \centering
  \caption{Horizontal sweeps. Resp. $\bar{\mathcal{U}}_{(3)}^{-1}$, $\bar{\mathcal{U}}_{(2)}^{-1}$}
  \includegraphics[scale = 0.2]{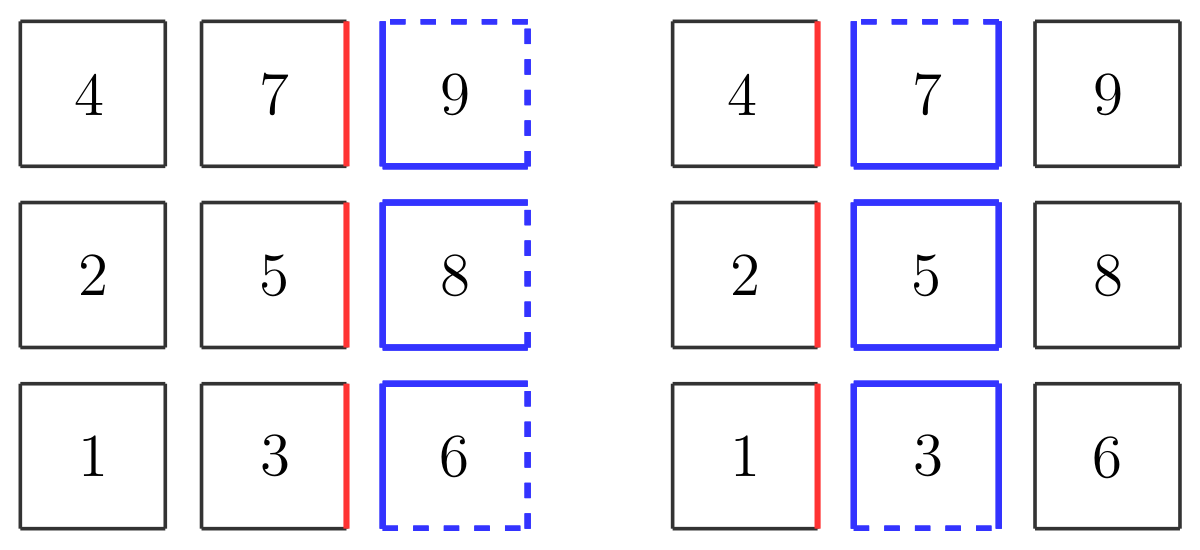} \\
  \label{fig:sub-illustration_SGS_CTA}
\end{subfigure}
\medskip \\
\begin{subfigure}{\textwidth}
  \centering
  \caption{Diagonal sweeps. Resp. $\bar{\mathcal{U}}_{(5)}^{-1}$, $\bar{\mathcal{U}}_{(4)}^{-1}$ and $\bar{\mathcal{U}}_{(3)}^{-1}$}
  \includegraphics[scale = 0.2]{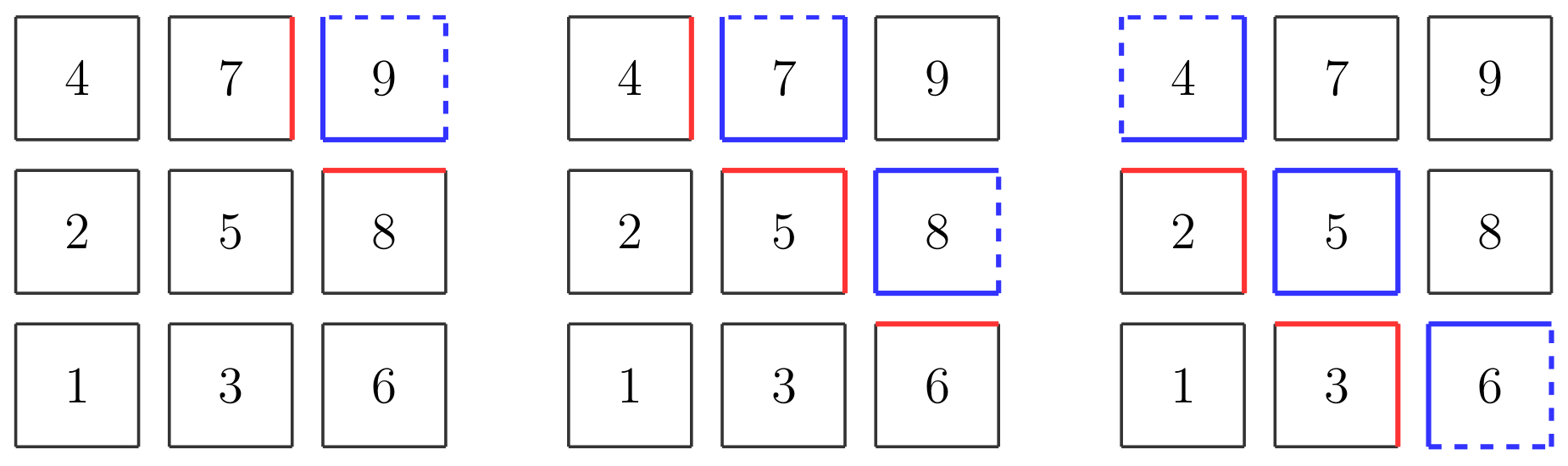} \\
  \label{fig:illustration_SGS_DTA}
\end{subfigure}
\caption{Illustration of the SGS preconditioner with horizontal sweeps and diagonal sweeps. The numbers correspond to the subdomain numbers. Red edges correspond to transmissions variables that are updated. Blue edges correspond to transmission variables that are used in the update formula.}
\label{fig:illustration_SGS}
\end{figure}

In Algorithm \ref{alg:GeneralSGS}, the application of $\bmat{F}_{[S+1]}^{[S]}$ on a vector $\mathbf{r}_{S}$ corresponds to solving subproblems defined on the subdomains of $\Omega_{[S]}$ with the transmission data contained in $\mathbf{r}_{S}$.
The result is used to update transmission variables in $\mathbf{r}_{S+1}$, associated to the subdomains of $\Omega_{[S+1]}$.
Therefore, the first loop, which corresponds to the application of $\bmat{L}^{-1}$, can be interpreted as a forward sweep, where information are propagated across groups of increasing number.
Similarily, the loop corresponding to the application of $\bmat{U}^{-1}$ can be interpreted as a backward sweep, where information are propagated across groups of decreasing number.
The sweeps are performed in the horizontal or diagonal direction not depending on the arrangement of the subdomains.
Iterations of the backward sweep are illustrated on Figure \ref{fig:illustration_SGS} for horizontal sweeps and diagonal sweeps.

In the SGS preconditioner, the forward and backward sweeps must be performed sequentially.
Nevertheless, in each iteration of both loops, the subproblems of a given group of subdomains can be solved in parallel.
With the horizontal sweeps, the update of the transmission variables inside each group have been avoided since the diagonal part $\bmat{D}$ has been replaced with an identity matrix.
With the diagonal sweeps, there is no coupling between the subdomains of the same group since there is no shared edge between them.

\subsection{Parallel Double Sweep (DS) preconditioner}
\label{section:DS}

With the parallel Double Sweep (DS) preconditioner, the $\bmat{L}$ and $\bmat{U}$ matrices are modified in such a way that applying one of them does not modify the elements of the vector used by the other.
The forward and backward sweeps can then be performed in parallel, without data race, which potentially reduces the runtime per iteration by a factor two in parallel environments.

\begin{figure}[!tb]
\centering
\captionsetup{font=small, labelfont=bf}
\small
\begin{subfigure}{\textwidth}
  \centering
  \caption{Horizontal sweeps. Resp. $\bar{\mathcal{U}}_{(3)}^{-1}$, $\bar{\mathcal{U}}_{(2)}^{-1}$}
  \includegraphics[scale = 0.2]{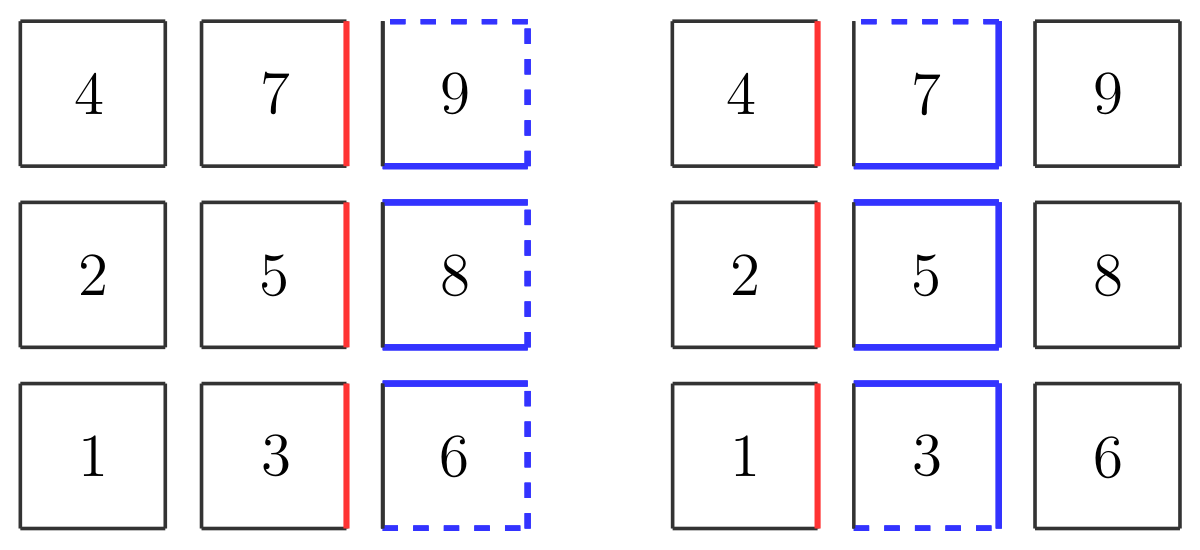} \\
  \label{fig:sub-illustration_DS_CTA}
\end{subfigure}
\medskip \\
\begin{subfigure}{\textwidth}
  \centering
  \caption{Diagonal sweeps. Resp. $\bar{\mathcal{U}}_{(5)}^{-1}$, $\bar{\mathcal{U}}_{(4)}^{-1}$ and $\bar{\mathcal{U}}_{(3)}^{-1}$}
  \includegraphics[scale = 0.2]{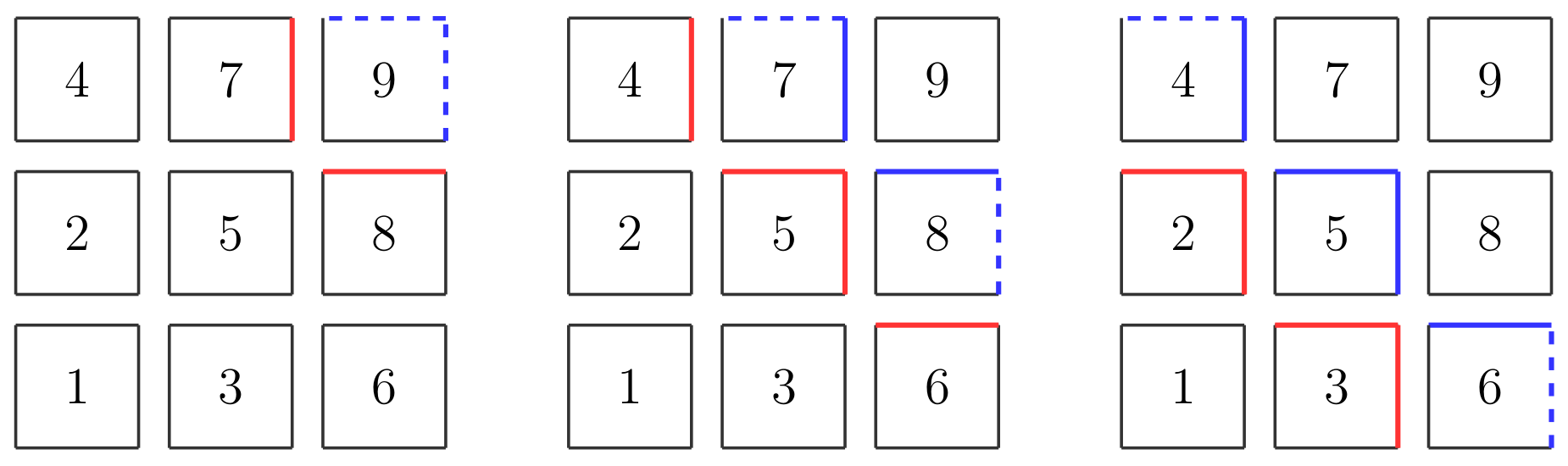} \\
  \label{fig:illustration_DS_DTA}
\end{subfigure}
\caption{Illustration of the DS preconditioner with horizontal sweeps and diagonal sweeps. The numbers correspond to the subdomain numbers. Red edges correspond to transmissions variables that are updated. Blue edges correspond to transmission variables that are used in the update formula.}
\label{fig:illustration:updateDS}
\end{figure}

The DS preconditioner can be written as $\widetilde{\bmat{F}}_{\text{DS}} = \widetilde{\bmat{L}}\,\widetilde{\bmat{U}} = \widetilde{\bmat{U}}\,\widetilde{\bmat{L}} = \widetilde{\bmat{L}} + \widetilde{\bmat{U}} - \bmat{I}$.
To remove the dependences between $\bmat{L}$ and $\bmat{U}$, the blocks are modified in such a way that, for a given group of subdomains $\Omega_{[S]}$, the forward sweep does not use transmission data from edges shared with a subdomain of $\Omega_{[S+1]}$ (these data are modified in the backward sweep), and the backward sweep does not use transmission data from edges shared with a subdomain of $\Omega_{[S-1]}$ (these data are modified in the forward sweep).
The effective update process is illustrated in Figure \ref{fig:illustration:updateDS} for both horizontal and diagonal sweeps.
The application of the DS preconditioner on a vector is detailed in Algorithm \ref{alg:DetailPrecond}.
The main difference with the SGS preconditioner is that one or more transmission variables are cancelled when solving the subdomains.
Thanks to this modification, the forward and backward sweeps can be performed in parallel.

In order to illustrate the modification of the blocks of $\bmat{L}$ and $\bmat{U}$, we consider the $3\times3$ domain partition with the diagonal-type arrangement.
The blocks corresponding to the coupling of $\Omega_{[2]}$ and $\Omega_{[3]}$ reads 
\begin{align}
\boldsymbol{\mathcal{F}}\!\:_{[3]}^{[2]}
=
\arraycolsep=3.6pt\def\arraystretch{1.2}
\left[
\begin{array}{cc}
  \boldsymbol{\mathcal{F}}\!\:_4^2 & 0 \\
  \boldsymbol{\mathcal{F}}\!\:_5^2 & \boldsymbol{\mathcal{F}}\!\:_5^3 \\
  0 & \boldsymbol{\mathcal{F}}\!\:_6^3
\end{array}
\right]
&=
\arraycolsep=3.6pt\def\arraystretch{1.2}
\left[
\begin{array}{cccc|cccc}
  \not\!0& \not\!0 & \not\!0 & \not\!0 & \not\!0 & \not\!0& \not\!0 & \not\!0 \\
  \not\!0& \mathcal{T}_{4}^{2} & \color{black!50}\mathcal{T}_{4}^{2} & \color{black!50}\mathcal{E}_{4}^{2} & 0 & \not\!0& 0 & 0 \\
  \not\!0& 0 & 0 & 0 & 0 & \not\!0& 0 & 0 \\ 
  \not\!0& \not\!0 & \not\!0 & \not\!0 & \not\!0 & \not\!0& \not\!0 & \not\!0 \\ \hline
  \not\!0& \mathcal{T}_{5}^{2} & \color{black!50}\mathcal{E}_{5}^{2} & \color{black!50}\mathcal{T}_{5}^{2} & 0 & \not\!0& 0 & 0 \\
  \not\!0& 0 & 0 & 0 & \mathcal{T}_{5}^{3} & \not\!0& \color{black!50}\mathcal{T}_{5}^{3} & \color{black!50}\mathcal{E}_{5}^{3} \\
  \not\!0& 0 & 0 & 0 & 0 & \not\!0& 0 & 0 \\
  \not\!0& 0 & 0 & 0 & 0 & \not\!0& 0 & 0 \\ \hline
  \not\!0& 0 & 0 & 0 & \mathcal{T}_{6}^{3} & \not\!0& \color{black!50}\mathcal{E}_{6}^{3} & \color{black!50}\mathcal{T}_{6}^{3} \\
  \not\!0& \not\!0 & \not\!0 & \not\!0 & \not\!0 & \not\!0& \not\!0 & \not\!0 \\
  \not\!0& \not\!0 & \not\!0 & \not\!0 & \not\!0 & \not\!0& \not\!0 & \not\!0 \\
  \not\!0& 0 & 0 & 0 & 0 & \not\!0& 0 & 0 \\
\end{array}
\right]
\end{align}
and
\begin{align}
\boldsymbol{\mathcal{F}}\!\:_{[2]}^{[3]}
=
\arraycolsep=3.6pt\def\arraystretch{1.2}
\left[
\begin{array}{ccc}
  \boldsymbol{\mathcal{F}}\!\:_2^4 & \boldsymbol{\mathcal{F}}\!\:_2^5 & 0\\
  0 & \boldsymbol{\mathcal{F}}\!\:_3^5 & \boldsymbol{\mathcal{F}}\!\:_3^6
\end{array}
\right]
&=
\arraycolsep=3.6pt\def\arraystretch{1.2}
\left[
\begin{array}{cccc|cccc|cccc}
  \not\!0& \not\!0 & \not\!0 & \not\!0 & \not\!0 & \not\!0& \not\!0 & \not\!0 & \not\!0 & \not\!0 & \not\!0 & \not\!0 \\
  \not\!0& 0 & 0 & \not\!0& 0 & 0 & 0 & 0 & 0 &\not\!0& \not\!0&  0 \\
  \not\!0& 0 & 0 & \not\!0& \color{black!50}\mathcal{E}_{2}^{5} & \color{black!50}\mathcal{T}_{2}^{5} & \mathcal{T}_{2}^{5} & \mathcal{T}_{2}^{5} & 0 &\not\!0& \not\!0&  0 \\
  \not\!0& \color{black!50}\mathcal{E}_{2}^{4} & \mathcal{T}_{2}^{4} & \not\!0& 0 & 0 & 0 & 0 & 0 &\not\!0& \not\!0&  0 \\ \hline
  \not\!0& 0 & 0 & \not\!0& 0 & 0 & 0 & 0 & \not\!0 &0& 0&  \not\!0 \\
  \not\!0& \not\!0 & \not\!0 & \not\!0 & \not\!0 & \not\!0& \not\!0 & \not\!0 & \not\!0 & \not\!0 & \not\!0 & \not\!0 \\
  \not\!0& 0 & 0 & \not\!0& 0 & 0 & 0 & 0 & \color{black!50}\mathcal{E}_{3}^{6} &\not\!0& \not\!0&  \mathcal{T}_{3}^{6} \\
  \not\!0& 0 & 0 & \not\!0& \color{black!50}\mathcal{T}_{3}^{5} & \color{black!50}\mathcal{E}_{3}^{5} & \mathcal{T}_{3}^{5} & \mathcal{T}_{3}^{5} & 0 &\not\!0& \not\!0&  0
\end{array}
\right],
\end{align}
where rows and columns with $\not\!0$ correspond to items on boundary edges which must be removed.
These blocks belong to $\bmat{L}$ and $\bmat{U}$, respectively.
The modified blocks $\widetilde{\boldsymbol{\mathcal{F}}}\!\:_{[2]}^{[3]}$ and $\widetilde{\boldsymbol{\mathcal{F}}}\!\:_{[3]}^{[2]}$ are obtained by removing the terms in gray.
In $\widetilde{\boldsymbol{\mathcal{F}}}\!\:_{[3]}^{[2]}$, we cancel the terms in the $3^{\text{rd}}$, $4^{\text{th}}$, $7^{\text{th}}$ and $8^{\text{th}}$ columns, corresponding to transmission variables on right edge and top edge of $\Omega_2$ and $\Omega_3$, shared with subdomains of the next group.
In $\widetilde{\boldsymbol{\mathcal{F}}}\!\:_{[2]}^{[3]}$, we cancel the terms in the $2^{\text{nd}}$, $5^{\text{th}}$, $6^{\text{th}}$ and $9^{\text{th}}$ columns, corresponding to transmission variables on bottom edge of $\Omega_4$, left edge and bottom edge of $\Omega_5$ and left edge of $\Omega_6$, shared with subdomains of the previous group.
The modified blocks verify $\widetilde{\boldsymbol{\mathcal{F}}}\!\:_{[2]}^{[3]}\,\widetilde{\boldsymbol{\mathcal{F}}}\!\:_{[3]}^{[2]} = 0$.

\subsection{Flexible preconditioners and parallel aspects}
\label{section:comput}

With both SGS and DS preconditioners, the forward and backward sweeps are performed in a direction that only depends on the blocks $-\boldsymbol{\mathcal{F}}_{[I]}^{[J]}$: a horizontal direction in which the blocks in the same column are performed parallelly, and a diagonal direction in which the blocks in the same diagonal are performed parallelly.
With the standard version of GMRES, only one sweeping direction must be selected.
Nevertheless, in practical situations, it could be advantageous to combine different sweeping directions in order to propagate information more rapidly in different zones of the computational domain.
This can be achieved thanks to the flexible version of GMRES, called F-GMRES \cite{Saad1993,saad2003iterative}, where a different preconditioner can be used at each iteration.
Therefore, the DS and SGS preconditioners can be used with sweeping directions that shall be modified in the course of the iterations, possibly accelerating the convergence of the iterative procedure.

The final computational procedure contains operations that can be performed simultaneously, allowing to use parallel computing architectures.
The forward and backward loops are intrinsically sequential, as the different groups of subdomains have to be treated successively in a specific order.
Nevertheless, parallelism can be found inside each iteration of these loops: Subproblems defined on subdomains of the same group can be solved in parallel with both preconditioners.
In addition, with the DS preconditioner, the forward and backward sweeps can be performed in parallel, as discussed in the previous section.
For applications requiring computations with multiple right-hand sides, strategies can also be used to accelete the computations by solving all the problems in parallel instead of successively.

In distributed-memory parallel environments, novel questions are raised, as the placement of the subdomains on the processors shall influence the communications and the parallel efficiency.
If only one subdomain is placed on each processor, all the processors will be waiting most of the time because of the sequential nature of the sweeping process.
Strategies to improve the parallel efficiency have been discussed in \cite{vion2018improved} for layered-type partitions.
Checkerboard partitions offer novel possibilities, as groups of subdomains can be placed on each processor, and different kind of groups can be chosen.
Placement strategies have been discussed in \cite{taus2020sweeps} in the context of the L-sweeps preconditioners.
Here, for instance, it can be advantageous to place one row of subdomains on each processor when using diagonal sweeps.
This strategy could improve the parallel efficiency by reducing the waiting time of processors.

\section{Computational results}
\label{Benchmarks}

In this section, the preconditionners are studied and compared by using several two-dimensional benchmarks solved with a high-order finite element method.
We consider scattering benchmarks with a single source (section \ref{sec:benchmark1}) and multiple sources (section \ref{sec:benchmark2}), the Marmousi problem (section \ref{sec:benchmark:Marmousi}) and acoustic radiation from engine intake (section \ref{sec:benchmark:radiation}).
The computational results presented in this article have been obtained with a single multi-core processor.
Only the solution of each subproblem was performed using shared-memory parallelism.
The numerical methods have been implemented in a dedicated research code\footnote{Repository: \url{https://gitlab.com/ruiyang/ddmwave}} written in \texttt{C}.
 Gmsh \cite{geuzaine2009gmsh} has been used for mesh generation, domain decomposition, and post-processing.


\subsection{Scattering problem with a single source}
\label{sec:benchmark1}


We consider the scattering of a plane wave by a sound-soft circular cylinder of radius equal to $1$.
The scattered field is computed on a two-dimensional square domain of size $12.5\times12.5$, which is partitioned into an grid of $5\times5$ square subdomains.
The scattering disk is placed in the middle of the subdomain that is at the left-down corner of the grid (first configuration) or the subdomain in the middle of the domain (second configuration).
The Dirichlet boundary condition $u(\mathbf{x}) = -e^{i k x}$ is prescribed on the boundary of the disk, and the Padé-type HABC is prescribed on the exterior border of the computational domain with compatibility conditions at the corners \cite{modave2020corner}.
The Padé-type HABC operator is also used in the transmission conditions prescribed at the interfaces between the subdomains with a suited cross-points treatment (see section~\ref{subsec:crossPointTreatmentHABC}).
Using a large number of auxiliary fields in the Padé-type HABC operator can improve the accuracy of the numerical solution and the efficiency of the DDM for scattering problems \cite{boubendir2012quasi,modave2020corner}.
The number of auxiliary fields $N=8$ and the parameter $\phi = \pi /3$ are used for both exterior and transmission conditions.
The wavenumber $k$ is $2 \pi$ and the wave length $\lambda$ is $1$.

The solution is computed using a standard high-order nodal finite element method on meshes made of triangles and generated with Gmsh \cite{geuzaine2009gmsh}.
Two numerical settings have been considered: P1 finite elements with 20 mesh vertices per wavelength (mesh elements of size $h = 1/20$), and P7 finite elements with 3 elements per wavelength ($h \approx 1/21$).
The meshes of these two settings are made of $97868$ nodes, $3966$ P7 triangles and $70619$ nodes, $139984$ P1 triangles, respectively.
The linear system resulting from the finite element discretization is solved using preconditioned versions of GMRES.
Both the symmetric Gauss-Seidel (SGS) and the parallel double sweep (DS) preconditioners have been tested with three strategies:
\begin{enumerate}
  \item Diagonal sweeps: The forward sweep goes from the left-down corner to the right-up corner of the partition, and the backward sweep does it the other way around.
  {\color{black!50}[SGS-D and DS-D]} 
  \item Horizontal sweeps: The forward sweep goes from the left to the right of the partition, and the backward sweep does it the other way around.
  {\color{black!50}[SGS-H and DS-H]} 
  \item Two diagonal sweeping directions are used in alternance with the flexible version of GMRES (FGMRES). The sweeps are performed between the left-down and right-up corners and between the left-up and right-down corners, alternatively.
  {\color{black!50}[SGS-2D and DS-2D]} 
\end{enumerate}


\subsubsection*{First configuration: Source close to the corner of the domain}

\begin{figure}[p]
  \centering
  \captionsetup{font=small, labelfont=bf}
  \begin{subfigure}{0.9\textwidth}
  \centering
  \includegraphics[scale = 0.135]{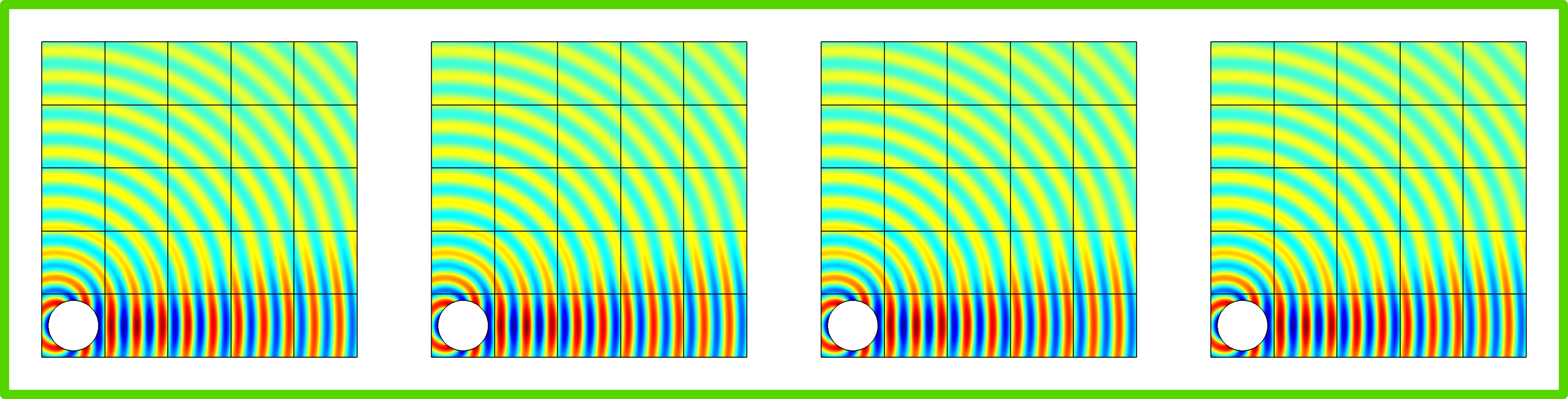} \\
  \caption{SGS preconditioner with diagonal sweeps {\color{black!50}[SGS-D]}}
  \end{subfigure} 
  \begin{subfigure}{0.9\textwidth}
  \centering
  \includegraphics[scale = 0.135]{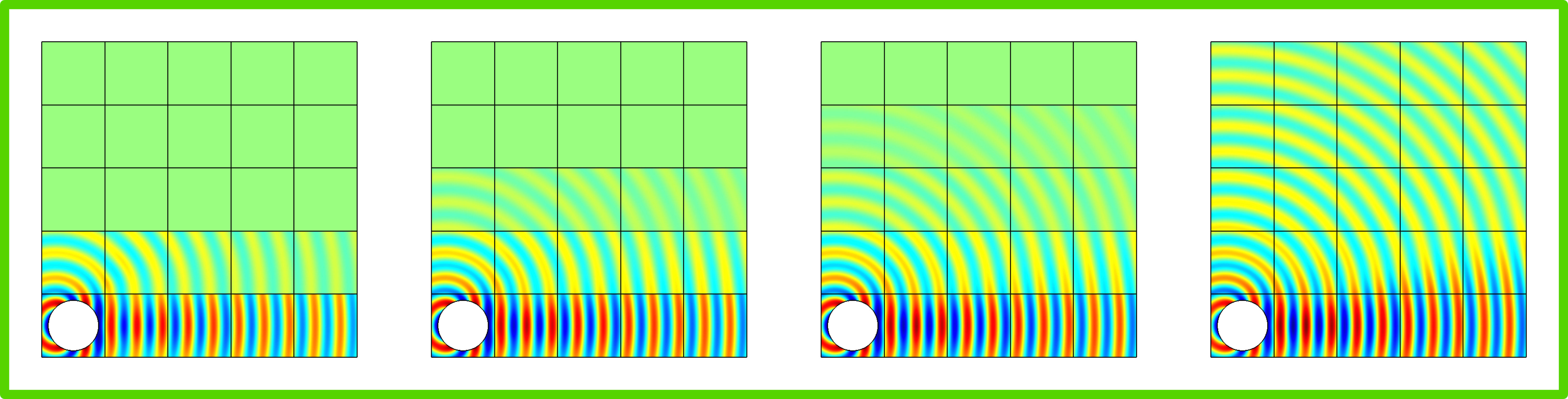} \\
  \caption{SGS preconditioner with horizontal sweeps {\color{black!50}[SGS-H]}}
  \end{subfigure} 
  \begin{subfigure}{0.9\textwidth}
  \centering
  \includegraphics[scale = 0.135]{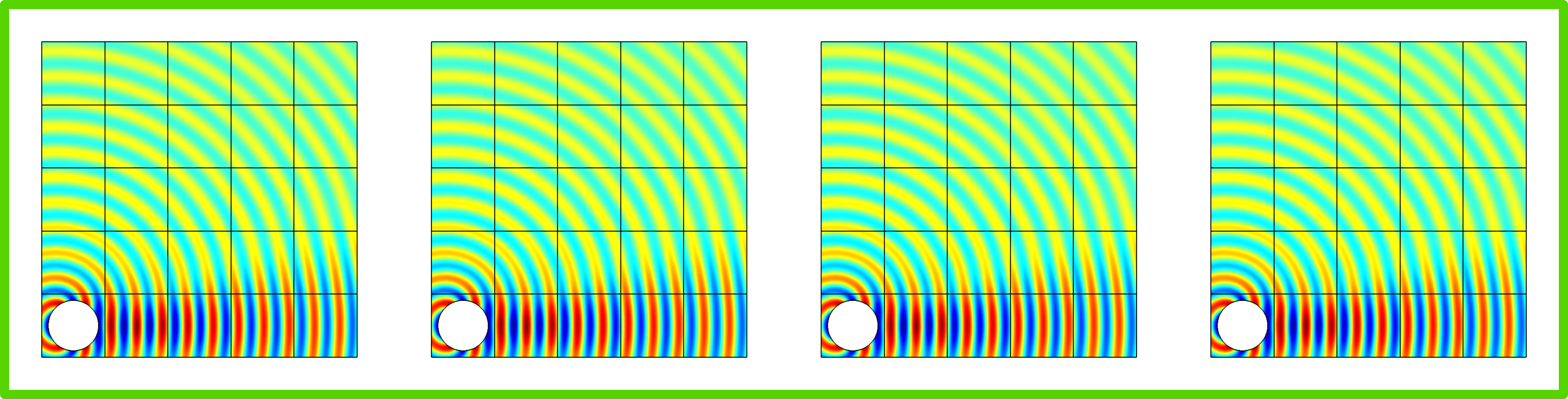} \\
  \caption{SGS preconditioner with F-GMRES and alternating diagonal sweeps {\color{black!50}[SGS-2D]}}
  \end{subfigure} 
  \begin{subfigure}{0.9\textwidth}
  \centering
  \includegraphics[scale = 0.135]{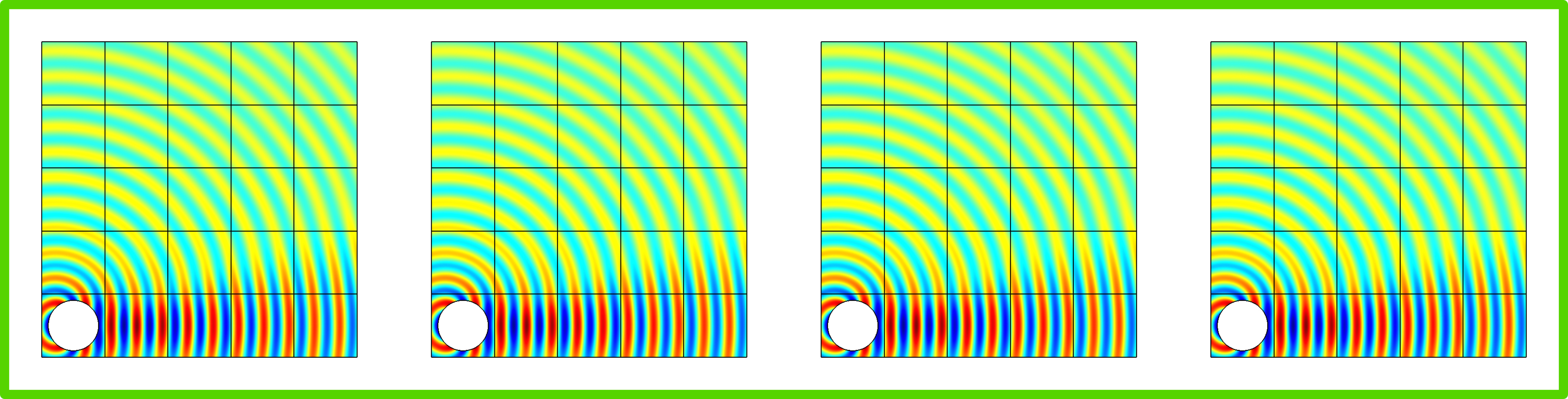} \\
  \caption{DS preconditioner with diagonal sweeps {\color{black!50}[DS-D]}}
  \end{subfigure}
  \begin{subfigure}{0.9\textwidth}
  \centering
  \includegraphics[scale = 0.135]{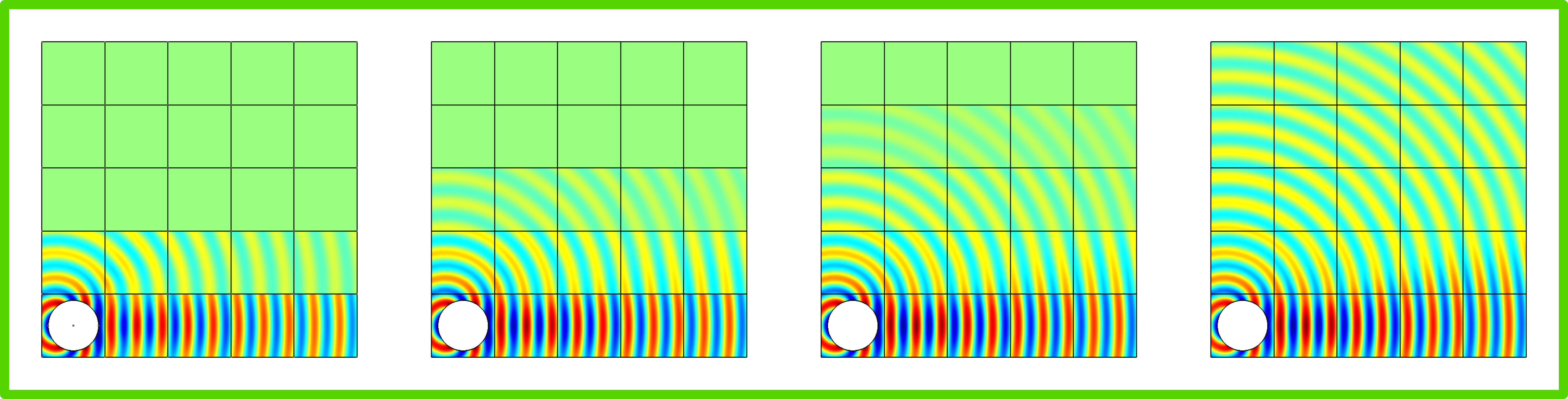} \\
  \caption{DS preconditioner with horizontal sweeps {\color{black!50}[DS-H]}}
  \end{subfigure}
  \begin{subfigure}{0.9\textwidth}
  \centering
  \includegraphics[scale = 0.135]{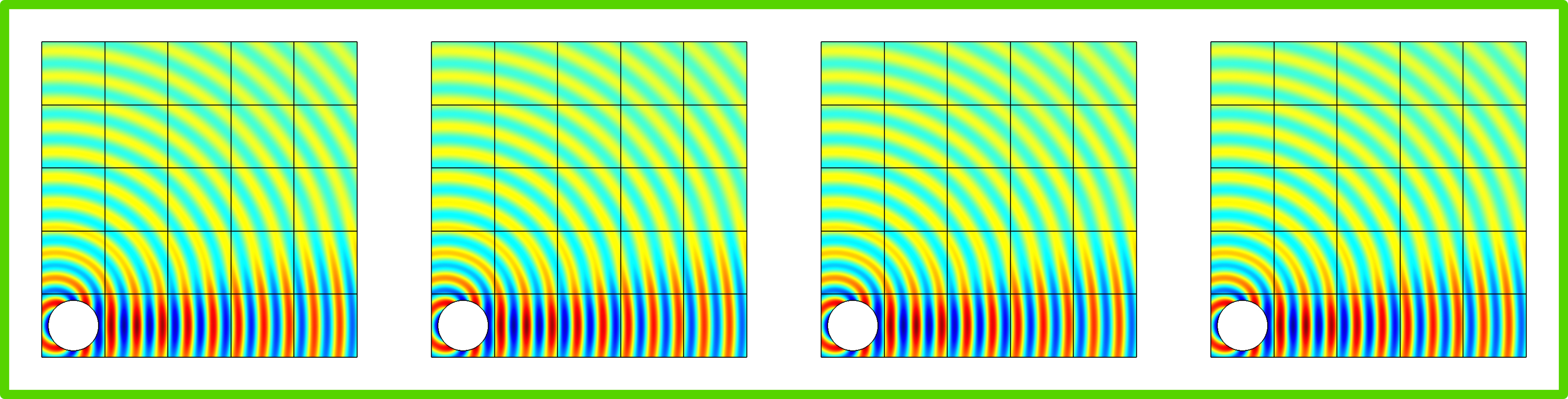} \\
  \caption{DS preconditioner with F-GMRES and alternating diagonal sweeps {\color{black!50}[DS-2D]}}
  \end{subfigure}
  \caption{Scattering problem with a single source (first configuration). Snapshot of the solution after 1, 2, 3 and 4 GMRES iterations with the different preconditioners.}
  \label{fig:single_corner}
\end{figure}

In the first configuration, the scattering disk is placed in the subdomain that is at the left-down corner of the grid.
Figure \ref{fig:single_corner} shows snapshots of the solutions after the first GMRES iterations with the different preconditioners.
For all the preconditioners with diagonal sweeps, the solution is already good after only one iteration, which is due to two successful strategies for this specific case.
First, the HABC operator used in the transmission conditions is particularly well-suited for scattering benchmarks.
Then, the first sweep over the subdomains goes from the left-bottom corner to right-up corner, which follows the natural behavior of waves in this benchmark.
Therefore, we get all correct information in all subdomains after the first iteration.
For the preconditioners with horizontal sweeps, relying on round-trips between left boundary and right boundary, we can see that we get the complete solution only at fourth iteration.

\begin{figure}[!tb]
  \centering
  \small
  \captionsetup{font=small, labelfont=bf}
  \begin{subfigure}[b]{0.47\textwidth}
  \includegraphics[scale = 0.4]{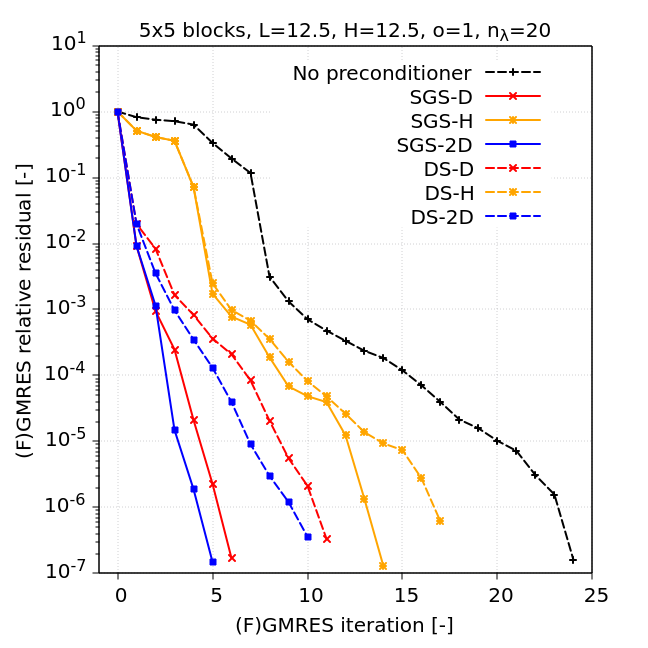} \\
  \end{subfigure}  
  \begin{subfigure}[b]{0.47\textwidth}
  \includegraphics[scale = 0.4]{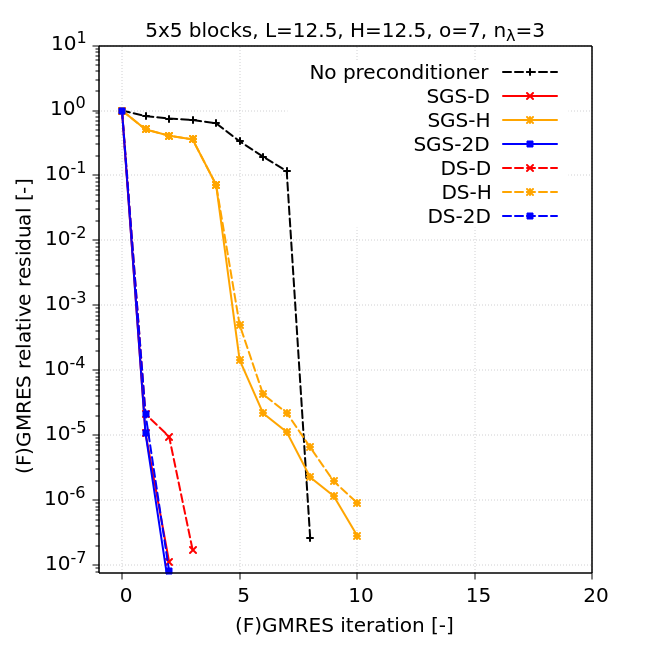} \\
  \end{subfigure}
  \caption{Scattering problem with a single source (first configuration). Residual history with/without preconditioner with P1 elements and 20 mesh vertices by wavelength (left) and with P7 elements and 3 elements by wavelength (right). Preconditioners with diagonal sweeps (SGS-D/DS-D), horizontal sweeps (SGS-H/DS-H) and alternating diagonal sweeps (SGS-2D/DS-2D) are considered.}
  \label{fig517}
\end{figure}

The residual histories obtained with the different preconditioners are shown in Figure \ref{fig517} for finite element schemes with P1 and P7 elements, respectively.
These results confirm the visual interpretations.
In both cases, the relative residual suddenly drops in residual at the first iteration when a preconditionner with diagonal sweeps is used (\textit{i.e.}~SGS-D, SGS-2D, DS-D and DS-2D), while it happens at the fourth iteration with horizontal sweeps (\textit{i.e.}~SGS-H and DS-H).
Without preconditioner, eight iterations are required to reach the sudden drop in residual because the waves have to go through eight subdomains to travel from the left-down corner to the right-up corner.

By comparing the results obtained with P1 and P7 elements, we observe that the residual histories are very similar before the sudden drop in residual with both kinds of finite elements.
Nevertheless, the sudden drops are much sharper and the residuals decrease more rapidly after the sudden drop with the P7 elements.
After the sudden drop in residual, we observe that the decay of the residual is nearly twice slower with the DS preconditioner than with the SGS preconditioner when diagonal sweeps is used.
This must be balanced with the fact that the SGS preconditioner is intrinsically a sequential procedure, while the DS preconditioner relies on two sweeps that can be done in parallel.


\subsubsection*{Second configuration: Source in the middle of the domain}

\begin{figure}[p]
  \centering
  \captionsetup{font=small, labelfont=bf}
  \begin{subfigure}{\textwidth}
  \centering
  \includegraphics[scale = 0.135]{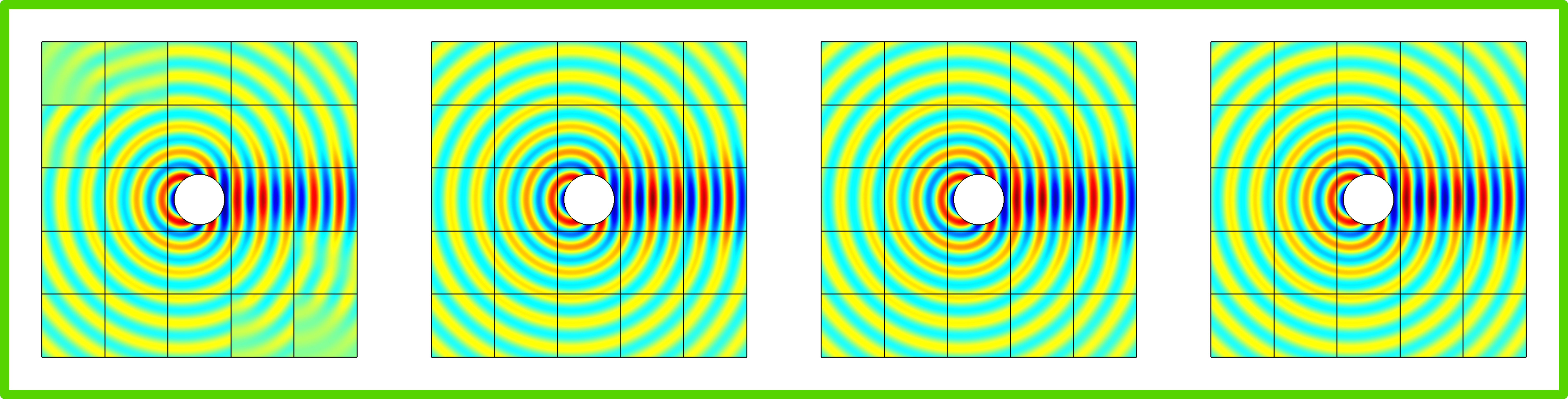} \\
  \caption{SGS preconditioner with diagonal sweeps {\color{black!50}[SGS-D]}}
  \end{subfigure}
  \begin{subfigure}{\textwidth}
  \centering
  \includegraphics[scale = 0.135]{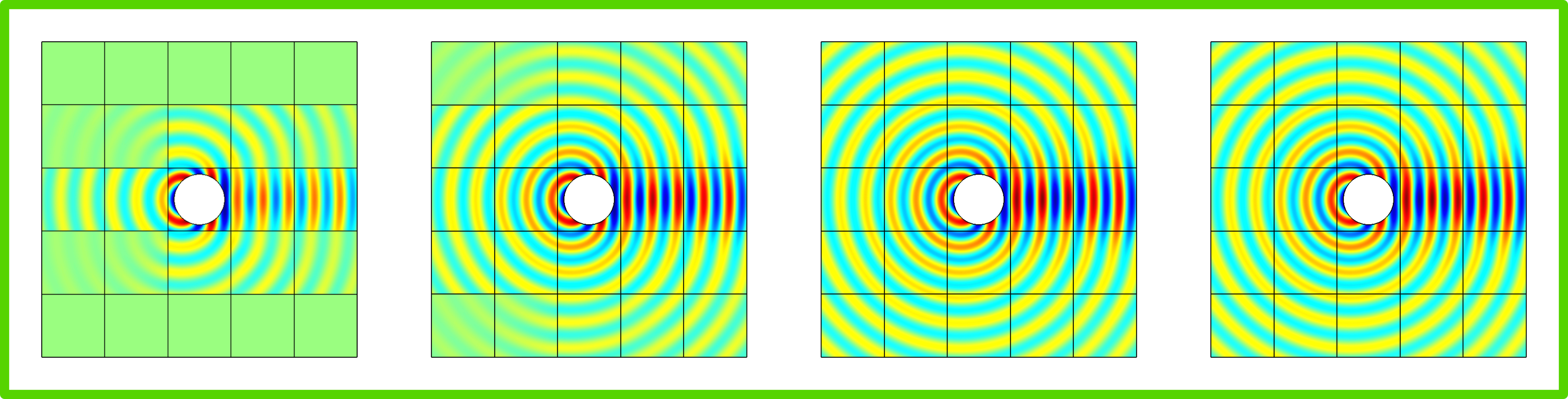} \\
  \caption{SGS preconditioner with horizontal sweeps {\color{black!50}[SGS-H]}}
  \end{subfigure}
  \begin{subfigure}{\textwidth}
  \centering
  \includegraphics[scale = 0.135]{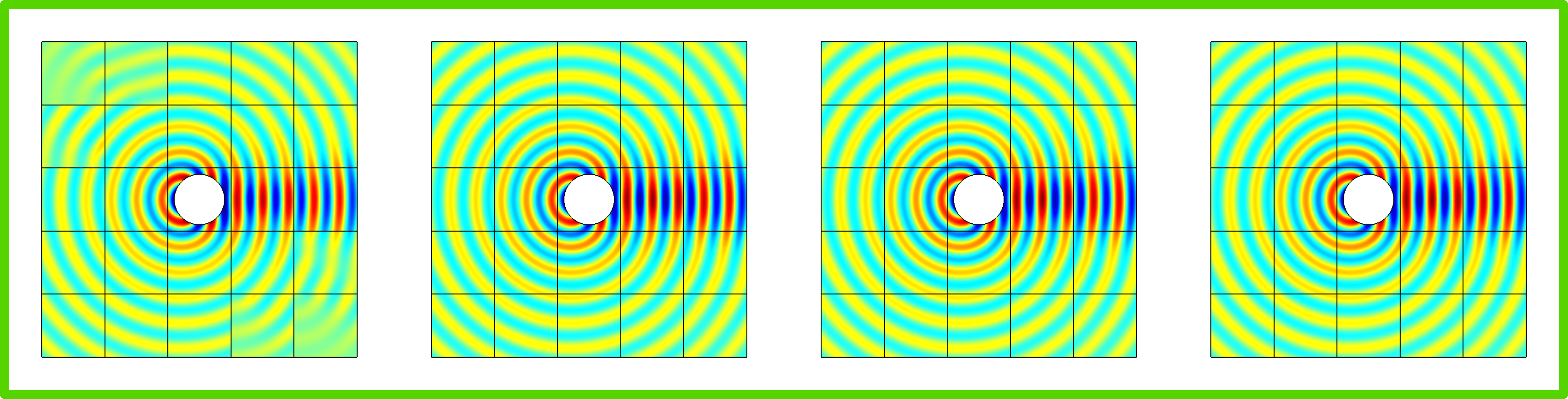} \\
  \caption{SGS preconditioner with F-GMRES and alternating diagonal sweeps {\color{black!50}[SGS-2D]}}
  \end{subfigure}
  \begin{subfigure}{\textwidth}
  \centering
  \includegraphics[scale = 0.135]{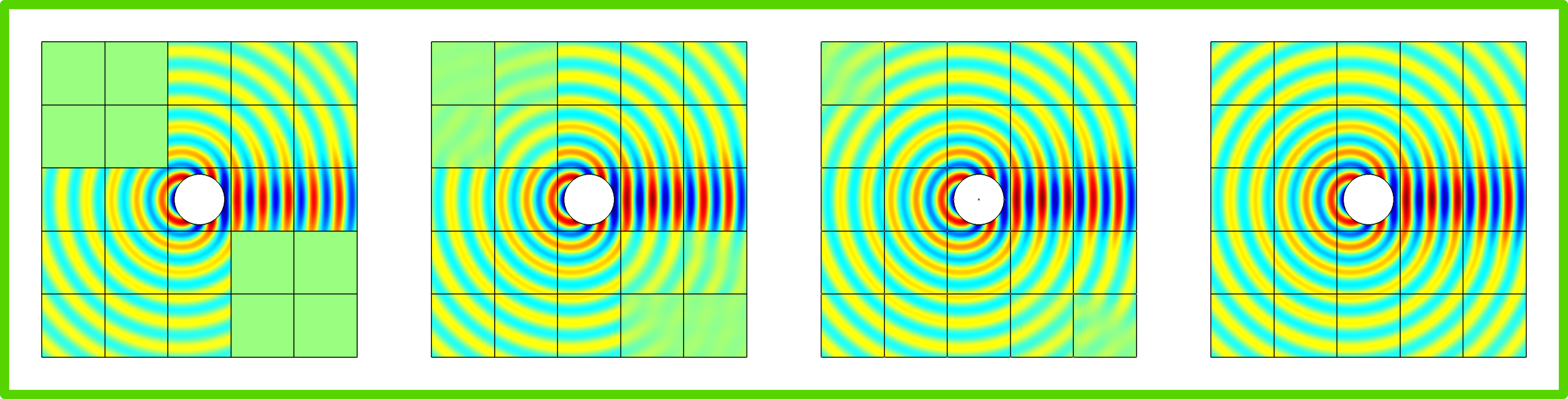} \\
  \caption{DS preconditioner with diagonal sweeps {\color{black!50}[DS-D]}}
  \end{subfigure}
  \begin{subfigure}{\textwidth}
  \centering
  \includegraphics[scale = 0.135]{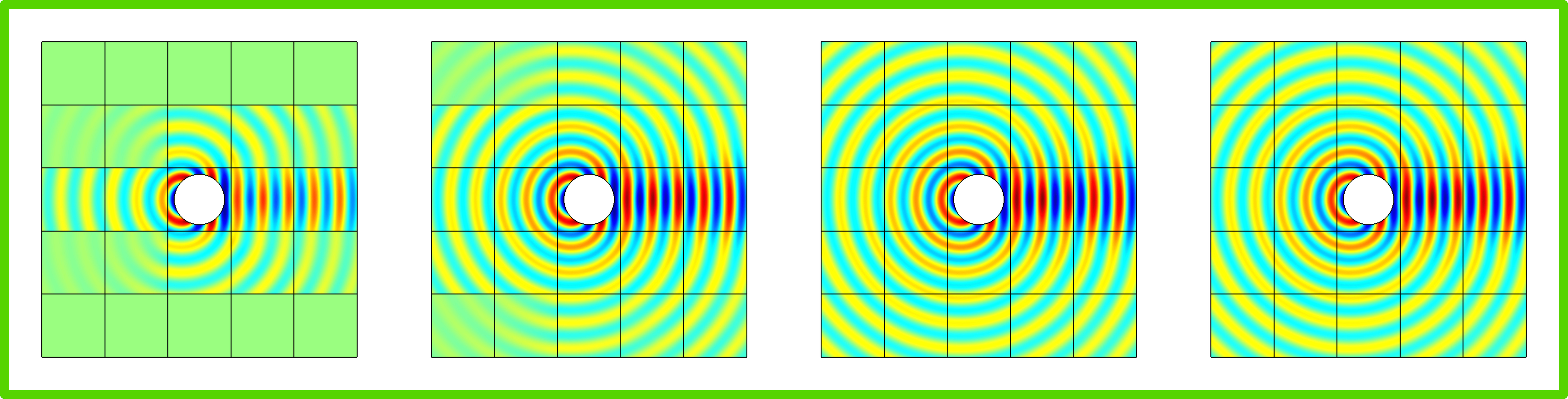} \\
  \caption{DS preconditioner with horizontal sweeps {\color{black!50}[DS-H]}}
  \end{subfigure}
  \begin{subfigure}{\textwidth}
  \centering
  \includegraphics[scale = 0.135]{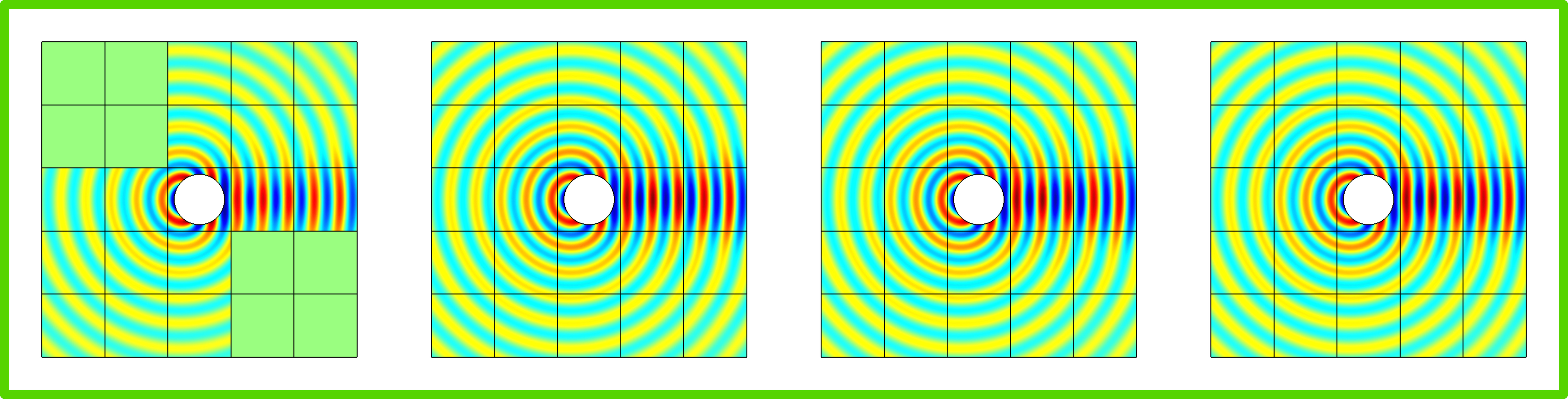} \\
  \caption{DS preconditioner with F-GMRES and alternating diagonal sweeps {\color{black!50}[DS-2D]}}
  \end{subfigure}
  \caption{Scattering problem with a single source (second configuration). Snapshot of the solution after 1, 2, 3 and 4 GMRES iterations with the different preconditioners.}
  \label{fig:single_arbitrary}
\end{figure}

In the second configuration, the scattering disk is placed in the middle of the computational domain.
Figure \ref{fig:single_arbitrary} shows snapshots of the solutions after the first GMRES iterations with the different preconditioners.
First, let us focus on the solutions obtained after the very first iteration.
On the snapshots, we see that the zone of influence of the source corresponds to subdomains that are mainly along the diagonal direction or along the horizontal direction starting from the center of the domain.
The propagation of the source in these subdomains is obviously related to the use of preconditioners with diagonal sweeps and horizontal sweeps of the subdomains, respectively. 

\begin{figure}[!tb]
  \centering
  \small
  \captionsetup{font=small, labelfont=bf}
  \begin{subfigure}[b]{0.47\textwidth}
  \includegraphics[scale = 0.4]{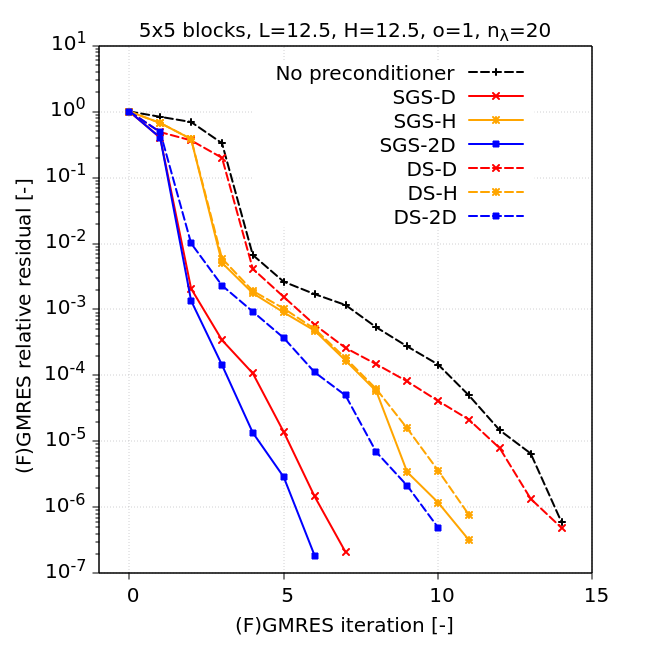} \\
  \end{subfigure}
  \begin{subfigure}[b]{0.47\textwidth}
  \includegraphics[scale = 0.4]{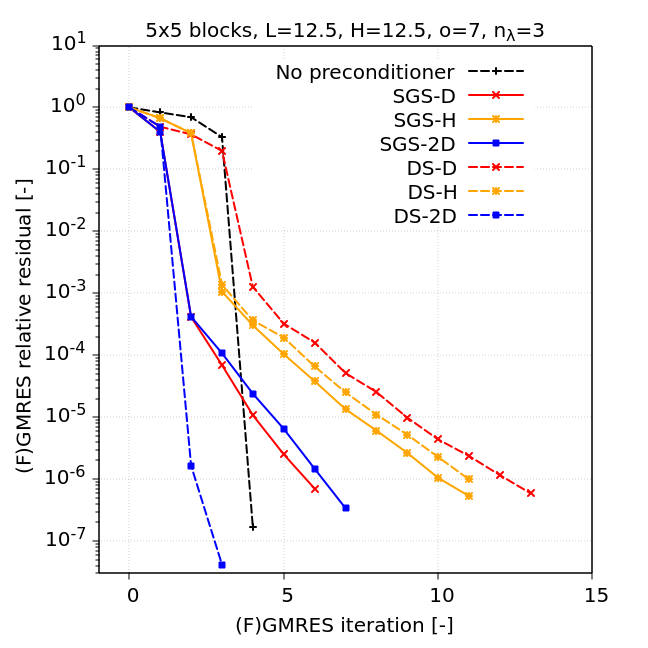} \\
  \end{subfigure}
  \caption{Scattering problem with a single source (second configuration). Residual history with/without preconditioner with P1 elements and 20 mesh vertices by wavelength (left) and with P7 elements and 3 elements by wavelength (right). Preconditioners with diagonal sweeps (SGS-D/DS-D), horizontal sweeps (SGS-H/DS-H) and alternating diagonal sweeps (SGS-2D/DS-2D) are considered.}
  \label{fig519}
\end{figure} 

By contrast with the first configuration, we have different results when the SGS and DS preconditioners are used with the diagonal sweeps or the F-GMRES strategy.
Because the forward sweeps and backward sweeps of DS-D do not affect each other, the subdomains on the right-up and left-down corners cannot be reached by the source after the first iteration, which is visible on Figures \ref{fig:single_arbitrary} (d) and (f).
By contrast, these subdomains can be reached during the backward sweep of the SGS preconditioner, which is performed after the forward sweep (Figures \ref{fig:single_arbitrary} (a) and (c)).
With diagonal sweeps, four iterations are required to get a good solution in all the subdomains with the DS preconditioner, while only two iterations are required with the SGS preconditioner.
When the strategy with F-GMRES is used, the solution is very good at the second iteration, even with the DS preconditioner, because the sweeps of successive iterations are performed along both diagonal directions (Figures \ref{fig:single_arbitrary} (b) and (e)).
 
The residual histories obtained with the different preconditioners are shown in Figure \ref{fig519} for finite element schemes with P1 and P7 elements, respectively.
In all the cases, the relative residual decreases slowly until a sudden drop, which happens when the source has been propagated in all the subdomains, and when the numerical solution is close to the converged solution.
The results confirm the visual observations.
With the SGS preconditioner, the sudden drop occurs at the second iteration if the diagonal sweeps or the stategy with F-GMRES is used, and a third iteration is required with horizontal sweeps.
With the DS preconditioner, four iterations are necessary with diagonal sweeps, while alternating diagonal sweeps realized by F-GMRES still requires two iterations.

When the source is placed in an arbitrary position in the computational domain, using flexible preconditioners, with several alternative sweeping directions, can be much more suitable than fixed preconditioners.
Here, both SGS and DS preconditioners perform very well with the alternating diagonal sweeps.
With P1 finite elements, the number of iterations to get the relative residual $10^{-6}$ is twice larger with the DS preconditioner than with the SGS preconditioner.
Because the sweeps of the DS preconditioner can be done in parallel, providing a speed up of 2, both SGS and DS approaches seem equivalent.
By contract, the number of iterations is similar with P7 finite elements.
The parallel DS preconditioner then is more interesting for that case.

\subsection{Scattering problem with multiple sources}
\label{sec:benchmark2}

\begin{figure}[!tb]
  \centering
  \small
  \includegraphics[scale = 0.20]{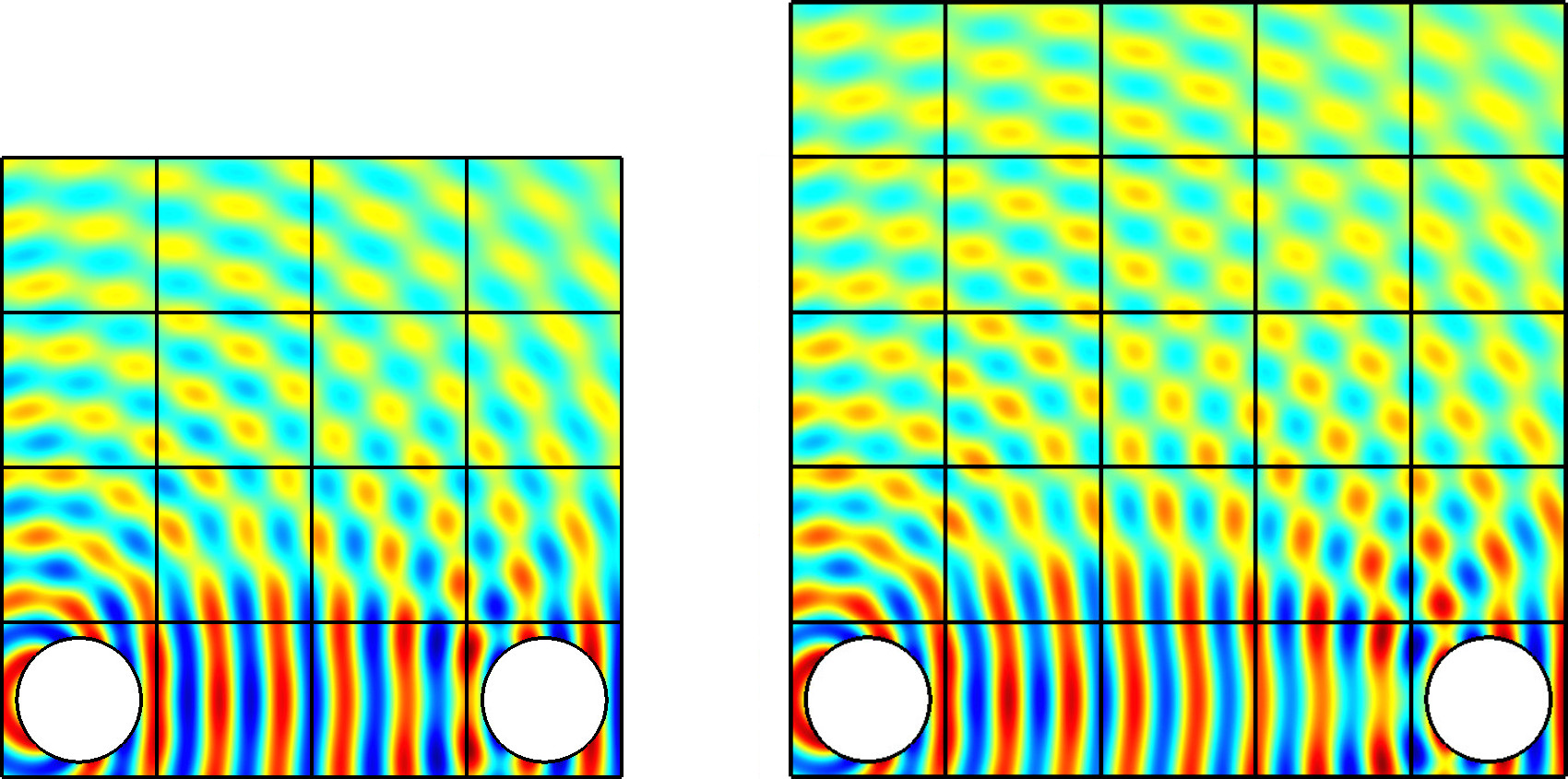} \\
  \captionsetup{font=small, labelfont=bf}
  \caption{Scattering problem with two sources. Snapshot of the solution for configurations with $4 \times 4$ and $5 \times 5$ subdomains.}
  \label{fig624}
\end{figure}

In this section, we consider the scattering of a plane wave by two sound-soft circular cylinders of unit radius.
This problem is more challenging for the DDM than the previous one because the multiple reflections between both obstacles can be complicated to capture.
The simulations are performed over square grids of $N_r\times N_c$ subdomains, with $N_r=N_c=4$, $8$, $12$, $16$ and $20$.
The dimension of each subdomain is $2.5\times2.5$.
For each configuration, the scattering disks are placed at the left-down corner and the right-down corner of the grid.
The physical and numerical parameters are the same as in the previous section.
The numerical solutions for the $4\times4$ and $5\times5$ configurations are shown in Figure \ref{fig624}.

\begin{figure}[!tb]
  \centering
  \small
  \captionsetup{font=small, labelfont=bf}
  \begin{subfigure}[b]{0.47\textwidth}
  \includegraphics[scale = 0.4]{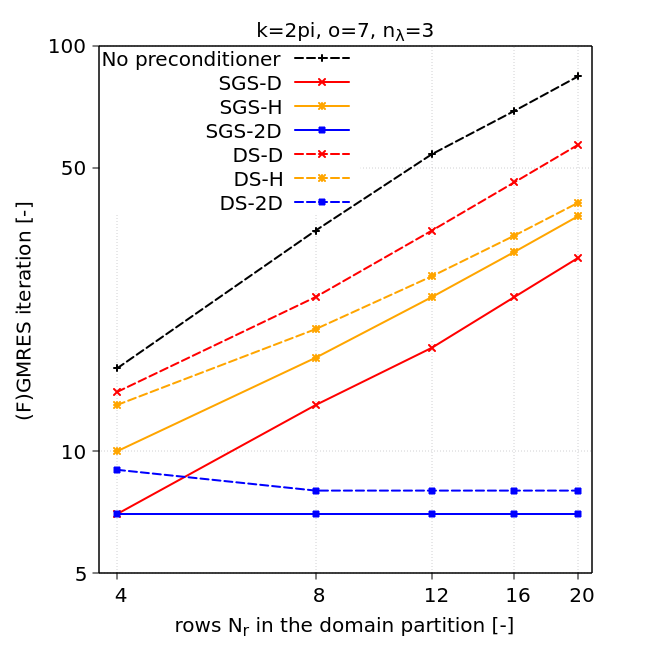}
  \end{subfigure} \quad
  \begin{subfigure}[b]{0.47\textwidth}
  \includegraphics[scale = 0.4]{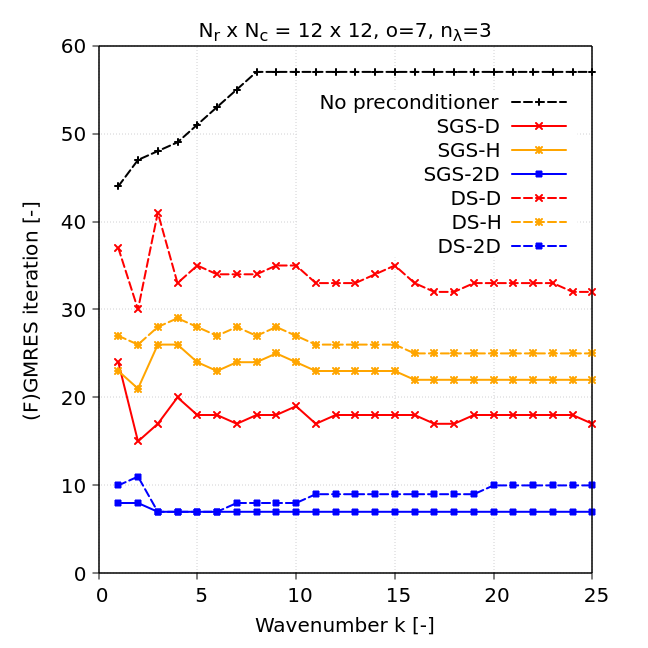}
  \end{subfigure}
  \caption{Scattering problem with two sources. Number of
    iterations with preconditioned GMRES and F-GMRES (without restart) to reach
    the relative residual $10^{-6}$ as a function of the number of rows $N_r$ in
    the domain decomposition, and as a function of the wavenumber $k$ for a
    given mesh density with 3 elements per wavelength.}
  \label{fig:nItVSnp_multObstacles}
\end{figure}

The numbers of GMRES iterations to reach a relative residual $10^{-6}$ with the different preconditioners are given in Figure \ref{fig:nItVSnp_multObstacles} for different domain partitions.
Contrary with the single obstacle case, the physical solution cannot be obtained in only one or two sweeps anymore because of the multiple reflections between both obstacles.
We observe that the number of iterations increases with the number of subdomains when the preconditioners are used with a fixed horizontal or diagonal sweeps.
By contrast, the strategy with the alternating diagonal sweepss keep the number of iterations constant, both with SGS and DS preconditioners.
This is mainly due to the position of the disks: each of the two alternating diagonal sweeps is well suited to one of the sources.
The number of iteration is slightly lower with the SGS preconditioner than with the DS preconditioner, but the latter will be more interesting in parallel environments considering that the forward and backward sweeps can be performed concurrently.

Next, we study the efficiency of the preconditioners with respect to the wavenumber $k$. The high frequency cases are very challenging.
The number of iterations to reach a relative residual $10^{-6}$ is plotted as a function of the wavenumber $k$ in Figure \ref{fig:nItVSnp_multObstacles}.
In all the cases, meshes with 3 elements per wavelength are used.
We observe that, for all preconditioners, the influence of wavenumber $k$ on the rate of convergence is not significant when $k$ varies from $5$ to $25$. The rate of convergence is not stable when $k < 5$, where the meshes are very coarse and the boundaries of the disks are not well represented.

\subsection{Marmousi benchmark}
\label{sec:benchmark:Marmousi}

\begin{figure}[!tb]
  \centering
  \small
  \captionsetup{font=small, labelfont=bf}
  \includegraphics[scale = 0.25]{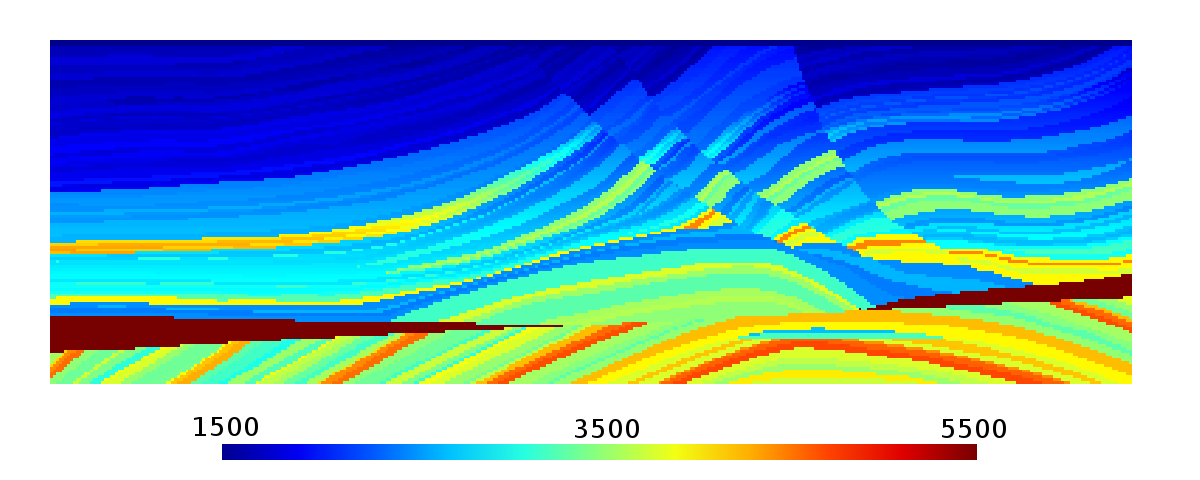} \\
  \caption{Marmousi benchmark. Velocity profile with values from 1500\:m/s to 5500\:m/s.}
  \label{fig625}
\end{figure}

\begin{figure}[!tb]
  \centering
  \small
  \captionsetup{font=small, labelfont=bf}
  \includegraphics[scale = 0.25]{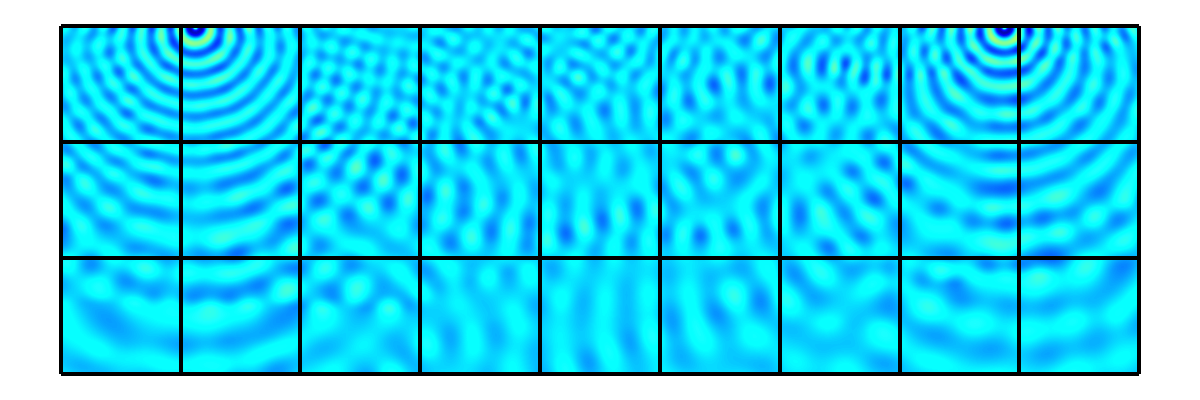} \\
  \caption{Marmousi benchmark. Numerical solution and domain partition with $3\times9$ subdomains.}
  \label{fig626}
\end{figure}

The Marmousi model is a 2D velocity model which is based on the geological structure of the Cuanza basin.
This model exhibits a complex velocity profile $c(\mathbf{x})$ with realistic features (see Figure \ref{fig625}).
It is frequently used to evaluate the performance of numerical solvers with heterogeneous media (\textit{e.g.}~\cite{stolk2013rapidly}).
The numerical simulations are performed over the computational domain\footnote{We assume that all the spatial dimensions are provided in the metric unit [m].} $[0,9192] \times [0,-2904]$.
The Helmholtz equation is solved over the domain, with the HABC prescribed on the boundary of the domain, and two point sources placed at coordinates $(9192/8, -10)$ and $(9192\times 7/8, -10)$, respectively.
The point sources can be placed on interfaces and cross-points (see \textit{e.g.}~\cite{modave2020non}, section 4.5.1., for a discussion of these configurations).
The angular frequency is $\omega = 20 \pi$, and the wavenumber is given by $k(\mathbf{x})=\omega/c(\mathbf{x})$.
Here, the maximum wavenumber is $k_{\text{max}} = 20 \pi / 1500 \approx 0.042$
and the minimum wavenumber is $k_{\text{min}} = 20 \pi / 4500 \approx 0.014$.
In this case, the wavenumber is much smaller than $1$ and the pollution term in \cite{Ihlenburg1997Finite} affects little.
Therefore, the numerical setting is only P1 triangular elements with 20 mesh vertices per wavelength.
The mesh of the computational domain is made of 137808 nodes and 274018 P1 triangles.
The parameters of the HABC operator are $N=4$ and $\phi = \pi /3$ for both the exterior boundary condition and the transmission conditions prescribed at the interfaces between the subdomains.
The numerical solution is shown in Figure \ref{fig626}.

\begin{figure}[!tb]
  \centering
  \small
  \captionsetup{font=small, labelfont=bf}
  \includegraphics[scale = 0.40]{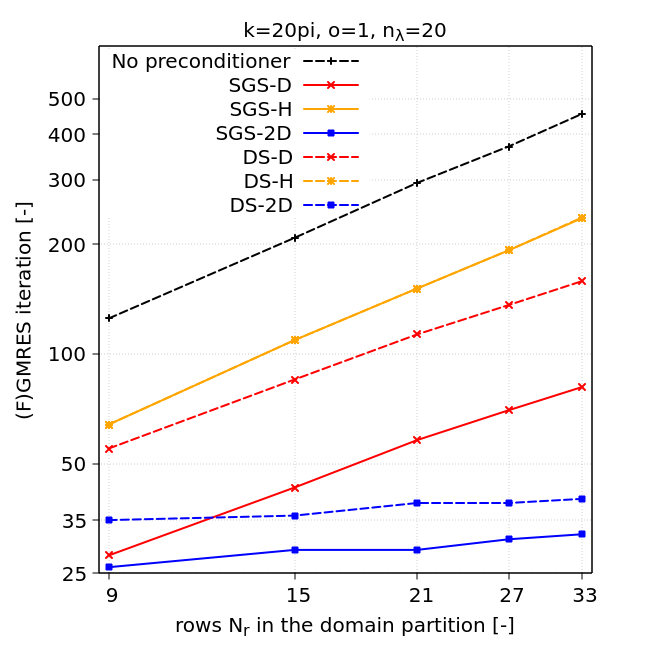} \\
  \caption{Marmousi benchmark. Number of iterations with preconditioned GMRES and F-GMRES (without restart) to reach the relative residual $10^{-6}$ for different domain partitions. The domain is partitioned into $N_r\times N_c$ subdomains, with $N_r=9$, $15$, $21$, $27$, $33$ and $N_c=3 N_r$.}
  \label{fig:nItVSnp_Marmousi}
\end{figure}

The numbers of GMRES iterations to reach a relative residual $10^{-6}$ with the different preconditioners are given in Figure \ref{fig:nItVSnp_Marmousi}.
We have considered domain decompositions into rectangular grids of $N_r \times N_c$ subdomains, where the number of rows is $N_r=3$, $9$, $15$, $21$, $27$, $33$ and the number of columns is $N_c=3 N_r$.
We observe that the number of iterations increases with the number of subdomains with all the preconditionners, but the increase is very slown with F-GMRES and the switching sweeping directions.
Between the coarsest domain partition ($3\times9$ subdomains) and the finest partition ($33\times99$ subdomains), the number of iterations has increased by a factor between $5$ and $9$ with the preconditioners with fixed sweeping directions (SGS-D, SGS-H, DS-D, DS-H), while the factor is smaller than $2$ with the switching sweeping directions (SGS-2D, DS-2D).
In nearly all the cases, the SGS preconditioner with switching sweeping directions requires the smallest number of iterations, but, considering the possible parallelisation of the forward/backward sweeps with the DS preconditioner, that strategy is the best if a parallel environment is used.

\subsection{Acoustic radiation from engine intake}
\label{sec:benchmark:radiation}


\subsubsection*{Description of the benchmark and domain partition}

\begin{figure}[!tb]
  \centering
  \small
  \includegraphics[height=53mm]{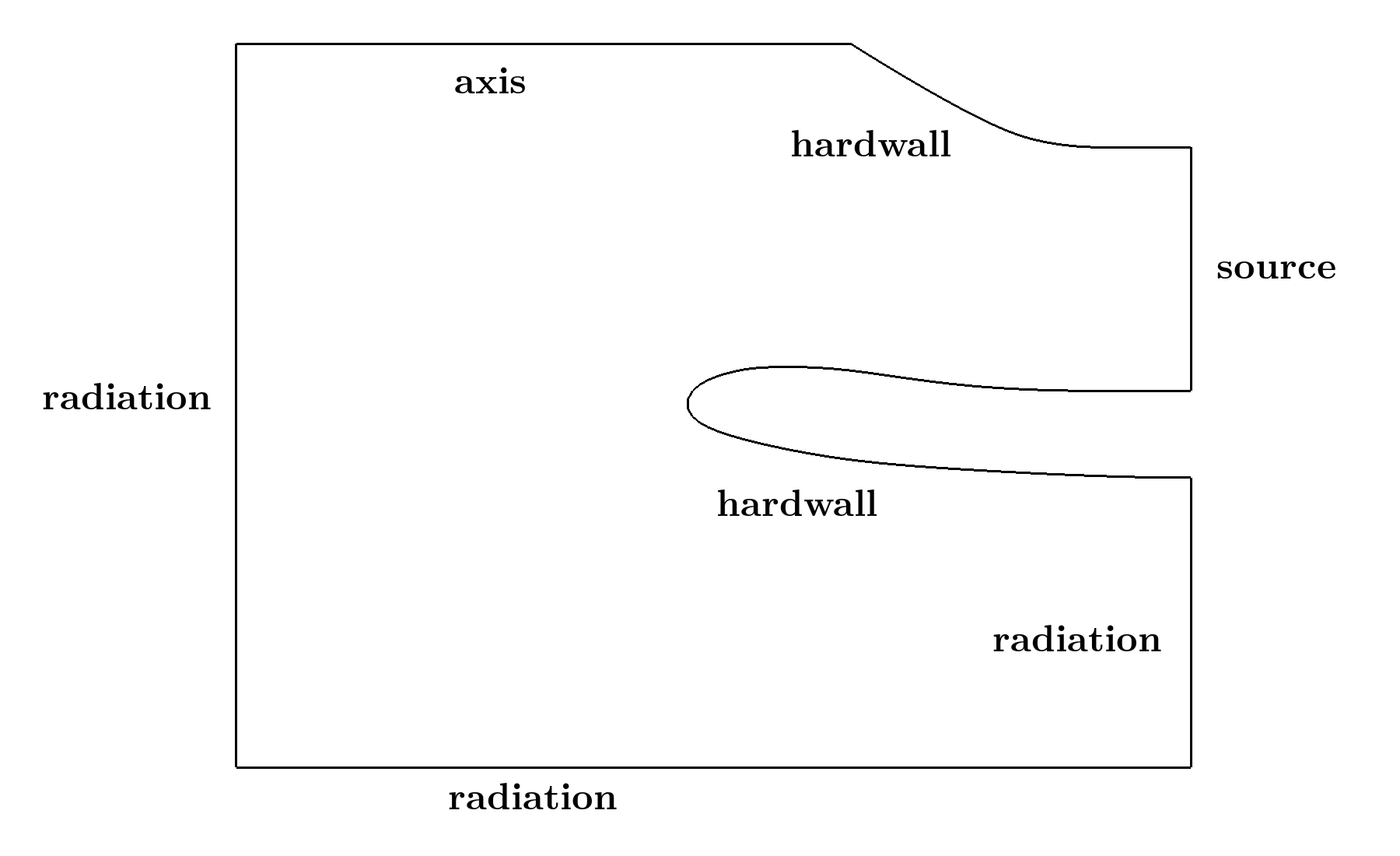} \\
  \captionsetup{font=small, labelfont=bf}
  \caption{Benchmark \textit{``Engine intake''}. Representation of the computational domain and boundary conditions.}
  \label{fig:engine_intake_boundary_conditions}
\end{figure}

\begin{figure}[!tb]
  \centering
  \small
  \captionsetup{font=small, labelfont=bf}
  \includegraphics[height=48mm]{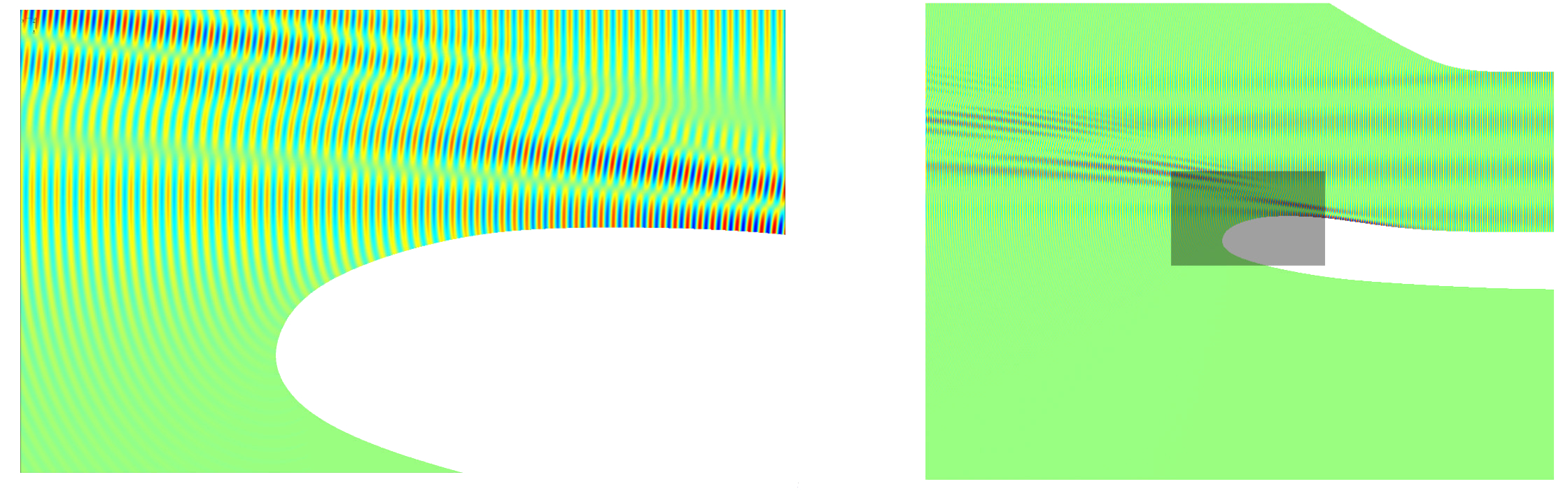} \\
  \caption{Benchmark \textit{``Engine intake''}. Snapshot of the numerical solution for wavenumber $k = 160 \pi$.}
  \label{fig:visu_noiseModel_restart_omega160pi}
\end{figure}

In the last benchmark, we address the computation of a time-harmonic acoustic field in a computational domain that is not rectangular.
It deals with the aeroacoustics of an idealized turbofan engine intake.
The geometry, shown in Figure \ref{fig:engine_intake_boundary_conditions},
is a cylindrical duct of slowly-varying cross-section.
The 2D Helmholtz equation is solved on this computational domain, which is included inside the rectangular region $[-0.3,3.0]\times[0.0,2.5]$.
We consider a Dirichlet condition on the source, $u|_{\text{source}} = \sin((2\pi / 0.50) \cdot y)$, homogeneous Neumann boundary condition on the hardwalls, and an HABC on the artificial borders where the waves must be radiated (see Figure \ref{fig:engine_intake_boundary_conditions}).
For the numerical solution, P5 finite elements are used with 5 elements per wavelength.
The mesh of the computational domain is made of about $3.5\times 10^7$ nodes and $2.8\times10^6$ P5 triangles for wavenumber $160\pi$ ($h \approx 1/2000$).
The parameters of the HABC operators are $N=4$ and $\phi=\pi/3$ for both interface and exterior edges.
The numerical solution corresponding to wavenumber $k=160\pi$ is shown on Figure \ref{fig:visu_noiseModel_restart_omega160pi}.

\begin{figure}[!tb]
  \centering
  \small
  \captionsetup{font=small, labelfont=bf}
  \includegraphics[height=35mm]{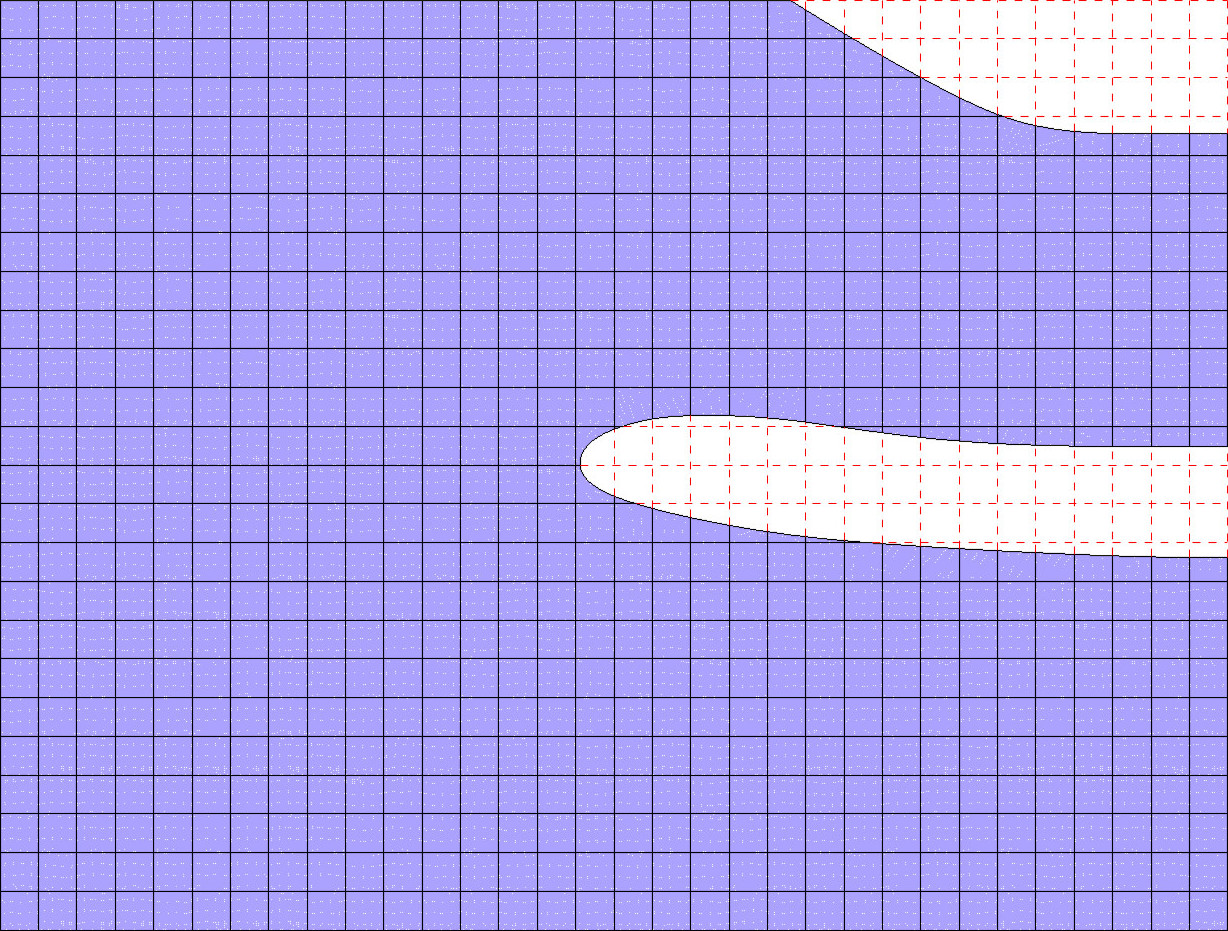} \quad
  \includegraphics[height=35mm]{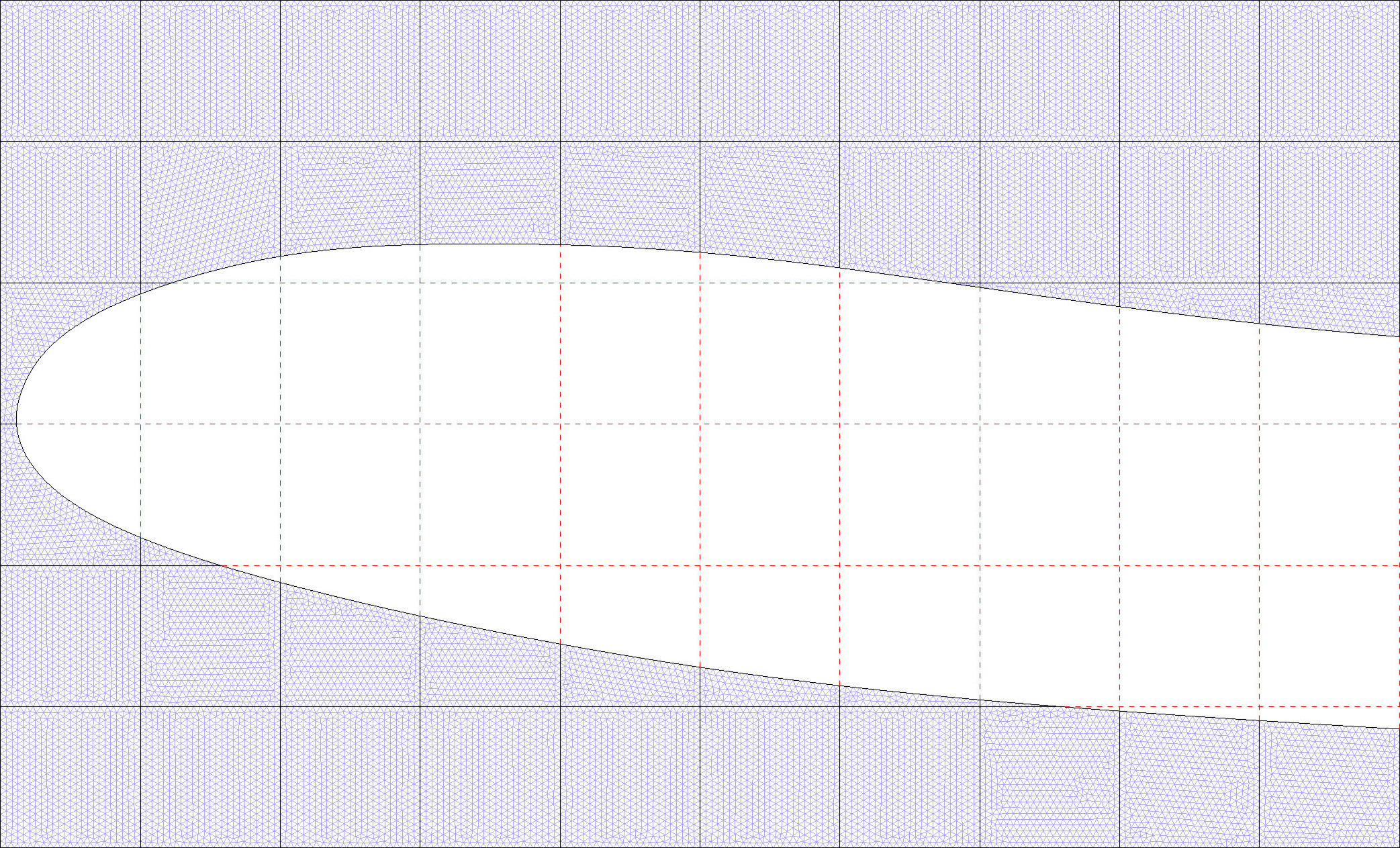} \quad
  \includegraphics[height=35mm]{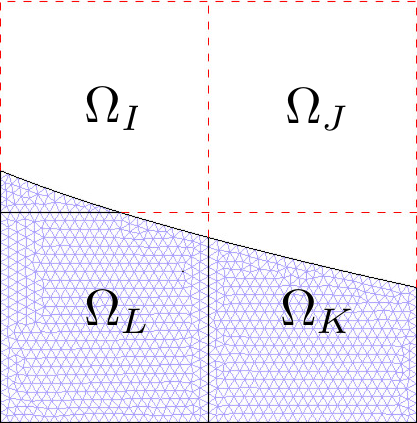}
  \caption{Benchmark \textit{``Engine intake''}. Mesh and example of domain partition.}
  \label{fig:engine_mesh}
\end{figure}

Because the domain is not rectangular, additional steps are required to apply the proposed sweeping preconditioners, which are designed a priori only for checkerboard partitions.
We have generated domain partitions of the rectangular region that contains the computational domain.
This process is performed with Gmsh before the mesh generation.
Then, every partition contains $N_r \times N_c$ rectangular subdomains (see Figure \ref{fig:engine_mesh}), but several ``null'' subdomains are fully outside the computational domain (\textit{e.g.}~subdomain $\Omega_J$ on Figure \ref{fig:engine_mesh}, right), and several subdomains are crossed by the border of the computational domain (\textit{e.g.}~subdomains $\Omega_I$, $\Omega_L$ and $\Omega_K$ on Figure \ref{fig:engine_mesh}, right).
After discretization, no unknown is associated to the null subdomains.
The number of unknowns is smaller for subdomains that are crossed by the domain border than for rectangular subdomains that are fully contained inside the computational domain.
Similarily, the number of discrete transmission variables is smaller on interface edges crossed by the domain boundary.
In practice, the solution procedure is performed by iterating over all the subdomains, wathever they are inside or outside the computational domain.
Dummy systems and dummy vectors of variables are associated to the null subdomains, and dummy transmission variables are associated to interface edges which do not belong to the computational domain (\textit{i.e.}~red dashed edges on Figure \ref{fig:engine_mesh}).
Therefore, the sweeping preconditioner and our computational code can be straightforwardly used for this benchmark.

\subsubsection*{Results}

\begin{figure}[!tb]
  \centering
  \small
  \captionsetup{font=small, labelfont=bf}
  \begin{subfigure}{0.9\textwidth}
  \centering
  \includegraphics[height=40mm]{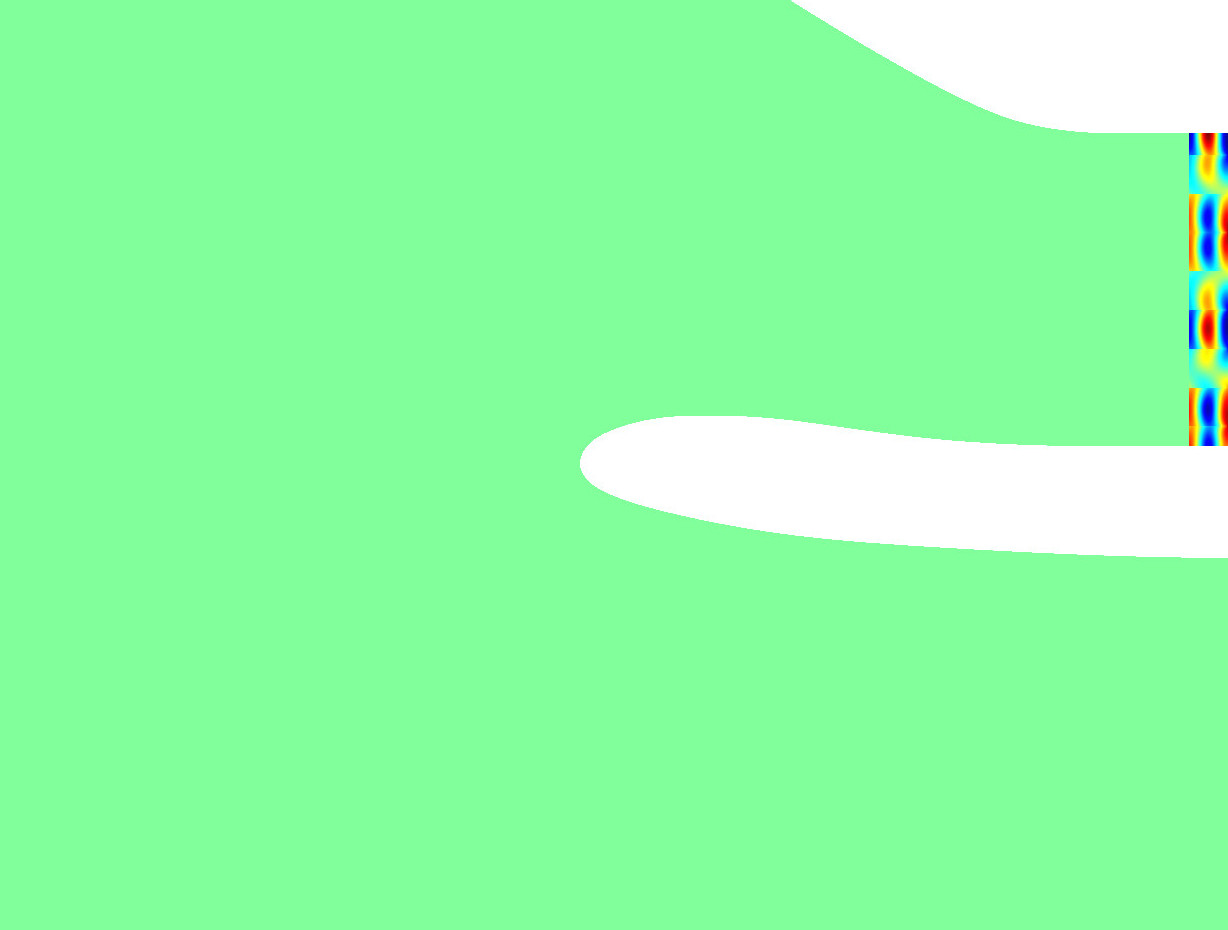}\hspace{1cm}
  \includegraphics[height=40mm]{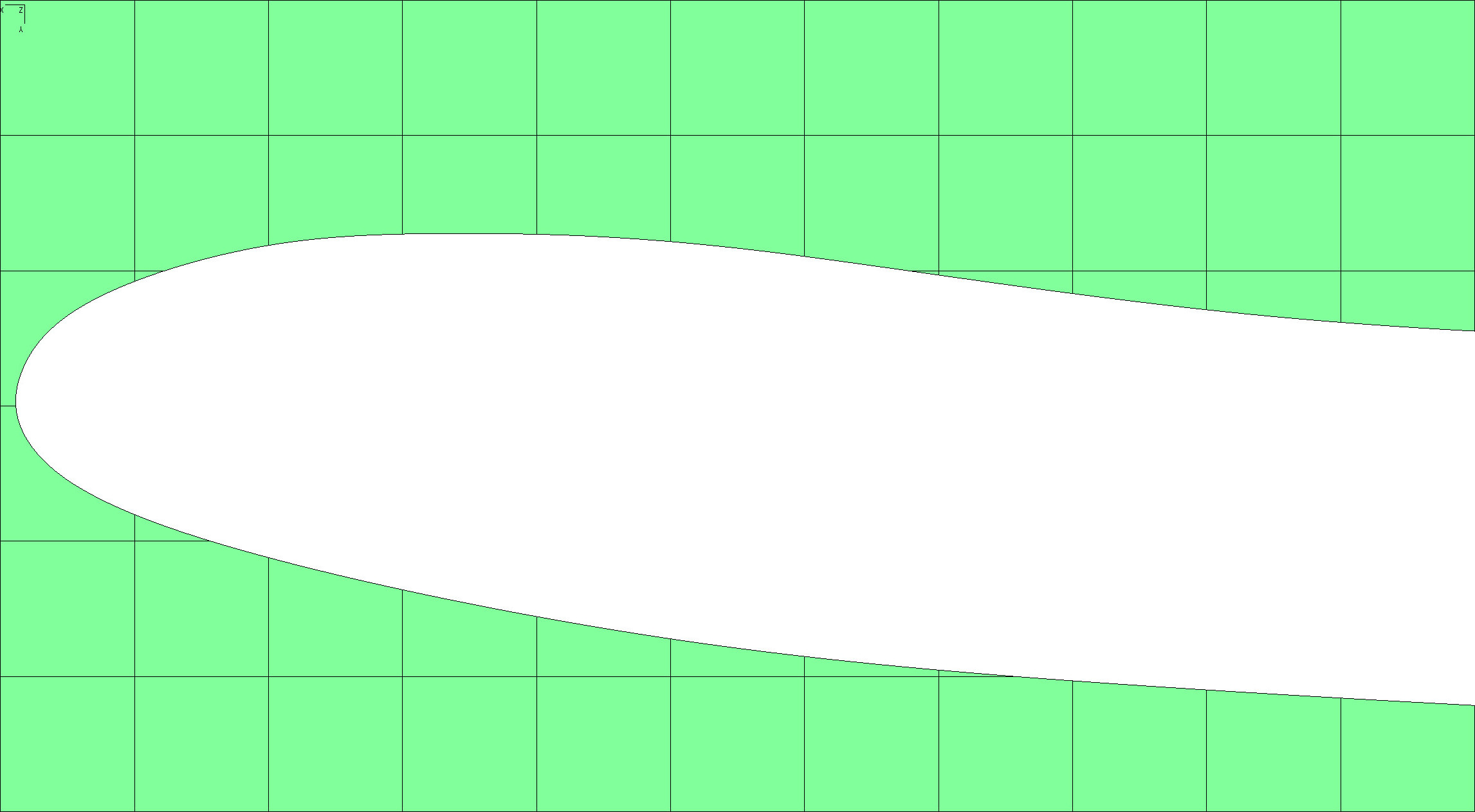}\\
  \caption{iteration 0.}
  \end{subfigure}
  \begin{subfigure}{0.9\textwidth}
  \centering
  \includegraphics[height=40mm]{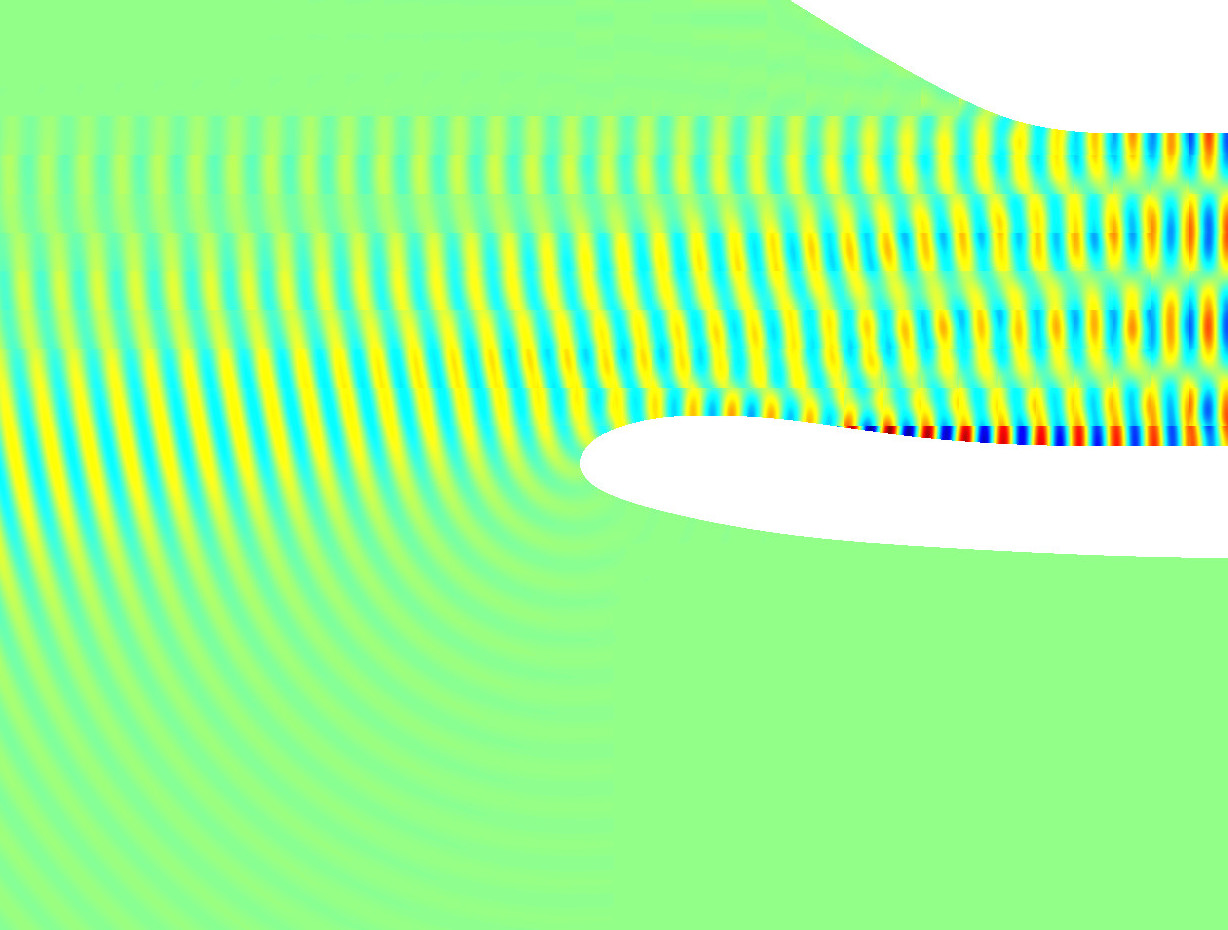}\hspace{1cm}
  \includegraphics[height=40mm]{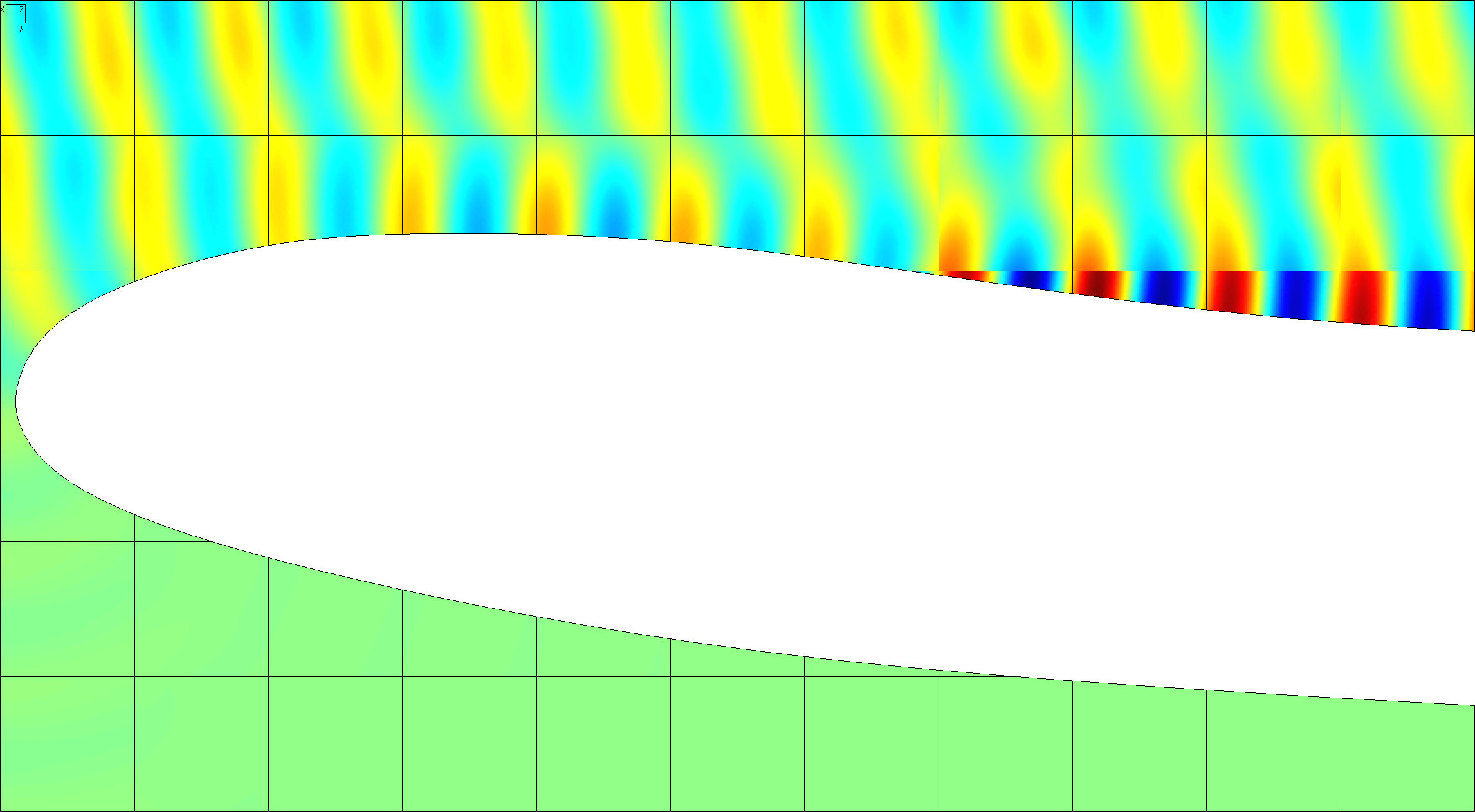}\\
  \caption{iteration 1.}
  \end{subfigure}
  \begin{subfigure}{0.9\textwidth}
  \centering
  \includegraphics[height=40mm]{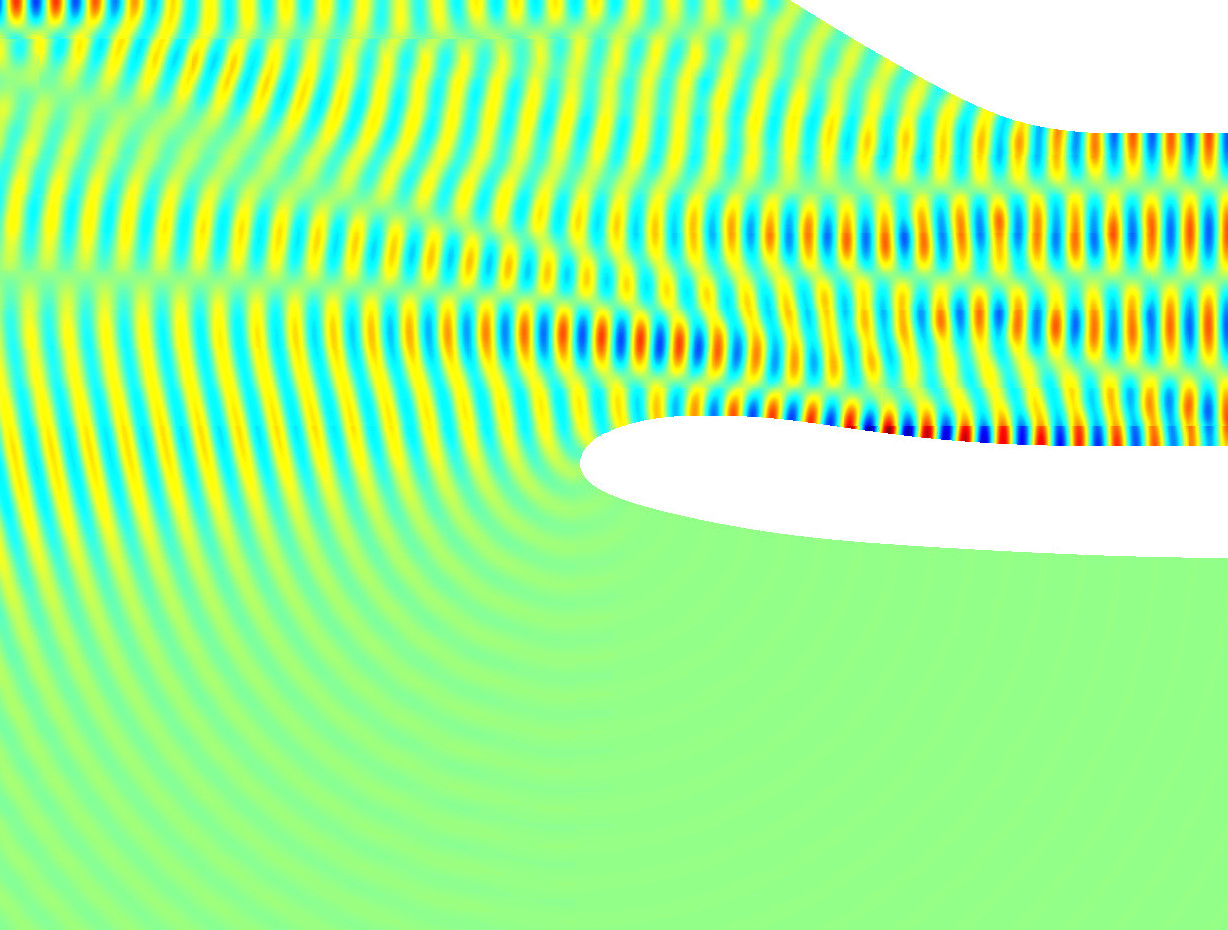}\hspace{1cm}
  \includegraphics[height=40mm]{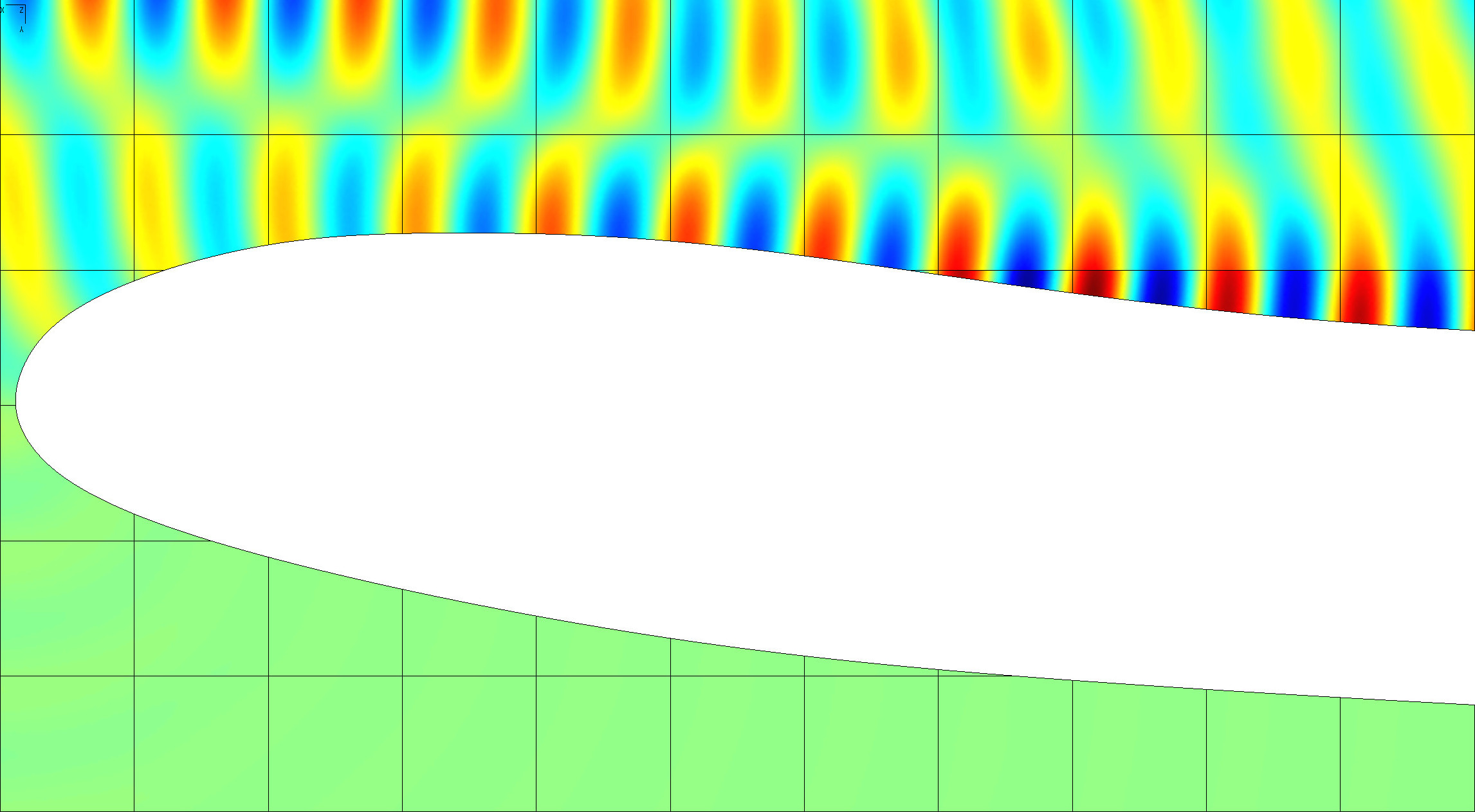}\\
  \caption{iteration 2.}
  \end{subfigure}
  \begin{subfigure}{0.9\textwidth}
  \centering
  \includegraphics[height=40mm]{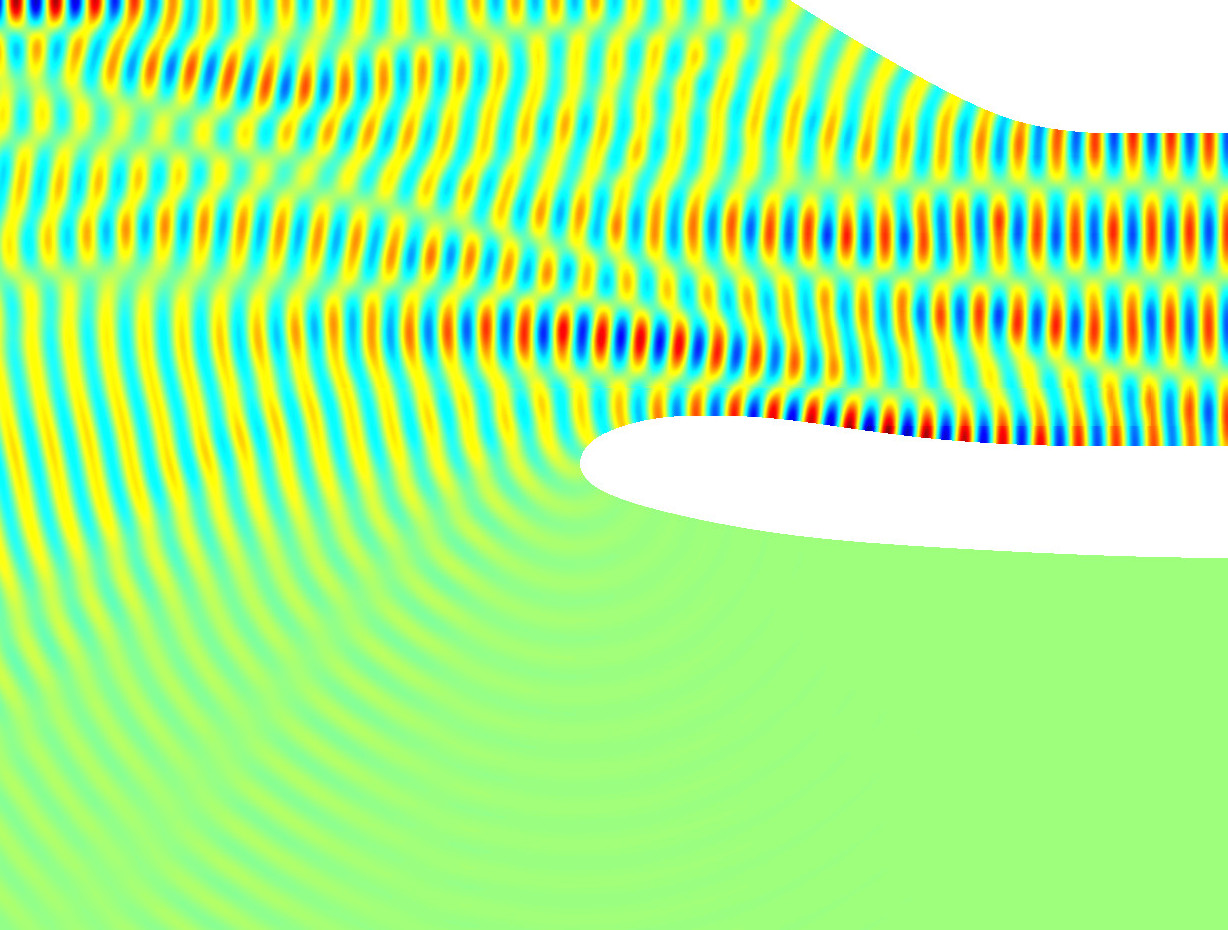}\hspace{1cm}
  \includegraphics[height=40mm]{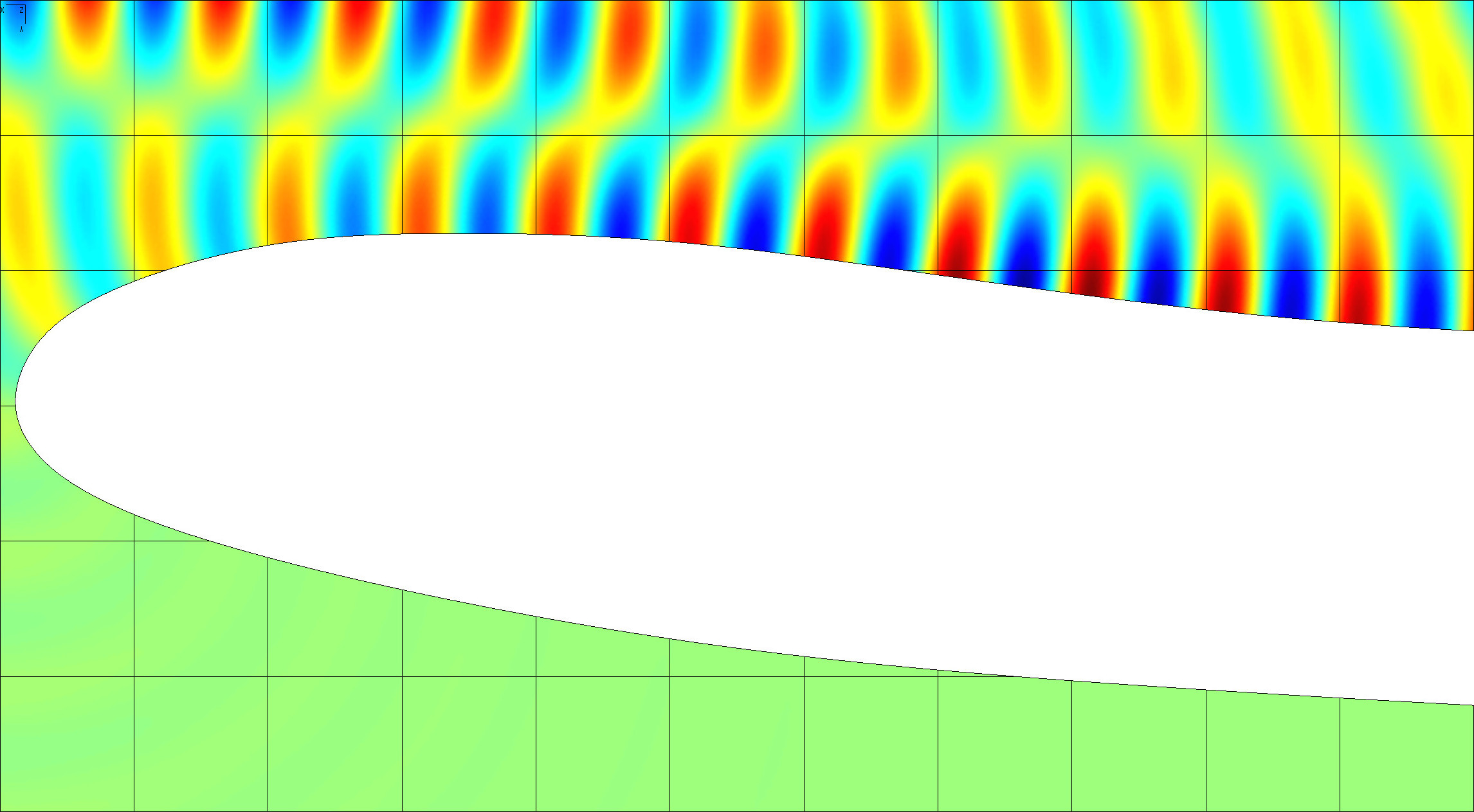}\\
  \caption{iteration 3.}
  \end{subfigure}
  \caption{Benchmark \textit{``Engine intake''}. Snapshot of the solution at the beginning of the procedure (iteration 0) and after 1, 2 and 3 iterations with SGS-2D. The wavenumber is $80\pi$.}
  \label{fig:engine_mesh_and_visu1to3}
\end{figure}

Figure \ref{fig:engine_mesh_and_visu1to3} presents the snapshot of the solution with $48 \times 36$ subdomains
at the beginning of the procedure from iteration 0 to iteration 3 with preconditioners SGS-2D.
The wavenumber $k$ is $80\pi$ in Figure \ref{fig:engine_mesh_and_visu1to3} instead of $160\pi$ considering that the visualisation is more clearly visible.
At the first iteration, partial information is obtained by sweeps that go from the left-down corner to the right-up corner and the other is missed at the top of the computational domain.
At the second iteration, the information at the top of the computational domain is got added, which is
contributed by alternating sweeps.
At the third iteration, more details in the information is presented in the computational domain.

\begin{figure}[!tb]
  \centering
  \small
  \captionsetup{font=small, labelfont=bf}
  \begin{subfigure}[b]{0.47\textwidth}
    \caption{Number of iterations}
    \includegraphics[scale = 0.40]{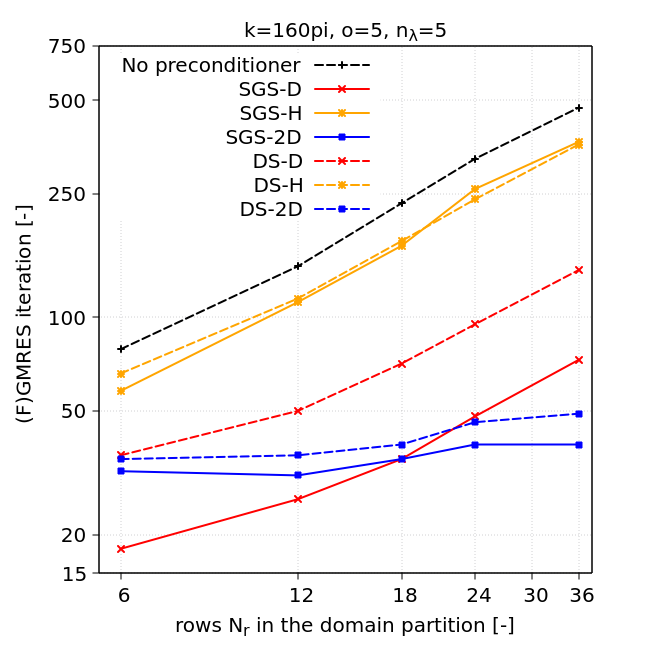}
  \end{subfigure} \quad
  \begin{subfigure}[b]{0.47\textwidth}
    \caption{Runtime in seconds}
    \includegraphics[scale = 0.40]{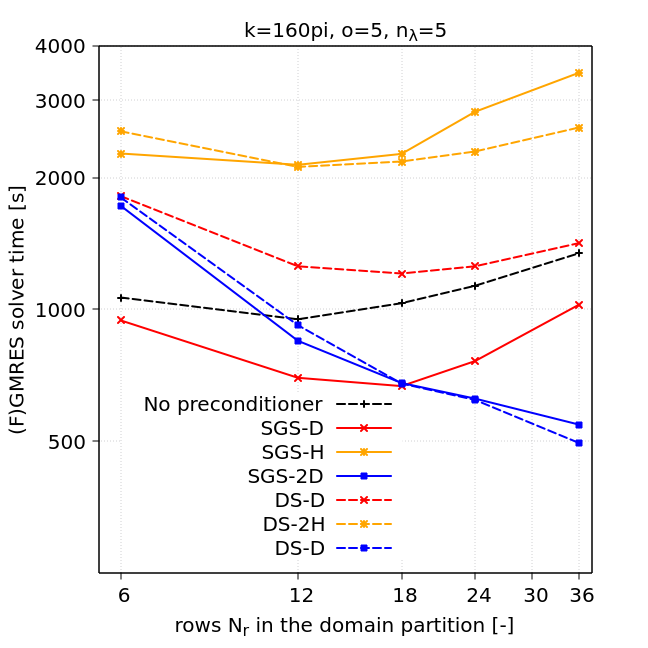}
  \end{subfigure}
  \caption{Benchmark \textit{``Engine intake''}. Number of iterations (a) and runtime in seconds (b) with preconditioned GMRES and F-GMRES (with restart) to reach the relative residual $10^{-6}$ for different domain partitions. The domain is partitioned into $N_r\times N_c$ subdomains, with $N_r=6$, $12$, $18$, $24$, $36$ and $N_c=4 N_r/3$.}
  \label{fig:noiseModel_restart_omega160pi}
\end{figure}

The number of iterations and the runtime to reach a relative residual $10^{-6}$ with the different preconditioners are given in Figure \ref{fig:noiseModel_restart_omega160pi}.
Simulations are carried out on a Intel Xeon Phi (CPU 7210@1.30GHz) and parallelized using OpenMP interface.
The runtime corresponds to the restarted (F)GMRES resolution phase (number of restart = 20). The number of threads is equal to the number of rows of subdomains.

The convergence rate with SGS-2D is the fastest.
The results are comparable with DS-2D. With both preconditioners, the number of iterations is stable.

Comparing flexible preconditioners and fixed preconditioners, we can also see that switching preconditioners improve robustness of DDMs. SGS-H and DS-H are not good enough. We can see that the number of iterations with both preconditioners increases with the number of rows of subdomains. This indicates that SGS-H and DS-H might be suitable for long geometries.

Comparing all the solvers, the DDMs with switching preconditioners perform best. When the number of subdomains increases, the runtime is smaller with switching preconditioners than with the others. The reason is that the number of iterations is smaller with switching preconditioners. There isn’t significant improvement on timing with fixed preconditioners SGS-D or DS-D. Although both preconditioners reduce the number of iterations, the inner steps of these preconditioners at each iteration are time-consuming.

\section{Conclusion}

We have proposed and compared multidirectional sweeping preconditioners for the finite element solution of Helmholtz problems with checkerboard domain decompositions.
The domain decomposition algorithm relies on high-order transmission conditions and a cross-point treatment proposed in \cite{modave2020non}.
This algorithm is well suited to Helmholtz problems with checkerboard partitions of the computational domain, but it cannot scale with the number of subdomains without an efficient preconditioning technique.

While most of the sweeping techniques have been studied for layered domain partitions, we have presented generalizations for checkerboard partitions, offering flexibility in the choice of the sweeping directions.
Horizontal, vertical and diagonal sweeping directions can be used with symmetric Gauss-Seidel and parallel double sweep preconditioners.
Several directions can be combined by using the flexible version of GMRES.
For applicative cases, these preconditioners provide an efficient way to rapidly transfer information in the different zones of the computational domain, then accelerating the convergence of iterative solution procedures with GMRES.
We have observed that the diagonal sweeping directions, with flipping between each iteration of the flexible GMRES, where particularly efficient in all the cases.

The multidirectional sweeping preconditioners can be straightforwardly extended to three dimensions.
Thanks to the block representation of the interface system (equation \eqref{eqn:bigTridiagSys}), every kind of sweeping preconditioner can be applied, such as the symmetric successive over-relaxation (SOR) preconditioner.
For instance, sweeping directions following the diagonals of a cuboid have been tested in \cite{dai2021generalized}.
A parallel quadruple sweep preconditioner, with four sweeps performing simultaneously in the $+x$, $-x$, $+y$ and $-y$ directions, can also be designed as a generalization of the parallel DS preconditioner.
The proposed preconditioners can also be used with other kinds of transmission conditions (\textit{e.g}~with low order conditions, conditions based on perfectly matched layers \cite{royer:hal-03416187}, etc.), since the approach only relies on the structure of the algebraic system.
They can be applied to other kinds of wave propagation problems as well, such as electromagnetic and elastic problems.
Note that, for these problems, cross-point treatments are not available for
high-order transmission conditions, but other kinds of conditions can be used instead.
In future works, these preconditioners will be tested to problems with multiple right-hand sides, and in distributed-memory parallel environments, where novel questions are raised for parallel efficiency.
Multi-dimensional sweeping strategies for automatic domain partitions will also be investigated.

\section{Acknowledgements}

This work was funded in part by the Communaut\'e Fran\c{c}aise de Belgique under
contract ARC WAVES 15/19-03 (``Large Scale Simulation of Waves in Complex
Media'') and by the F.R.S.-FNRS under grant PDR 26104939 (``Fast Helmholtz
Solvers on GPUs'').

\small
\setlength{\bibsep}{2pt plus 0ex}
\bibliographystyle{abbrv}
\bibliography{./sections/myrefs}

\end{document}